\documentclass[12 pt, draftcls, onecolumn]{IEEEtran}	
\usepackage{amsmath, amssymb, amsthm} 					
\usepackage{bm}																	
\usepackage{accents}
\usepackage[noadjust]{cite}
\usepackage[usenames,dvipsnames]{color}
\usepackage{setspace}	
\usepackage{graphicx}													

\usepackage{url}																

\usepackage[plainpages=false]{hyperref}					

\hypersetup{
 colorlinks=true,
 citecolor=Maroon,
 linkcolor=Maroon,
 urlcolor=Maroon}

\usepackage{pst-all}                            

\usepackage{nicefrac}														

\newcommand{\ubar}[1]{\underaccent{\bar}{#1}}

\newcommand{\Prob}{\mbox{Prob}}
\DeclareMathOperator*{\argmax}{\arg\!\max}
\DeclareMathOperator*{\argmin}{\arg\!\min}

\newcommand{\E}{\mathbb{E}}

\newcommand{\Var}{\ensuremath{\textrm{Var}}}

\newtheorem{Theorem}{Theorem}
\newtheorem{Corollary}{Corollary}

\newtheorem{Lemma}{Lemma}
\newtheorem*{Proof*}{Proof}

\newtheorem{Remark}{Remark}

\newtheorem{Property}{Property}

\definecolor{snow}{rgb}{1.000,0.980,0.980}
\definecolor{snow1}{rgb}{1.000,0.980,0.980}
\definecolor{snow2}{rgb}{0.933,0.914,0.914}
\definecolor{snow3}{rgb}{0.804,0.788,0.788}
\definecolor{snow4}{rgb}{0.545,0.537,0.537}
\definecolor{GhostWhite}{rgb}{0.973,0.973,1.000}
\definecolor{WhiteSmoke}{rgb}{0.961,0.961,0.961}
\definecolor{gainsboro}{rgb}{0.863,0.863,0.863}
\definecolor{FloralWhite}{rgb}{1.000,0.980,0.941}
\definecolor{OldLace}{rgb}{0.992,0.961,0.902}
\definecolor{linen}{rgb}{0.980,0.941,0.902}
\definecolor{AntiqueWhite}{rgb}{0.980,0.922,0.843}
\definecolor{PapayaWhip}{rgb}{1.000,0.937,0.835}
\definecolor{BlanchedAlmond}{rgb}{1.000,0.922,0.804}
\definecolor{bisque}{rgb}{1.000,0.894,0.769}
\definecolor{PeachPuff}{rgb}{1.000,0.855,0.725}
\definecolor{NavajoWhite}{rgb}{1.000,0.871,0.678}
\definecolor{moccasin}{rgb}{1.000,0.894,0.710}
\definecolor{cornsilk}{rgb}{1.000,0.973,0.863}
\definecolor{ivory}{rgb}{1.000,1.000,0.941}
\definecolor{LemonChiffon}{rgb}{1.000,0.980,0.804}
\definecolor{seashell}{rgb}{1.000,0.961,0.933}
\definecolor{honeydew}{rgb}{0.941,1.000,0.941}
\definecolor{MintCream}{rgb}{0.961,1.000,0.980}
\definecolor{azure}{rgb}{0.941,1.000,1.000}
\definecolor{AliceBlue}{rgb}{0.941,0.973,1.000}
\definecolor{lavender}{rgb}{0.902,0.902,0.980}
\definecolor{LavenderBlush}{rgb}{1.000,0.941,0.961}
\definecolor{MistyRose}{rgb}{1.000,0.894,0.882}
\definecolor{white}{rgb}{1.000,1.000,1.000}
\definecolor{black}{rgb}{0.000,0.000,0.000}
\definecolor{DarkSlateGray}{rgb}{0.184,0.310,0.310}
\definecolor{DimGray}{rgb}{0.412,0.412,0.412}
\definecolor{SlateGray}{rgb}{0.439,0.502,0.565}
\definecolor{LightSlateGray}{rgb}{0.467,0.533,0.600}
\definecolor{gray}{rgb}{0.745,0.745,0.745}
\definecolor{LightGray}{rgb}{0.827,0.827,0.827}
\definecolor{MidnightBlue}{rgb}{0.098,0.098,0.439}
\definecolor{navy}{rgb}{0.000,0.000,0.502}
\definecolor{NavyBlue}{rgb}{0.000,0.000,0.502}
\definecolor{CornflowerBlue}{rgb}{0.392,0.584,0.929}
\definecolor{DarkSlateBlue}{rgb}{0.282,0.239,0.545}
\definecolor{SlateBlue}{rgb}{0.416,0.353,0.804}
\definecolor{MediumSlateBlue}{rgb}{0.482,0.408,0.933}
\definecolor{LightSlateBlue}{rgb}{0.518,0.439,1.000}
\definecolor{MediumBlue}{rgb}{0.000,0.000,0.804}
\definecolor{RoyalBlue}{rgb}{0.255,0.412,0.882}
\definecolor{blue}{rgb}{0.000,0.000,1.000}
\definecolor{DodgerBlue}{rgb}{0.118,0.565,1.000}
\definecolor{DeepSkyBlue}{rgb}{0.000,0.749,1.000}
\definecolor{SkyBlue}{rgb}{0.529,0.808,0.922}
\definecolor{LightSkyBlue}{rgb}{0.529,0.808,0.980}
\definecolor{SteelBlue}{rgb}{0.275,0.510,0.706}
\definecolor{LightSteelBlue}{rgb}{0.690,0.769,0.871}
\definecolor{LightBlue}{rgb}{0.678,0.847,0.902}
\definecolor{PowderBlue}{rgb}{0.690,0.878,0.902}
\definecolor{PaleTurquoise}{rgb}{0.686,0.933,0.933}
\definecolor{DarkTurquoise}{rgb}{0.000,0.808,0.820}
\definecolor{MediumTurquoise}{rgb}{0.282,0.820,0.800}
\definecolor{turquoise}{rgb}{0.251,0.878,0.816}
\definecolor{cyan}{rgb}{0.000,1.000,1.000}
\definecolor{LightCyan}{rgb}{0.878,1.000,1.000}
\definecolor{CadetBlue}{rgb}{0.373,0.620,0.627}
\definecolor{MediumAquamarine}{rgb}{0.400,0.804,0.667}
\definecolor{aquamarine}{rgb}{0.498,1.000,0.831}
\definecolor{DarkGreen}{rgb}{0.000,0.392,0.000}
\definecolor{DarkOliveGreen}{rgb}{0.333,0.420,0.184}
\definecolor{DarkSeaGreen}{rgb}{0.561,0.737,0.561}
\definecolor{SeaGreen}{rgb}{0.180,0.545,0.341}
\definecolor{MediumSeaGreen}{rgb}{0.235,0.702,0.443}
\definecolor{LightSeaGreen}{rgb}{0.125,0.698,0.667}
\definecolor{PaleGreen}{rgb}{0.596,0.984,0.596}
\definecolor{SpringGreen}{rgb}{0.000,1.000,0.498}
\definecolor{LawnGreen}{rgb}{0.486,0.988,0.000}
\definecolor{green}{rgb}{0.000,1.000,0.000}
\definecolor{chartreuse}{rgb}{0.498,1.000,0.000}
\definecolor{MediumSpringGreen}{rgb}{0.000,0.980,0.604}
\definecolor{GreenYellow}{rgb}{0.678,1.000,0.184}
\definecolor{LimeGreen}{rgb}{0.196,0.804,0.196}
\definecolor{YellowGreen}{rgb}{0.604,0.804,0.196}
\definecolor{ForestGreen}{rgb}{0.133,0.545,0.133}
\definecolor{OliveDrab}{rgb}{0.420,0.557,0.137}
\definecolor{DarkKhaki}{rgb}{0.741,0.718,0.420}
\definecolor{khaki}{rgb}{0.941,0.902,0.549}
\definecolor{PaleGoldenrod}{rgb}{0.933,0.910,0.667}
\definecolor{LightGoldenrodYellow}{rgb}{0.980,0.980,0.824}
\definecolor{LightYellow}{rgb}{1.000,1.000,0.878}
\definecolor{yellow}{rgb}{1.000,1.000,0.000}
\definecolor{gold}{rgb}{1.000,0.843,0.000}
\definecolor{LightGoldenrod}{rgb}{0.933,0.867,0.510}
\definecolor{goldenrod}{rgb}{0.855,0.647,0.125}
\definecolor{DarkGoldenrod}{rgb}{0.722,0.525,0.043}
\definecolor{RosyBrown}{rgb}{0.737,0.561,0.561}
\definecolor{IndianRed}{rgb}{0.804,0.361,0.361}
\definecolor{SaddleBrown}{rgb}{0.545,0.271,0.075}
\definecolor{sienna}{rgb}{0.627,0.322,0.176}
\definecolor{peru}{rgb}{0.804,0.522,0.247}
\definecolor{burlywood}{rgb}{0.871,0.722,0.529}
\definecolor{beige}{rgb}{0.961,0.961,0.863}
\definecolor{wheat}{rgb}{0.961,0.871,0.702}
\definecolor{SandyBrown}{rgb}{0.957,0.643,0.376}
\definecolor{tan}{rgb}{0.824,0.706,0.549}
\definecolor{chocolate}{rgb}{0.824,0.412,0.118}
\definecolor{firebrick}{rgb}{0.698,0.133,0.133}
\definecolor{brown}{rgb}{0.647,0.165,0.165}
\definecolor{DarkSalmon}{rgb}{0.914,0.588,0.478}
\definecolor{salmon}{rgb}{0.980,0.502,0.447}
\definecolor{LightSalmon}{rgb}{1.000,0.627,0.478}
\definecolor{orange}{rgb}{1.000,0.647,0.000}
\definecolor{DarkOrange}{rgb}{1.000,0.549,0.000}
\definecolor{coral}{rgb}{1.000,0.498,0.314}
\definecolor{LightCoral}{rgb}{0.941,0.502,0.502}
\definecolor{tomato}{rgb}{1.000,0.388,0.278}
\definecolor{OrangeRed}{rgb}{1.000,0.271,0.000}
\definecolor{red}{rgb}{1.000,0.000,0.000}
\definecolor{HotPink}{rgb}{1.000,0.412,0.706}
\definecolor{DeepPink}{rgb}{1.000,0.078,0.576}
\definecolor{pink}{rgb}{1.000,0.753,0.796}
\definecolor{LightPink}{rgb}{1.000,0.714,0.757}
\definecolor{PaleVioletRed}{rgb}{0.859,0.439,0.576}
\definecolor{maroon}{rgb}{0.690,0.188,0.376}
\definecolor{MediumVioletRed}{rgb}{0.780,0.082,0.522}
\definecolor{VioletRed}{rgb}{0.816,0.125,0.565}
\definecolor{magenta}{rgb}{1.000,0.000,1.000}
\definecolor{violet}{rgb}{0.933,0.510,0.933}
\definecolor{plum}{rgb}{0.867,0.627,0.867}
\definecolor{orchid}{rgb}{0.855,0.439,0.839}
\definecolor{MediumOrchid}{rgb}{0.729,0.333,0.827}
\definecolor{DarkOrchid}{rgb}{0.600,0.196,0.800}
\definecolor{DarkViolet}{rgb}{0.580,0.000,0.827}
\definecolor{BlueViolet}{rgb}{0.541,0.169,0.886}
\definecolor{purple}{rgb}{0.627,0.125,0.941}
\definecolor{MediumPurple}{rgb}{0.576,0.439,0.859}
\definecolor{thistle}{rgb}{0.847,0.749,0.847}
\definecolor{seashell1}{rgb}{1.000,0.961,0.933}
\definecolor{seashell2}{rgb}{0.933,0.898,0.871}
\definecolor{seashell3}{rgb}{0.804,0.773,0.749}
\definecolor{seashell4}{rgb}{0.545,0.525,0.510}
\definecolor{AntiqueWhite1}{rgb}{1.000,0.937,0.859}
\definecolor{AntiqueWhite2}{rgb}{0.933,0.875,0.800}
\definecolor{AntiqueWhite3}{rgb}{0.804,0.753,0.690}
\definecolor{AntiqueWhite4}{rgb}{0.545,0.514,0.471}
\definecolor{bisque1}{rgb}{1.000,0.894,0.769}
\definecolor{bisque2}{rgb}{0.933,0.835,0.718}
\definecolor{bisque3}{rgb}{0.804,0.718,0.620}
\definecolor{bisque4}{rgb}{0.545,0.490,0.420}
\definecolor{PeachPuff1}{rgb}{1.000,0.855,0.725}
\definecolor{PeachPuff2}{rgb}{0.933,0.796,0.678}
\definecolor{PeachPuff3}{rgb}{0.804,0.686,0.584}
\definecolor{PeachPuff4}{rgb}{0.545,0.467,0.396}
\definecolor{NavajoWhite1}{rgb}{1.000,0.871,0.678}
\definecolor{NavajoWhite2}{rgb}{0.933,0.812,0.631}
\definecolor{NavajoWhite3}{rgb}{0.804,0.702,0.545}
\definecolor{NavajoWhite4}{rgb}{0.545,0.475,0.369}
\definecolor{LemonChiffon1}{rgb}{1.000,0.980,0.804}
\definecolor{LemonChiffon2}{rgb}{0.933,0.914,0.749}
\definecolor{LemonChiffon3}{rgb}{0.804,0.788,0.647}
\definecolor{LemonChiffon4}{rgb}{0.545,0.537,0.439}
\definecolor{cornsilk1}{rgb}{1.000,0.973,0.863}
\definecolor{cornsilk2}{rgb}{0.933,0.910,0.804}
\definecolor{cornsilk3}{rgb}{0.804,0.784,0.694}
\definecolor{cornsilk4}{rgb}{0.545,0.533,0.471}
\definecolor{ivory1}{rgb}{1.000,1.000,0.941}
\definecolor{ivory2}{rgb}{0.933,0.933,0.878}
\definecolor{ivory3}{rgb}{0.804,0.804,0.757}
\definecolor{ivory4}{rgb}{0.545,0.545,0.514}
\definecolor{honeydew1}{rgb}{0.941,1.000,0.941}
\definecolor{honeydew2}{rgb}{0.878,0.933,0.878}
\definecolor{honeydew3}{rgb}{0.757,0.804,0.757}
\definecolor{honeydew4}{rgb}{0.514,0.545,0.514}
\definecolor{LavenderBlush1}{rgb}{1.000,0.941,0.961}
\definecolor{LavenderBlush2}{rgb}{0.933,0.878,0.898}
\definecolor{LavenderBlush3}{rgb}{0.804,0.757,0.773}
\definecolor{LavenderBlush4}{rgb}{0.545,0.514,0.525}
\definecolor{MistyRose1}{rgb}{1.000,0.894,0.882}
\definecolor{MistyRose2}{rgb}{0.933,0.835,0.824}
\definecolor{MistyRose3}{rgb}{0.804,0.718,0.710}
\definecolor{MistyRose4}{rgb}{0.545,0.490,0.482}
\definecolor{azure1}{rgb}{0.941,1.000,1.000}
\definecolor{azure2}{rgb}{0.878,0.933,0.933}
\definecolor{azure3}{rgb}{0.757,0.804,0.804}
\definecolor{azure4}{rgb}{0.514,0.545,0.545}
\definecolor{SlateBlue1}{rgb}{0.514,0.435,1.000}
\definecolor{SlateBlue2}{rgb}{0.478,0.404,0.933}
\definecolor{SlateBlue3}{rgb}{0.412,0.349,0.804}
\definecolor{SlateBlue4}{rgb}{0.278,0.235,0.545}
\definecolor{RoyalBlue1}{rgb}{0.282,0.463,1.000}
\definecolor{RoyalBlue2}{rgb}{0.263,0.431,0.933}
\definecolor{RoyalBlue3}{rgb}{0.227,0.373,0.804}
\definecolor{RoyalBlue4}{rgb}{0.153,0.251,0.545}
\definecolor{blue1}{rgb}{0.000,0.000,1.000}
\definecolor{blue2}{rgb}{0.000,0.000,0.933}
\definecolor{blue3}{rgb}{0.000,0.000,0.804}
\definecolor{blue4}{rgb}{0.000,0.000,0.545}
\definecolor{DodgerBlue1}{rgb}{0.118,0.565,1.000}
\definecolor{DodgerBlue2}{rgb}{0.110,0.525,0.933}
\definecolor{DodgerBlue3}{rgb}{0.094,0.455,0.804}
\definecolor{DodgerBlue4}{rgb}{0.063,0.306,0.545}
\definecolor{SteelBlue1}{rgb}{0.388,0.722,1.000}
\definecolor{SteelBlue2}{rgb}{0.361,0.675,0.933}
\definecolor{SteelBlue3}{rgb}{0.310,0.580,0.804}
\definecolor{SteelBlue4}{rgb}{0.212,0.392,0.545}
\definecolor{DeepSkyBlue1}{rgb}{0.000,0.749,1.000}
\definecolor{DeepSkyBlue2}{rgb}{0.000,0.698,0.933}
\definecolor{DeepSkyBlue3}{rgb}{0.000,0.604,0.804}
\definecolor{DeepSkyBlue4}{rgb}{0.000,0.408,0.545}
\definecolor{SkyBlue1}{rgb}{0.529,0.808,1.000}
\definecolor{SkyBlue2}{rgb}{0.494,0.753,0.933}
\definecolor{SkyBlue3}{rgb}{0.424,0.651,0.804}
\definecolor{SkyBlue4}{rgb}{0.290,0.439,0.545}
\definecolor{LightSkyBlue1}{rgb}{0.690,0.886,1.000}
\definecolor{LightSkyBlue2}{rgb}{0.643,0.827,0.933}
\definecolor{LightSkyBlue3}{rgb}{0.553,0.714,0.804}
\definecolor{LightSkyBlue4}{rgb}{0.376,0.482,0.545}
\definecolor{SlateGray1}{rgb}{0.776,0.886,1.000}
\definecolor{SlateGray2}{rgb}{0.725,0.827,0.933}
\definecolor{SlateGray3}{rgb}{0.624,0.714,0.804}
\definecolor{SlateGray4}{rgb}{0.424,0.482,0.545}
\definecolor{LightSteelBlue1}{rgb}{0.792,0.882,1.000}
\definecolor{LightSteelBlue2}{rgb}{0.737,0.824,0.933}
\definecolor{LightSteelBlue3}{rgb}{0.635,0.710,0.804}
\definecolor{LightSteelBlue4}{rgb}{0.431,0.482,0.545}
\definecolor{LightBlue1}{rgb}{0.749,0.937,1.000}
\definecolor{LightBlue2}{rgb}{0.698,0.875,0.933}
\definecolor{LightBlue3}{rgb}{0.604,0.753,0.804}
\definecolor{LightBlue4}{rgb}{0.408,0.514,0.545}
\definecolor{LightCyan1}{rgb}{0.878,1.000,1.000}
\definecolor{LightCyan2}{rgb}{0.820,0.933,0.933}
\definecolor{LightCyan3}{rgb}{0.706,0.804,0.804}
\definecolor{LightCyan4}{rgb}{0.478,0.545,0.545}
\definecolor{PaleTurquoise1}{rgb}{0.733,1.000,1.000}
\definecolor{PaleTurquoise2}{rgb}{0.682,0.933,0.933}
\definecolor{PaleTurquoise3}{rgb}{0.588,0.804,0.804}
\definecolor{PaleTurquoise4}{rgb}{0.400,0.545,0.545}
\definecolor{CadetBlue1}{rgb}{0.596,0.961,1.000}
\definecolor{CadetBlue2}{rgb}{0.557,0.898,0.933}
\definecolor{CadetBlue3}{rgb}{0.478,0.773,0.804}
\definecolor{CadetBlue4}{rgb}{0.325,0.525,0.545}
\definecolor{turquoise1}{rgb}{0.000,0.961,1.000}
\definecolor{turquoise2}{rgb}{0.000,0.898,0.933}
\definecolor{turquoise3}{rgb}{0.000,0.773,0.804}
\definecolor{turquoise4}{rgb}{0.000,0.525,0.545}
\definecolor{cyan1}{rgb}{0.000,1.000,1.000}
\definecolor{cyan2}{rgb}{0.000,0.933,0.933}
\definecolor{cyan3}{rgb}{0.000,0.804,0.804}
\definecolor{cyan4}{rgb}{0.000,0.545,0.545}
\definecolor{DarkSlateGray1}{rgb}{0.592,1.000,1.000}
\definecolor{DarkSlateGray2}{rgb}{0.553,0.933,0.933}
\definecolor{DarkSlateGray3}{rgb}{0.475,0.804,0.804}
\definecolor{DarkSlateGray4}{rgb}{0.322,0.545,0.545}
\definecolor{aquamarine1}{rgb}{0.498,1.000,0.831}
\definecolor{aquamarine2}{rgb}{0.463,0.933,0.776}
\definecolor{aquamarine3}{rgb}{0.400,0.804,0.667}
\definecolor{aquamarine4}{rgb}{0.271,0.545,0.455}
\definecolor{DarkSeaGreen1}{rgb}{0.757,1.000,0.757}
\definecolor{DarkSeaGreen2}{rgb}{0.706,0.933,0.706}
\definecolor{DarkSeaGreen3}{rgb}{0.608,0.804,0.608}
\definecolor{DarkSeaGreen4}{rgb}{0.412,0.545,0.412}
\definecolor{SeaGreen1}{rgb}{0.329,1.000,0.624}
\definecolor{SeaGreen2}{rgb}{0.306,0.933,0.580}
\definecolor{SeaGreen3}{rgb}{0.263,0.804,0.502}
\definecolor{SeaGreen4}{rgb}{0.180,0.545,0.341}
\definecolor{PaleGreen1}{rgb}{0.604,1.000,0.604}
\definecolor{PaleGreen2}{rgb}{0.565,0.933,0.565}
\definecolor{PaleGreen3}{rgb}{0.486,0.804,0.486}
\definecolor{PaleGreen4}{rgb}{0.329,0.545,0.329}
\definecolor{SpringGreen1}{rgb}{0.000,1.000,0.498}
\definecolor{SpringGreen2}{rgb}{0.000,0.933,0.463}
\definecolor{SpringGreen3}{rgb}{0.000,0.804,0.400}
\definecolor{SpringGreen4}{rgb}{0.000,0.545,0.271}
\definecolor{green1}{rgb}{0.000,1.000,0.000}
\definecolor{green2}{rgb}{0.000,0.933,0.000}
\definecolor{green3}{rgb}{0.000,0.804,0.000}
\definecolor{green4}{rgb}{0.000,0.545,0.000}
\definecolor{chartreuse1}{rgb}{0.498,1.000,0.000}
\definecolor{chartreuse2}{rgb}{0.463,0.933,0.000}
\definecolor{chartreuse3}{rgb}{0.400,0.804,0.000}
\definecolor{chartreuse4}{rgb}{0.271,0.545,0.000}
\definecolor{OliveDrab1}{rgb}{0.753,1.000,0.243}
\definecolor{OliveDrab2}{rgb}{0.702,0.933,0.227}
\definecolor{OliveDrab3}{rgb}{0.604,0.804,0.196}
\definecolor{OliveDrab4}{rgb}{0.412,0.545,0.133}
\definecolor{DarkOliveGreen1}{rgb}{0.792,1.000,0.439}
\definecolor{DarkOliveGreen2}{rgb}{0.737,0.933,0.408}
\definecolor{DarkOliveGreen3}{rgb}{0.635,0.804,0.353}
\definecolor{DarkOliveGreen4}{rgb}{0.431,0.545,0.239}
\definecolor{khaki1}{rgb}{1.000,0.965,0.561}
\definecolor{khaki2}{rgb}{0.933,0.902,0.522}
\definecolor{khaki3}{rgb}{0.804,0.776,0.451}
\definecolor{khaki4}{rgb}{0.545,0.525,0.306}
\definecolor{LightGoldenrod1}{rgb}{1.000,0.925,0.545}
\definecolor{LightGoldenrod2}{rgb}{0.933,0.863,0.510}
\definecolor{LightGoldenrod3}{rgb}{0.804,0.745,0.439}
\definecolor{LightGoldenrod4}{rgb}{0.545,0.506,0.298}
\definecolor{LightYellow1}{rgb}{1.000,1.000,0.878}
\definecolor{LightYellow2}{rgb}{0.933,0.933,0.820}
\definecolor{LightYellow3}{rgb}{0.804,0.804,0.706}
\definecolor{LightYellow4}{rgb}{0.545,0.545,0.478}
\definecolor{yellow1}{rgb}{1.000,1.000,0.000}
\definecolor{yellow2}{rgb}{0.933,0.933,0.000}
\definecolor{yellow3}{rgb}{0.804,0.804,0.000}
\definecolor{yellow4}{rgb}{0.545,0.545,0.000}
\definecolor{gold1}{rgb}{1.000,0.843,0.000}
\definecolor{gold2}{rgb}{0.933,0.788,0.000}
\definecolor{gold3}{rgb}{0.804,0.678,0.000}
\definecolor{gold4}{rgb}{0.545,0.459,0.000}
\definecolor{goldenrod1}{rgb}{1.000,0.757,0.145}
\definecolor{goldenrod2}{rgb}{0.933,0.706,0.133}
\definecolor{goldenrod3}{rgb}{0.804,0.608,0.114}
\definecolor{goldenrod4}{rgb}{0.545,0.412,0.078}
\definecolor{DarkGoldenrod1}{rgb}{1.000,0.725,0.059}
\definecolor{DarkGoldenrod2}{rgb}{0.933,0.678,0.055}
\definecolor{DarkGoldenrod3}{rgb}{0.804,0.584,0.047}
\definecolor{DarkGoldenrod4}{rgb}{0.545,0.396,0.031}
\definecolor{RosyBrown1}{rgb}{1.000,0.757,0.757}
\definecolor{RosyBrown2}{rgb}{0.933,0.706,0.706}
\definecolor{RosyBrown3}{rgb}{0.804,0.608,0.608}
\definecolor{RosyBrown4}{rgb}{0.545,0.412,0.412}
\definecolor{IndianRed1}{rgb}{1.000,0.416,0.416}
\definecolor{IndianRed2}{rgb}{0.933,0.388,0.388}
\definecolor{IndianRed3}{rgb}{0.804,0.333,0.333}
\definecolor{IndianRed4}{rgb}{0.545,0.227,0.227}
\definecolor{sienna1}{rgb}{1.000,0.510,0.278}
\definecolor{sienna2}{rgb}{0.933,0.475,0.259}
\definecolor{sienna3}{rgb}{0.804,0.408,0.224}
\definecolor{sienna4}{rgb}{0.545,0.278,0.149}
\definecolor{burlywood1}{rgb}{1.000,0.827,0.608}
\definecolor{burlywood2}{rgb}{0.933,0.773,0.569}
\definecolor{burlywood3}{rgb}{0.804,0.667,0.490}
\definecolor{burlywood4}{rgb}{0.545,0.451,0.333}
\definecolor{wheat1}{rgb}{1.000,0.906,0.729}
\definecolor{wheat2}{rgb}{0.933,0.847,0.682}
\definecolor{wheat3}{rgb}{0.804,0.729,0.588}
\definecolor{wheat4}{rgb}{0.545,0.494,0.400}
\definecolor{tan1}{rgb}{1.000,0.647,0.310}
\definecolor{tan2}{rgb}{0.933,0.604,0.286}
\definecolor{tan3}{rgb}{0.804,0.522,0.247}
\definecolor{tan4}{rgb}{0.545,0.353,0.169}
\definecolor{chocolate1}{rgb}{1.000,0.498,0.141}
\definecolor{chocolate2}{rgb}{0.933,0.463,0.129}
\definecolor{chocolate3}{rgb}{0.804,0.400,0.114}
\definecolor{chocolate4}{rgb}{0.545,0.271,0.075}
\definecolor{firebrick1}{rgb}{1.000,0.188,0.188}
\definecolor{firebrick2}{rgb}{0.933,0.173,0.173}
\definecolor{firebrick3}{rgb}{0.804,0.149,0.149}
\definecolor{firebrick4}{rgb}{0.545,0.102,0.102}
\definecolor{brown1}{rgb}{1.000,0.251,0.251}
\definecolor{brown2}{rgb}{0.933,0.231,0.231}
\definecolor{brown3}{rgb}{0.804,0.200,0.200}
\definecolor{brown4}{rgb}{0.545,0.137,0.137}
\definecolor{salmon1}{rgb}{1.000,0.549,0.412}
\definecolor{salmon2}{rgb}{0.933,0.510,0.384}
\definecolor{salmon3}{rgb}{0.804,0.439,0.329}
\definecolor{salmon4}{rgb}{0.545,0.298,0.224}
\definecolor{LightSalmon1}{rgb}{1.000,0.627,0.478}
\definecolor{LightSalmon2}{rgb}{0.933,0.584,0.447}
\definecolor{LightSalmon3}{rgb}{0.804,0.506,0.384}
\definecolor{LightSalmon4}{rgb}{0.545,0.341,0.259}
\definecolor{orange1}{rgb}{1.000,0.647,0.000}
\definecolor{orange2}{rgb}{0.933,0.604,0.000}
\definecolor{orange3}{rgb}{0.804,0.522,0.000}
\definecolor{orange4}{rgb}{0.545,0.353,0.000}
\definecolor{DarkOrange1}{rgb}{1.000,0.498,0.000}
\definecolor{DarkOrange2}{rgb}{0.933,0.463,0.000}
\definecolor{DarkOrange3}{rgb}{0.804,0.400,0.000}
\definecolor{DarkOrange4}{rgb}{0.545,0.271,0.000}
\definecolor{coral1}{rgb}{1.000,0.447,0.337}
\definecolor{coral2}{rgb}{0.933,0.416,0.314}
\definecolor{coral3}{rgb}{0.804,0.357,0.271}
\definecolor{coral4}{rgb}{0.545,0.243,0.184}
\definecolor{tomato1}{rgb}{1.000,0.388,0.278}
\definecolor{tomato2}{rgb}{0.933,0.361,0.259}
\definecolor{tomato3}{rgb}{0.804,0.310,0.224}
\definecolor{tomato4}{rgb}{0.545,0.212,0.149}
\definecolor{OrangeRed1}{rgb}{1.000,0.271,0.000}
\definecolor{OrangeRed2}{rgb}{0.933,0.251,0.000}
\definecolor{OrangeRed3}{rgb}{0.804,0.216,0.000}
\definecolor{OrangeRed4}{rgb}{0.545,0.145,0.000}
\definecolor{red1}{rgb}{1.000,0.000,0.000}
\definecolor{red2}{rgb}{0.933,0.000,0.000}
\definecolor{red3}{rgb}{0.804,0.000,0.000}
\definecolor{red4}{rgb}{0.545,0.000,0.000}
\definecolor{DeepPink1}{rgb}{1.000,0.078,0.576}
\definecolor{DeepPink2}{rgb}{0.933,0.071,0.537}
\definecolor{DeepPink3}{rgb}{0.804,0.063,0.463}
\definecolor{DeepPink4}{rgb}{0.545,0.039,0.314}
\definecolor{HotPink1}{rgb}{1.000,0.431,0.706}
\definecolor{HotPink2}{rgb}{0.933,0.416,0.655}
\definecolor{HotPink3}{rgb}{0.804,0.376,0.565}
\definecolor{HotPink4}{rgb}{0.545,0.227,0.384}
\definecolor{pink1}{rgb}{1.000,0.710,0.773}
\definecolor{pink2}{rgb}{0.933,0.663,0.722}
\definecolor{pink3}{rgb}{0.804,0.569,0.620}
\definecolor{pink4}{rgb}{0.545,0.388,0.424}
\definecolor{LightPink1}{rgb}{1.000,0.682,0.725}
\definecolor{LightPink2}{rgb}{0.933,0.635,0.678}
\definecolor{LightPink3}{rgb}{0.804,0.549,0.584}
\definecolor{LightPink4}{rgb}{0.545,0.373,0.396}
\definecolor{PaleVioletRed1}{rgb}{1.000,0.510,0.671}
\definecolor{PaleVioletRed2}{rgb}{0.933,0.475,0.624}
\definecolor{PaleVioletRed3}{rgb}{0.804,0.408,0.537}
\definecolor{PaleVioletRed4}{rgb}{0.545,0.278,0.365}
\definecolor{maroon1}{rgb}{1.000,0.204,0.702}
\definecolor{maroon2}{rgb}{0.933,0.188,0.655}
\definecolor{maroon3}{rgb}{0.804,0.161,0.565}
\definecolor{maroon4}{rgb}{0.545,0.110,0.384}
\definecolor{VioletRed1}{rgb}{1.000,0.243,0.588}
\definecolor{VioletRed2}{rgb}{0.933,0.227,0.549}
\definecolor{VioletRed3}{rgb}{0.804,0.196,0.471}
\definecolor{VioletRed4}{rgb}{0.545,0.133,0.322}
\definecolor{magenta1}{rgb}{1.000,0.000,1.000}
\definecolor{magenta2}{rgb}{0.933,0.000,0.933}
\definecolor{magenta3}{rgb}{0.804,0.000,0.804}
\definecolor{magenta4}{rgb}{0.545,0.000,0.545}
\definecolor{orchid1}{rgb}{1.000,0.514,0.980}
\definecolor{orchid2}{rgb}{0.933,0.478,0.914}
\definecolor{orchid3}{rgb}{0.804,0.412,0.788}
\definecolor{orchid4}{rgb}{0.545,0.278,0.537}
\definecolor{plum1}{rgb}{1.000,0.733,1.000}
\definecolor{plum2}{rgb}{0.933,0.682,0.933}
\definecolor{plum3}{rgb}{0.804,0.588,0.804}
\definecolor{plum4}{rgb}{0.545,0.400,0.545}
\definecolor{MediumOrchid1}{rgb}{0.878,0.400,1.000}
\definecolor{MediumOrchid2}{rgb}{0.820,0.373,0.933}
\definecolor{MediumOrchid3}{rgb}{0.706,0.322,0.804}
\definecolor{MediumOrchid4}{rgb}{0.478,0.216,0.545}
\definecolor{DarkOrchid1}{rgb}{0.749,0.243,1.000}
\definecolor{DarkOrchid2}{rgb}{0.698,0.227,0.933}
\definecolor{DarkOrchid3}{rgb}{0.604,0.196,0.804}
\definecolor{DarkOrchid4}{rgb}{0.408,0.133,0.545}
\definecolor{purple1}{rgb}{0.608,0.188,1.000}
\definecolor{purple2}{rgb}{0.569,0.173,0.933}
\definecolor{purple3}{rgb}{0.490,0.149,0.804}
\definecolor{purple4}{rgb}{0.333,0.102,0.545}
\definecolor{MediumPurple1}{rgb}{0.671,0.510,1.000}
\definecolor{MediumPurple2}{rgb}{0.624,0.475,0.933}
\definecolor{MediumPurple3}{rgb}{0.537,0.408,0.804}
\definecolor{MediumPurple4}{rgb}{0.365,0.278,0.545}
\definecolor{thistle1}{rgb}{1.000,0.882,1.000}
\definecolor{thistle2}{rgb}{0.933,0.824,0.933}
\definecolor{thistle3}{rgb}{0.804,0.710,0.804}
\definecolor{thistle4}{rgb}{0.545,0.482,0.545}
\definecolor{gray5}{rgb}{0.051,0.051,0.051}
\definecolor{gray10}{rgb}{0.102,0.102,0.102}
\definecolor{gray15}{rgb}{0.149,0.149,0.149}
\definecolor{gray20}{rgb}{0.200,0.200,0.200}
\definecolor{gray25}{rgb}{0.251,0.251,0.251}
\definecolor{gray30}{rgb}{0.302,0.302,0.302}
\definecolor{gray35}{rgb}{0.349,0.349,0.349}
\definecolor{gray40}{rgb}{0.400,0.400,0.400}
\definecolor{gray45}{rgb}{0.451,0.451,0.451}
\definecolor{gray50}{rgb}{0.498,0.498,0.498}
\definecolor{gray55}{rgb}{0.549,0.549,0.549}
\definecolor{gray60}{rgb}{0.600,0.600,0.600}
\definecolor{gray65}{rgb}{0.651,0.651,0.651}
\definecolor{gray70}{rgb}{0.702,0.702,0.702}
\definecolor{gray75}{rgb}{0.749,0.749,0.749}
\definecolor{gray80}{rgb}{0.800,0.800,0.800}
\definecolor{gray85}{rgb}{0.851,0.851,0.851}
\definecolor{gray90}{rgb}{0.898,0.898,0.898}
\definecolor{gray95}{rgb}{0.949,0.949,0.949}
\definecolor{gray100}{rgb}{1.000,1.000,1.000}
\definecolor{DarkGray}{rgb}{0.663,0.663,0.663}
\definecolor{DarkBlue}{rgb}{0.000,0.000,0.545}
\definecolor{DarkCyan}{rgb}{0.000,0.545,0.545}
\definecolor{DarkMagenta}{rgb}{0.545,0.000,0.545}
\definecolor{DarkRed}{rgb}{0.545,0.000,0.000}
\definecolor{LightGreen}{rgb}{0.565,0.933,0.565}

\newcommand{\hs}{\hspace{-1pt}}
\newcommand{\hhs}{\hspace{-3pt}}

\newcommand{\vs}{\vspace{-0pt}}

\usepackage[percent]{overpic}
\begin{document}
\allowdisplaybreaks
\title{\LARGE Local Optimality of Almost Piecewise-Linear Quantizers for Witsenhausen's Problem\vspace{-6pt}}



\author{Amir Ajorlou$^\dagger$$^\ddagger$ and Ali Jadbabaie$^\dagger$$^\ddagger$, \emph{Fellow}, \emph{IEEE}
\thanks{$^\dagger$Laboratory for Information and Decision Systems, Massachusetts Institute of Technology (MIT), Cambridge, MA 02139, USA. E-mail: \url{ajorlou@mit.edu}}
\thanks{$^\ddagger$Institute for Data, Systems, and Society, Massachusetts Institute of Technology (MIT), Cambridge, MA 02139, USA. E-mail: \url{jadbabai@mit.edu}}

\thanks{This work was supported by ARO MURI  W911NF-12-1-0509.}\vspace{-12pt}}
\maketitle
\thispagestyle{empty}

\begin{abstract}
We pose Witsenhausen's problem as a  leader-follower game of incomplete information. The follower makes a noisy observation of the leader's action (who moves first) and chooses an action minimizing her expected deviation from the leader's action. Knowing this, leader who observes the realization of the state, chooses an action that minimizes her distance to the state of the world and the ex-ante expected deviation from the follower's action. We study the perfect Bayesian equilibria of the game and identify a class of ``near piecewise-linear equilibria'' when leader cares much more about being close to the follower than the state, and the state is highly volatile.
As a major consequence of this result, we prove the existence of a set of local minima for Witsenhausen's problem in form of \textit{slopey} quantizers, which are at most a constant factor away from the optimal cost.
\end{abstract}

\vspace{-12pt}
\section{Introduction}
In his seminal work \cite{Witsen_68},  Witsenhausen constructed a  simple two-stage Linear-Quadratic-Gaussian (LQG) decentralized control problem where the optimal controller happens to be nonlinear.
This example showed for the first time that linear quadratic Gaussian team problems can have nonlinear solutions.
Using this counterexample, \cite{Basar_TAC_76}
produced an example showing that the standard decentralized static output feedback optimal control problem of linear deterministic systems could also admit nonlinear solutions.\footnote{Another relevant setting in which nonlinear equilibria may emerge is the seminal work of Crawford and Sobel (\cite{Sobel_82}) on signaling games where misaligned objectives of a sender and a receiver can result in quantized equilibrium strategies. Recently, authors in \cite{Serdar_TAC_17} consider an extension of this model to a noisy channel setup and show that for a Gaussian source and scalar signals
the equilibrium encoder is linear.}
For nearly half a century, this counterexample has been a subject of intense research across multiple communities (\cite{Ho_TAC_80, Basar_TAC_87, Sanjoy_99, Verdu_CDC_11, Grover_IT_15, Basar_SIAM_15}).

The endogenous information structure of  Witsenhausen's counterexample, where the signal observed in the second stage is a noisy version of the control action in the first stage,
gives rise to a nonclassical information structure. While the problem looks deceptively simple with quadratic cost, it is actually a very complicated, nonconvex, functional optimization problem.  This counterexample has shed light on intricacies  of optimal decisions in stochastic team optimization problems with similar information structure. Naturally, this problem has given rise to a large body of literature. For example, \cite{Basar_Birkhauser_13} provides a variant of Witsenhausen's counterexample with discrete primitive random variables and finite support, where no optimal solution exists. Another interesting variant, with the same information structure but different cost function, is the Gaussian test channel (\cite{Basar_TAC_87, Basar_CDC_08}) where the linear strategies can be shown to be optimal. Interestingly, \cite{Rotkowitz2006} shows that if the objective function in \cite{Witsen_68} is changed to a worst case induced norm, the linear controllers dominate nonlinear policies.

Although the optimal strategy and optimal cost for  Witsenhausen's counterexample are still unknown, it can be shown that carefully designed nonlinear strategies can largely outperform the linear strategies (see, e.g., the multi-point quantization strategies proposed by \cite{Sanjoy_99}). This result, in particular, implies the fragility of the comparative statics and policies solely derived based on the linear strategies in problems with similar setting.
A relevant line of research is to provide error bounds on the proximity to optimality for approximate solutions. \cite{Grover_IJSCC_2010,Grover_TAC_13} use information theoretic techniques and vector versions of the original problem to provide such bounds. In \cite{Yuksel_arxiv_16}, authors provide a general result on when one can approximate a continuous team decision problem with a finite one through quantized approximations, using which they show that quantized policies are asymptotically optimal for Witsenhausen's counterexample.
There are also several works aiming to approximate the optimal solution. \cite{Shamma_CDC,parisini_TAC_01,Lee_TAC_01,Rose_ISIT_2014} employ different heuristic approaches, all confirming what one might intuitively call an almost piecewise-linear form for the optimal controller. However, a complete optimality proof for such strategies has been elusive.

Although Witsenhausen's counterexample has been around for half a century, a little is known about
the topological properties of its optimal solution (e.g., whether it is continuous or not, the number of its fixed points, etc.).
A breakthrough in this direction is \cite{Verdu_CDC_11}, where  authors view Witsenhausen's setup in an optimal transport theory framework. This enables them not only to prove the existence of the optimal solution in a much more condense fashion, but also to derive some important characteristics of the optimal solution.
In particular, they show that the optimal controller is a strictly increasing function with a real analytic left inverse.\footnote{Note that this does not imply the continuity of the optimal solution.} As a consequence, no piecewise-linear strategy can be optimal, though not ruling out the possibility of a ``near piecewise-linear'' optimal solution.

In this paper, following Witsenhausen's original intuition, we view the problem as a leader-follower coordination game in which the action of the leader is corrupted by an additive noise, before reaching the follower. The leader  aims to coordinate with the follower while staying close to the observed state, recognizing that her action is not observed perfectly. As a result, she needs to  signal the follower in a manner that can be decoded efficiently.
More than a mere academic counterexample, the above setup  could model a scenario where coordination happens across generations and
the insights of the leader who is from a  different generation is corrupted/lost by the time the message reaches the future generations.
If the leader can internalize the fact that her actions will not be observed perfectly, how should she act to make sure coordination occur?
When the leader cares far more about coordination with the follower than staying ``on the message'', the  near piecewise-linear equilibrium strategy of the leader coarsens the observation in well-spaced intervals, rather than merely broadcasting a linearly scaled version of the observed state (as the linear strategy would suggest). 

To this end, we analyze the perfect Bayesian equilibria of this game and show that strong complementarity\footnote{Games of strategic complementarities are those in which the best response of each player
is increasing in actions of others \cite{Vives_90}.} between the leader and the follower combined with a prior with poor enough precision can give rise to nonlinear equilibria, and in particular,  equilibria in form of what has been deemed in the literature as {\it slopey quantizers} \cite{Grover_WiOPT_09}. We subsequently show that these equilibria are indeed local minima of the original Witsenhausen's problem.  Using some related results from asymptotic quantization theory (\cite{Panter_IRE_51, Loyd_IT_82,Na_IT_01}) together with analytical lower bounds on the optimal cost of Witsenhausen's problem derived in \cite{Grover_TAC_13}, we further show that these local minima include near-optimal solutions in the sense that their corresponding cost is at most a constant factor away from the optimal one.
Our work thus provides an analytical support for the local optimality of slopey quantization strategies for Witsenhausen's counterexample for a highly volatile state.

The main idea  behind  the proof is to carefully construct a class of what we informally refer to as {\em near piecewise-linear} or {\em slopey quantization} strategies for the leader that stays invariant under the best response operator.
These strategies can be viewed as \emph{small-slope variations} of a \emph{fixed-rate scalar quantizer} minimizing the mean squared quantization error (\cite{Panter_IRE_51, Loyd_IT_82,Na_IT_01}). Such an optimal quantizer is characterized by optimality conditions on the threshold levels which determine the boundaries of the quantization cells (or segments) and quantization levels: i) quantization levels must be the centroid of the segments, and ii) thresholds in between two adjacent quantization levels must be equidistant from them.
For any fixed number of segments, we consider the strategies whose segments are in a vicinity of the optimal MSE quantizer, have a unique fixed point in each segment close to the quantization level, and are almost linear within each segment with a near-zero derivative.
For such strategies, leader's actions remain very close to fixed points  of the strategy in each segment.
Therefore, well-spaced fixed points (combined with appropriate relative prior of the state in different segments) reveal the leader's actions to the follower with high probability, making the ``signal'' easily decodable. As a consequence, we can characterize the best response of the follower to leader's strategy. Using this characterization, we show that the best response of the leader to follower's strategy also varies very little, essentially remaining  near piecewise-linear over most of the range of the observed signals.

A key challenge in deriving the invariance property for this set of strategies for the leader is to bound and tightly control the displacement in the fixed points and endpoints of the segments of leader's strategy under the best response operator.
One major observation here is that the fixed points of the leader's best responses are {\it local minimizers} of the expected deviation  of the leader's action from the follower (which is known to be a non-convex functional \cite{Verdu_CDC_11}).
This insight allows us to show that the fixed points of the leader's best response lie in a tight neighborhood of the fixed points of the follower's strategy.
We then show that the fixed points of the follower's strategy in turn lie in a vicinity of a convex combination of the leader's fixed points and the expected value of the state of the world within each segment. Combining the two, we can derive an approximate dynamics for the displacement in the fixed points and endpoints of the segments in leader's strategy under the best response.
Using this approximate dynamics, we then characterize an invariant set of fixed points and interval endpoints for leader's strategy, which we can then use in order to prove the existence of a near piecewise-linear equilibrium strategy for the leader.

\section{Model}
\label{sec::model}
We view Witsenhausen's problem (\cite{Witsen_68}) as a game between a leader $L$ and a follower $F$. Before the agents act, the state of the world $\theta$ is drawn from a normal distribution with zero mean and variance $\sigma^2$. The  leader can observe the realization of $\theta$ and acts first. The payoff of the leader is given as follows
\begin{equation}\vs
\label{eq::uL}
u_L=-r_L (\theta-a_L)^2-(1-r_L)(a_F- a_L)^2,
\vs\end{equation}
where $a_F$ is the action of the follower and $0<r_L<1$. The follower makes a private, noisy observation of the leader's action, $s=a_L+\delta$ where $\delta\sim N(0,1)$. The payoff of the follower is given by
\begin{equation}\vs
\label{eq::uF}
u_F=-(a_L-a_F)^2.
\vs\end{equation}

We consider the perfect Bayesian equilibria of the game and show that they reduce to the Bayes Nash equilibria due to the Gaussian noise in the observation.\footnote{See, e.g., \cite{Tirole_91} for a definition of perfect Bayesian equilibrium and Bayes Nash equilibrium.} Denote with $a_L^*(\theta)$ and $a_F^*(s)$ the equilibrium strategies, and with $\nu^*(\cdot|s)$ the follower's belief on leader's action given $s$. Due to the normal noise in the observation, $\nu^*(\cdot|s)$ is fully determined by $a_L^*(\theta)$ and the prior as there are  no  off-equilibrium-path information sets. Equilibrium strategies should thus satisfy
\begin{align}\vs
\label{eq::PBE}
a_F^*(s)=& \E_{\nu^*}[a_L^*|s]=\int_{-\infty}^{\infty}{a_L\nu^*(a_L|s)da_L},\nonumber\\
a_L^*(\theta)=&\argmax_{a_L} -r_L (\theta-a_L)^2\nonumber\\
&-(1-r_L)\int_{-\infty}^{\infty}(a_F^*(s)- a_L)^2\phi(s-a_L)ds,
\vs\end{align}
where $\phi(\cdot)$ denotes the standard normal density function.

Our model yields the original setup in \cite{Witsen_68} by choosing $\frac{r_L}{1-r_L}=k^2$. The expected control cost then maps to the (negated) expected payoff of the leader. It is a simple exercise to find the optimal solution to Witsenhausen's problem in the class of linear strategies (see Lemma~11 in \cite{Witsen_68}), which is also an equilibrium of the game described above. Witsenhausen (\cite{Witsen_68}) showed that, for sufficiently large $\sigma$, this linear solution is not optimal.
In fact, the linear solution can be extremely suboptimal in the sense that the asymptotic ratio of the corresponding cost to the optimal one is infinity (\cite{Sanjoy_99}).
Our objective in this paper is to characterize a set of local minima for the problem in \cite{Witsen_68}, with a near piecewise-linear strategy for the leader and a cost within a constant factor of the optimal one, given a sufficiently large $\sigma$.
To this end, we analyze the equilibria of the game described above in regime $\frac{1}{2}\leq r_L\sigma^2\leq1$ and sufficiently large $\sigma$.\footnote{This clearly covers the case $k^2\sigma^2=1$.}

\vspace{-8pt}
\section{Nonlinear Equilibria}
\label{sec::nonlineareq}
We first prove the existence of a collection of equilibria with a near piecewise-linear  strategy for the leader for sufficiently large values of $\sigma$. Our approach is to identify a set of such strategies for the leader which is invariant under the best response operator. We characterize such a set in the next section.

Given $m\in\mathbb{N}$, consider a partition of the normal distribution $N(0,\sigma^2)$ into $2m+1$ segments $\cup_{k=-m}^{m} B_k^{\rm Q}$, with $B_k^{\rm Q}=[b_k^{\rm Q},b_{k+1}^{\rm Q})$ for $k\in\mathbb{N}_m$, $B_0^{\rm Q}=(b_{-1}^{\rm Q},b_1^{\rm Q})$, and $B_{-k}^{\rm Q}=(b_{-k-1}^{\rm Q},b_{-k}^{\rm Q}]$, with $b_{-k}^{\rm Q}=-b_k^{\rm Q}$ and $b_{m+1}^{\rm Q}=-b_{-m-1}^{\rm Q}=+\infty$. Denote with $c_k^{\rm Q}$ the centroid of segment $B_k^{\rm Q}$, that is, $c_k^{\rm Q}=\E_{N(0,\sigma^2)}[\theta|\theta\in B_k^{\rm Q}]$. Clearly, $c_0^{\rm Q}=0$ and $c_{-k}^{\rm Q}=-c_k^{\rm Q}$ for $k\in\mathbb{N}_m$.
We now specifically focus on a partition where the interval endpoints $b_k^{\rm Q}$ are equidistant from the centroids adjacent to them, i.e., $b_k^{\rm Q}=\frac{c_{k-1}^{\rm Q}+c_k^{\rm Q}}{2}$ for $k\in\mathbb{N}_m$. We can show that such a partition exists and is unique. This partition in fact corresponds to the $(2m+1)$-level fixed-rate scalar quantizer  that minimizes the mean-square distortion for a source characterized by $\theta\sim N(0,\sigma^2)$ (\cite{Panter_IRE_51, Loyd_IT_82,Na_IT_01}). The properties of this quantizer as $m\rightarrow\infty$ are extensively studied in asymptotic quantization theory, as will be discussed and used in analyzing the asymptotic performance of our proposed local minima in Section~\ref{sec::localOptima}.


\begin{figure}
\vspace{-12pt}
\centering
\begin{overpic}[scale=0.48]{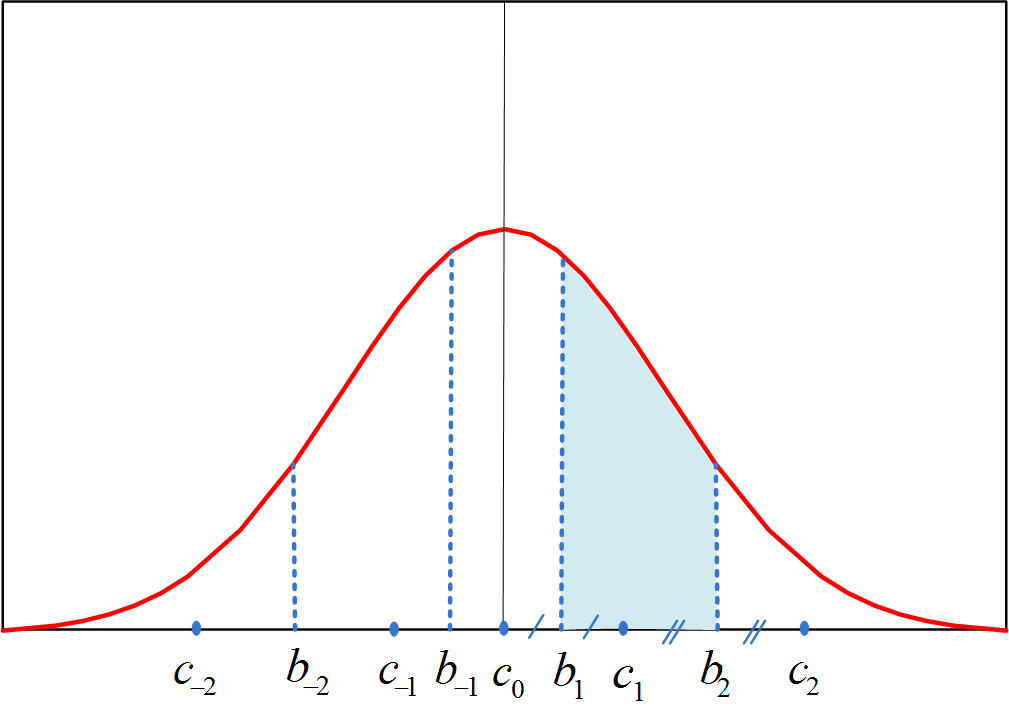}
\put(55.3,15){{$c_1^{\rm Q}=\atop{
\E[\theta|\theta\in B_1^{\rm Q}]}$}}
\put(56.5,4.0){${}^{\rm Q}$}
\put(62.3,4.0){${}^{\rm Q}$}
\put(71.3,4.0){${}^{\rm Q}$}
\put(79.8,4.0){${}^{\rm Q}$}

\put(50.3,4.0){${}^{\rm Q}$}
\put(45.3,4.0){${}^{\rm Q}$}
\put(39.3,4.0){${}^{\rm Q}$}
\put(30.3,4.0){${}^{\rm Q}$}
\put(18.8,4.0){${}^{\rm Q}$}

\end{overpic}
\vspace{-6pt}
\caption{Partition of the normal distribution in a $(2m+1)$-level optimal MSE quantizer, for $m=2$.}
\vspace{-10pt}
\end{figure}
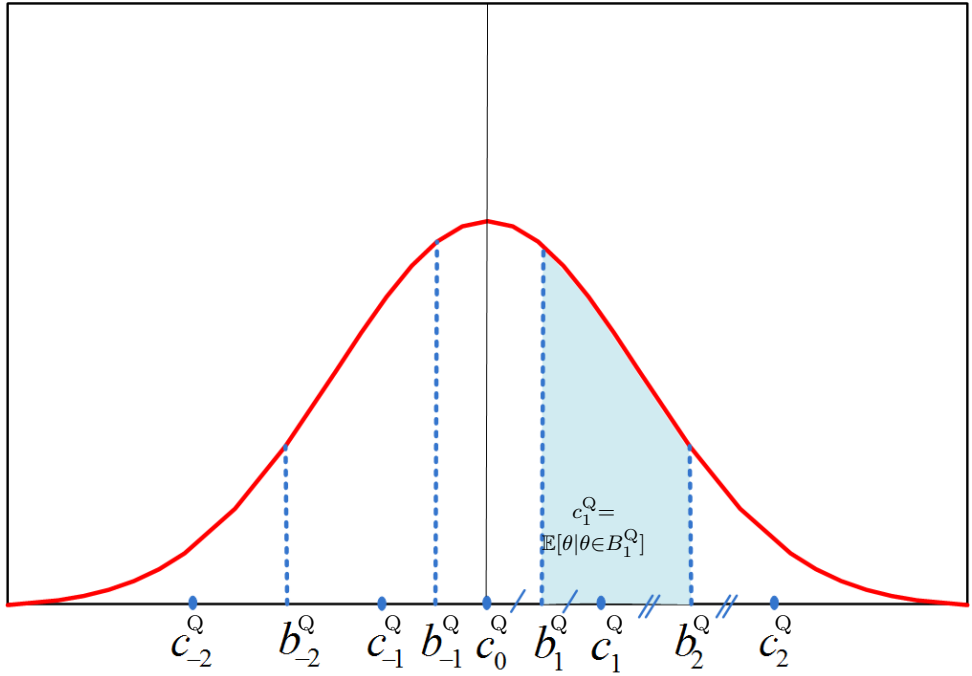

Roughly speaking, the set of strategies we propose for the leader are a class of $(2m+1)$-segmented strategies with segments being close to $B_k^{\rm Q}$ ($-m\leq k\leq m$), with a fixed point in each segment in a certain vicinity of $c_k^{\rm Q}$ ($-m\leq k\leq m$), and almost linear with a slope close to $r_L$ over each segment.
Before proceeding further, we present some (non-asymptotic) properties of this base configuration, which will facilitate the proof of the invariance property for the proposed set of strategies.
\begin{Lemma}
\label{lemma::baseconfig}
Given $m\in\mathbb{N}$, consider the partition of the normal distribution $N(0,\sigma^2)$ in a $(2m+1)$-level optimal MSE quantizer as described above. Define the $k$-th half-step of the quantizer as $x_k^{\rm Q}=\frac{c_k^{\rm Q}-c_{k-1}^{\rm Q}}{2}$, for $1\leq k\leq m$. Then,
\begin{enumerate}
\item[i)] $1-(\frac{x_m^{\rm Q}}{\sigma})^2\leq \frac{c_m^{\rm Q} x_m^{\rm Q}}{\sigma^2}\leq1$. If $m\geq2$, then $\frac{3}{4}\leq \frac{c_m^{\rm Q} x_m^{\rm Q}}{\sigma^2}\leq1$.
\item[ii)] For $1\leq k<m$,
\begin{align}\vs
\label{eq::baseii}
    \frac{\phi(\frac{b_k^{\rm Q}}{\sigma})}{\phi(\frac{c_k^{\rm Q}}{\sigma})}\leq(\frac{x_{k+1}^{\rm Q}}{x_k^{\rm Q}})^2\leq\frac{\phi(\frac{b_k^{\rm Q}}{\sigma})}{\phi(\frac{b_{k+1}^{\rm Q}}{\sigma})}.
\vs\end{align}
As a result, $1\leq\frac{x_{k+1}^{\rm Q}}{x_k^{\rm Q}}\leq e$, for $1\leq k<m$.
\item[iii)] For $0\leq j\leq k\leq m$, we have $\frac{\Prob[\theta|\theta\in B_k^{\rm Q}]}{\Prob[\theta|\theta\in B_j^{\rm Q}]}\leq\frac{1+e}{2}$.
\item[iv)] For any $m\geq5$,
\begin{align}\vs
 \frac{\sqrt{3\pi}}{e\sqrt{10}m}\leq&\frac{x_1^{\rm Q}}{\sigma}\leq\frac{\sqrt{2\pi e}}{2m},\nonumber\\
 1.1\sqrt{\ln m}-\frac{1}{1.1\sqrt{\ln m}}\leq&\frac{c_m^{\rm Q}}{\sigma}\leq2\sqrt{2\ln m+1.4}+1.45,\nonumber\\
  \frac{1}{2\sqrt{2\ln m+6}}\leq&\frac{x_m^{\rm Q}}{\sigma}\leq\frac{1}{1.1\sqrt{\ln m}}.
\vs\end{align}
\end{enumerate}
\end{Lemma}
\noindent\textit{Proof.} \hyperref[proof::lemma::baseconfig]{See the appendix.}$\hfill\blacksquare$

Next, we construct a set of $(2m+1)$-segmented increasing odd functions, denoted by $A_L^m(r_L,\sigma)$ satisfying the following properties:
\begin{Property}
\label{propty1}
For every $a_L(\theta)\in A_L^m(r_L,\sigma)$, there exist $2m+1$ segments $B_k=[b_k,b_{k+1})$, for $k\in\mathbb{N}_m$, $B_0=(-b_1,b_1)$, and $B_{-k}=(b_{-k-1},b_{-k}]$, with $b_{m+1}=-b_{-m-1}=+\infty$ such that:
\begin{itemize}
\item $a_L(\theta)$ is increasing and odd (i.e., $a_L(-\theta)=-a_L(\theta)$), and is smooth over each interval.
\item $a_L(\theta)$ has a unique fixed point in each segment.  That is, for each interval $B_k$, $(-m\leq k \leq m)$, there exists a unique $c_k\in B_k$  such that $a_L(c_k)=c_k$, with $c_0=0$.
\end{itemize}
\end{Property}

We also impose the constraint that interval endpoints $b_k$ remain close to midpoints of $[c_{k-1},c_{k}]$ and that fixed points $c_k$ remain within certain vicinity of $c_k^{\rm Q}$'s.
\begin{Property}
\label{propty2}
For every $k\in\mathbb{N}_m$, $|b_k-\frac{c_{k-1}+c_k}{2}|\leq 0.1r_L$.
Moreover,
$|c_k-c_k^{\rm Q}|\leq2.9$.
\end{Property}
From the above property, if we define $\bar x_k=x_k^{\rm Q}+3$ and $\ubar x_{k}=x_k^{\rm Q}-3$ for $1\leq k\leq m$, then $\bar x_k$ and $\ubar x_k$ represent upper and lower bounds on the lengths of both half-intervals $[c_{k-1}, b_k]$ and $[b_k, c_{k+1}]$.

Finally, we impose a constraint on the slope of $a_L(\theta)$ in each interval, keeping the slope very close to $r_L$, as well as a linear bound on $a_L(\theta)$ in the tail. We impose the following property:
\begin{Property}
\label{propty3}
For every $-m<k<m$ and $\theta\in B_k$,
$\ubar{r}\leq \frac{d}{d\theta}a_L(\theta)\leq\bar{r}$,
where $\ubar r=r_L(1-0.5r_L^2\sigma^2)$ and $\bar r=r_L(1+0.5r_L^2\sigma^2)$.
For the tail interval $B_m$,
$\ubar{r}\leq \frac{d}{d\theta}a_L(\theta)\leq\bar{r}$ for $b_m<\theta<c_m+\sqrt{e}\sigma\bar x_{m}$. For $\theta>c_m+\sqrt{e}\sigma\bar x_{m}$  we have $a_L(\theta)\leq c_m+3r_L(\theta-c_m)$.\footnote{We state the properties (and in many cases the analysis) only for $\theta\geq0$. The counterpart for $\theta\leq0$ is immediate since the function is odd.}
\end{Property}

For any $\sigma>0$, define $M(\sigma)=\{m\in\mathbb{N}|x_1^{\rm Q}>2\sqrt{2\ln\sigma}+5,~m\geq 25\}$.\footnote{ As we will see in Section~\ref{sec::localOptima}, this ensures inclusion of local minima with an expected cost within a constant factor of the optimal cost.}
We then claim that the set of strategies $A_L^m(r_L,\sigma)$ for $m\in M(\sigma)$, characterized by Property~\ref{propty1}-\ref{propty3},
is invariant under the best response operator for sufficiently large values of $\sigma$ in the regime $\frac{1}{2}\leq r_L\sigma^2\leq1$.
We formally state this result here, and provide the proof which is based on the best response analysis carried out in {Section~\ref{sec::Bestresponse}, in the appendix.}

\begin{Theorem}
\label{theorem::invariance}
Consider the regime $\frac{1}{2}\leq r_L\sigma^2\leq1$ with $\sigma>0$. Then, the set of $(2m+1)$-segmented strategies $A_L^m(r_L,\sigma)$ for the leader characterized by Property~\ref{propty1}-\ref{propty3}, where $m\in M(\sigma)=\{m\in\mathbb{N}|x_1^{\rm Q}>2\sqrt{2\ln\sigma}+5,~m\geq 25\}$ and $\sigma\geq 300$, is nonempty and invariant under the best response operator.\footnote{Recall that $x_1^{\rm Q}$ is the first half-step in a $(2m+1)$-level optimal MSE quantizer.}  Moreover, the game described in Section~\ref{sec::model} has an equilibrium for which:
\begin{itemize}
\item[i)] $a_L^*(\theta,r_L,\sigma)\in A_L^m(r_L,\sigma)$, and
\item[ii)] $(a_L^*(\theta,r_L,\sigma),a_F^*(s)=\E_\delta[a_L^*|s])$ maximizes the expected payoff of the leader over all pair of strategies $(a_L(\theta,r_L,\sigma),a_F(s)=\E_\delta[a_L|s])$ where  $a_L(\theta,r_L,\sigma)\in A_L^m(r_L,\sigma)$.
\end{itemize}
\end{Theorem}
\noindent\textit{Proof.} \hyperref[proof::theorem::invariance]{See the appendix.}$\hfill\blacksquare$

We expect the above theorem to hold for much smaller values of $\sigma$ and much smaller number of levels ($2m+1$), by optimizing/tightening the bounds used in deriving this result.
While this is in principle doable, we believe this would sacrifice clarity given the tedious calculations required, and would substantially increase the length of the paper.

\vspace{-12pt}
\section{Local Minima and Asymptotic Performance Guarantees}
\label{sec::localOptima}
Let
\begin{align}\vs
\label{eq::cost1}
&U\hspace{-1pt}(a_L,a_F\hspace{-3pt})\hspace{-3pt}=\hspace{-3pt}-\E_{\theta}[u_L(\theta,a_L,a_F)]=\hspace{-3pt}r_L\hspace{-3pt}\hspace{-3pt}\int_{-\infty}^{\infty}\hspace{-3pt}(\theta-a_L(\theta))^2\frac{\phi(\frac{\theta}{\sigma})}{\sigma}d\theta\nonumber\\
&\hspace{-3pt}+\hspace{-3pt}(1-r_L)\hspace{-3pt}\int_{-\infty}^{\infty}\hspace{-3pt}\int_{-\infty}^{\infty}\hspace{-3pt}(a_F(s)- a_L(\theta))^2\phi(s-a_L(\theta))\frac{\phi(\frac{\theta}{\sigma})}{\sigma}dsd\theta,
\vs\end{align}
for any two measurable functions $a_L,a_F:\mathbb{R}\to\mathbb{R}$. As discussed in Section~\ref{sec::model}, $U(a_L,a_F)$ defined above maps to the expected cost of the original Witsenhausen's problem in \cite{Witsen_68}. The aim of this section is to study the performance of the equilibrium strategies characterized by Theorem~\ref{theorem::invariance} in view of the above cost function.
Being an equilibrium implies that the cost cannot be improved by changing one of the strategies $a_L^*$ or $a_F^*$ while keeping the other fixed, although this does not rule out the possibility of obtaining a lower cost by simultaneously changing both strategies.\footnote{We thank Anant Sahai for bringing this point into the authors' attention.}
With a bit of manipulation, however, we can show that the image of the strategy obtained from an infinitesimal variation in $a_L^*$ also lies within $A_L^m$, using which we can show that $(a_L^*,a_F^*)$ is indeed a local minimum of $U$.

\begin{Lemma}
\label{lemma::localoptima}
Any pair of equilibrium strategies $(a_L^*,a_F^*)$ characterized by Theorem~\ref{theorem::invariance}, where $a_L^*\in A_L^m(r_L,\sigma)$ and $a_F^*(s)=\E_\delta[a_L^*|s]$\footnote{Recall that $\delta\sim N(0,1)$ is the noise in the follower's observation.}, is a local minimum of the cost functional $U$  in \eqref{eq::cost1}.
\end{Lemma}
\noindent\textit{Proof.} \hyperref[proof::lemma::localoptima]{See the appendix.}$\hfill\blacksquare$

We now have the first main result of the paper: a near piecewise-linear strategy for the
leader (or first controller) that leads to a local minimum of Witsenhausen's problem. Although the importance of such strategies has been already noticed in the literature (\cite{Shamma_CDC,parisini_TAC_01,Lee_TAC_01,Rose_ISIT_2014}), no analytical result concerning the optimality of such strategies is reported in the literature.
Theorem~\ref{theorem::invariance} also
presents an important result in the context of two-stage games of incomplete information.

We next aim to evaluate the asymptotic performance of these local minima with respect to the optimal cost. According to Theorem~\ref{theorem::invariance}, $(a_L^*,a_F^*)$ is a minimizer of $U$ over the pair of strategies $(a_L,a_F)$ with $a_L\in A_L^m(r_L,\sigma)$. Therefore, we can use any other pair of strategies with the leader's strategy being in $A_L^m$ (for which it is easier to evaluate the cost) to find an upper bound for $U(a_L^*,a_F^*)$. For this purpose we use $U(a_L^*,a_F^*)\leq U(a_L^{\rm Q},a_F^{\rm Q})$, where $a_L^{\rm Q}$ is the piecewise-linear strategy with segments $B_k^{\rm Q}$ and fixed points $c_k^{\rm Q}$ specified in the base configuration in Section~\ref{sec::nonlineareq} and $\frac{d}{d\theta}a_L^{\rm Q}(\theta)=r_L$ over each interval, and $a_F^{\rm Q}$ is the optimal $(2m+1)$-level MSE quantizer (i.e., constant value of $c_k^{\rm Q}$ over segment $B_k^{\rm Q}$).
It is easy to see that $(\theta-a_L^{\rm Q}(\theta))^2=(1-r_L)^2(\theta-a_F^{\rm Q}(\theta))^2$. We can thus write $U(a_L^{\rm Q},a_F^{\rm Q})=r_L(1-r_L)^2D_L^{\rm Q}+(1-r_L)D_F^{\rm Q}$, with
\begin{align}\vs
\label{eq::DLF}
D_L^{\rm Q}=&\int_{-\infty}^{\infty}(\theta-a_F^{\rm Q}(\theta))^2\frac{\phi(\frac{\theta}{\sigma})}{\sigma}d\theta,\nonumber\\
D_F^{\rm Q}=&\int_{-\infty}^{\infty}\int_{-\infty}^{\infty}(a_F^{\rm Q}(s)- a_L^{\rm Q}(\theta))^2\phi(s-a_L^{\rm Q}(\theta))\frac{\phi(\frac{\theta}{\sigma})}{\sigma}dsd\theta.
\vs\end{align}
$D_F^{\rm Q}$ can be upper-bounded as $D_F^{\rm Q}\leq4\sqrt{\frac{2}{e}}\frac{(2-r_L)^2}{(1-r_L)^2}\phi(\frac{x_1^{\rm Q}}{\sqrt{2}})+r_L^2 D_L^{\rm Q}$ (see the proof of Lemma~\ref{lemma::DLF}). We can find the exact asymptotic value of $D_L^{\rm Q}$ using results from asymptotic quantization theory (\cite{Panter_IRE_51, Loyd_IT_82,Na_IT_01}): $D_L^{\rm Q}$ is the mean-square error of an optimal ($2m+1$)-level MSE quantizer for a source $\theta\sim N(0,\sigma^2)$ (see, e.g.,  \cite{Na_IT_01}). It is known that for large $m$, $D_L^{\rm Q}\approx \frac{c_\infty}{(2m+1)^2}$, where $c_\infty$ is the \emph{Panter-Dite constant} of a normal  source given by
\begin{align}\vs
c_\infty=\frac{1}{12}\left(\int_{-\infty}^{\infty} (\frac{\phi(\frac{\theta}{\sigma})}{\sigma})^{\frac{1}{3}}d\theta\right)^3=\frac{\sqrt{3}\pi}{2}\sigma^2.
\vs\end{align}
Another interesting exact asymptotic equality is $(2m+1)\frac{x_1^{\rm Q}}{\sigma}\approx\frac{\sqrt{6\pi}}{2}$ using which we can alternatively write $D_L^{\rm Q}\approx\frac{(x_1^{\rm Q})^2}{\sqrt{3}}$ as $m\to\infty$. We have the following lemma.

\begin{Lemma}
\label{lemma::DLF}
For the pair of equilibrium strategies $(a_L^*,a_F^*)$ characterized by Theorem~\ref{theorem::invariance}, where $a_L^*\in A_L^m(r_L,\sigma)$, $a_F^*(s)=\E_\delta[a_L^*|s]$, and $m \in M(\sigma)$ we have
\begin{align}\vs
\liminf_{m\to\infty}\frac{\frac{r_L(1-r_L)(x_1^{\rm Q})^2}{\sqrt{3}}+4\sqrt{\frac{2}{e}}\frac{(2-r_L)^2}{(1-r_L)}\phi(\frac{x_1^{\rm Q}}{\sqrt{2}})}{U(a_L^*,a_F^*)}\geq1.
\vs\end{align}
\end{Lemma}
\noindent\textit{Proof.} \hyperref[proof::lemma::DLF]{See the appendix.}$\hfill\blacksquare$

The above asymptotic upper bound on $U(a_L^*,a_F^*)$ is minimized when $x_1^{\rm Q}\approx 2\sqrt{2\ln\sigma}$ for large $\sigma$, with number of levels $(2m+1)\approx\frac{\sqrt{3\pi}\sigma}{4\sqrt{\ln\sigma}}$ and
yielding a cost $\approx\frac{8r_L\ln\sigma}{\sqrt{3}}$. Recalling that $M(\sigma)=\{m\in\mathbb{N}|x_1^{\rm Q}>2\sqrt{2\ln\sigma}+5, m\geq 25\}$, this implies the existence of a local minimum with near piecewise-linear strategy for the leader with a cost asymptotically as low as $\frac{8r_L\ln\sigma}{\sqrt{3}}$. To compare with the optimal solution, we  use the lower bounds on the optimal cost of Witsenhausen's problem derived in \cite{Grover_TAC_13}. The following lemma is an immediate result of Theorem~4 in \cite{Grover_TAC_13}.
\begin{Lemma}
\label{lemma::Sahai}
Denote with $U^*(\sigma)$ the minimum value of the cost functional $U(a_L,a_F)$ given by \eqref{eq::cost1} in the regime $r_L\sigma^2=1$. Then
\begin{align}\vs
\limsup_{\sigma\to\infty}\frac{\frac{\ln\sigma}{6\sigma^2}}{U^*(\sigma)}\leq1.
\vs\end{align}
\end{Lemma}
\noindent\textit{Proof.} \hyperref[proof::lemma::Sahai]{See the appendix.}$\hfill\blacksquare$

This lower bound is quite loose (as also pointed out by the authors in \cite{Grover_TAC_13}), but still serves our purpose of showing that our proposed set of local optima include solutions that are only a constant factor away from the optimal cost as $\sigma\to\infty$.\footnote{The ratio between the upper and lower bounds in \cite{Grover_TAC_13} is almost $100$. For the well-known case of $\sigma=5$ and $k^2\sigma^2=1$, the lowest known cost $\approx0.167$ is only 12.5 times the value obtained from the lower bound.} We summarize the main findings of this section in the theorem below.
\begin{Theorem}
\label{theorem::performance}
Any pair of equilibrium strategies $(a_L^*,a_F^*)$ characterized by Theorem~\ref{theorem::invariance}, where $a_L^*\in A_L^m(r_L,\sigma)$ and $a_F^*(s)=\E_\delta[a_L^*|s]$, is a local minimum of the cost functional $U$  in \eqref{eq::cost1}. Moreover, 
\begin{align}\vs
\label{eq::th2}
\liminf_{\sigma\to\infty}\frac{\frac{8r_L\ln\sigma}{\sqrt{3}}}{\min_{m\in M(\sigma)} U(a_L^*,a_F^*)}\geq1.
\vs\end{align}
In the regime $r_L\sigma^2=1$,
at least one of these local minima are less than 27.8 times away from the optimal value, as $\sigma\to\infty$.
\end{Theorem}
\noindent\textit{Proof.} \hyperref[proof::theorem::performance]{See the appendix.}$\hfill\blacksquare$

\vspace{-8pt}
\section{Best Response Analysis}
\label{sec::Bestresponse}
The objective of this section is to prove Theorem~\ref{theorem::invariance} on the invariance of the set of strategies $A_L^m(r_L,\sigma)$ (for $m\in M(\sigma)$), and the existence of an equilibrium  with the leader's strategy in this set. The first step in verifying the invariance of $A_L^m(r_L,\sigma)$ is to
characterize the best response of the follower $a_F(s)$ to the leader's strategy $a_L(\theta)\in A_L^m(r_L,\sigma)$. We can then use these properties to find the updated best response of the leader to $a_F(s)$, denoted by $\tilde a_L(\theta)$  and enforce its inclusion in $A_L^m(r_L,\sigma)$.

The follower's best response to the strategy of the leader $a_L(\theta)$ is the
expected action of the leader given the observation $s=a_L+\delta$, that is $a_F(s)=\E_{\delta}[a_L|s]$. Following a simple application of Bayes rule we can obtain
\begin{equation}\vs
\label{eq::aF}
a_F(s)=\frac{\int_{-\infty}^{\infty}{a_L(\theta)\phi(s-a_L(\theta))\phi(\frac{\theta}{\sigma})d\theta}}{\int_{-\infty}^{\infty}{\phi(s-a_L(\theta))\phi(\frac{\theta}{\sigma})d\theta}}.
\vs\end{equation}
Using this, we can easily show that $a_F(s)$ is analytic and increasing, with $\frac{d}{ds}a_F(s)=\Var[a_L|s]$ (see \cite{Witsen_68} for a proof).

In order to characterize $a_F(s)$, we start by estimating the expected action of the leader and its variance conditioned on the interval to which $\theta$ belongs.
Actions of the leader in interval $B_k$ ($k\neq\pm m$) are well-concentrated around $c_k$. In fact $a_L(\theta)\in[c_k-\bar{r}\bar x_{k},c_k+\bar{r}\bar x_{k+1}]$ for $\theta\in B_k$, from which the lemma below follows immediately.

\begin{Lemma}
\label{lemma::aFlocal}
For $0\leq k<m$, $|\E[a_L(\theta)|s,\theta\in B_k]-c_k|\leq \bar r \bar{x}_{k+1}$
and $\Var[a_L(\theta)|s,\theta\in B_k]\leq \bar r^2 (\frac{\bar x_{k}+\bar x_{k+1}}{2})^2$.
\end{Lemma}
\noindent\textit{Proof.} \hyperref[proof::lemma::aFlocal]{See the appendix.}$\hfill\blacksquare$

The analysis is a bit involved in the tail, since for $\theta>c_m$ the leader's actions are not in a bounded vicinity of $c_m$ anymore. However, we can derive several useful properties for the tail as well.

\begin{Lemma}
\label{lemma::aFlocal-tail}
Consider a tail observation by the leader (i.e., $\theta\in B_m$). Then,
\begin{equation}\vs
\E[a_L(\theta)|s,\theta\in B_m]-c_m\leq\bar r\bar {x}_{m+1},
\vs\end{equation}
for $s\leq c_m+\bar x_{m+1}$, where $\bar x_{m+1}=\sqrt{e}\bar x_m$. For $s>c_m+\bar x_{m+1}$, we have
\begin{equation}\vs
\E[a_L(\theta)|s,\theta\in B_m]-c_m\leq3r_L\sigma(s-c_m+1).
\vs\end{equation}
Also, $\E[a_L(\theta)|s,\theta\in B_m]-c_m\geq -\bar r\bar x_m$.
As for the variance,
\begin{equation}\vs
\Var[a_L\hhs(\theta)|s,\theta\hhs\in\hhs B_m]\hhs\leq\hhs\begin{cases}
\hhs\frac{1}{3}, \text{ for } s<c_{m-1}&\\
\hhs\frac{3}{4}\bar r^2(\hhs\frac{\bar x_m+\bar x_{m+1}}{2}\hhs)^2\hhs,\text{ for } c_{m-1}\hhs\leq\hhs s\hhs\leq\hhs c_m\hhs+\hhs\bar x_{m+1}&\\
\hhs2.5r_L^2\sigma^2(s\hhs-\hhs c_m)^2, \text{ for } s\hhs>c_m+\bar x_{m+1}.&
\end{cases}
\vs\end{equation}
\end{Lemma}
\noindent\textit{Proof.} \hyperref[proof::lemma::aFlocal-tail]{See the appendix.}$\hfill\blacksquare$

Let the signal observed by the follower be between $c_k$ and $c_{k+1}$, i.e., $s=c_k+\delta$ with $0\leq\delta\leq c_{k+1}-c_k$.
Then, we claim that the follower's posterior on $\theta$ given $s$ has a negligible probability out of the neighboring intervals $B_k\cup B_{k+1}$.

\begin{Lemma}
\label{lemma::p(theta|s)}
Let the observed signal by the follower be $s=c_k+\delta$, where $0\leq\delta\leq c_{k+1}-c_k$, with $k\geq0$. Then, for any $j\geq 1$,
\begin{align}\vs
\frac{\Prob[\theta\in B_{k-j}|s]}{\Prob[\theta\in B_{k}|s]}\leq  e^{-\frac{(c_{k}-c_{k-j})^2}{2}+3j+1}.
\vs\end{align}
Similarly,
\begin{align}\vs
\frac{\Prob[\theta\in B_{k+j+1}|s]}{\Prob[\theta\in B_{k+1}|s]}\leq e^{-\frac{(c_{k+j+1}-c_{k+1})^2}{2}+3j+1}.
\vs\end{align}
\end{Lemma}
\noindent\textit{Proof.} \hyperref[proof::lemma::p(theta|s)]{See the appendix.}$\hfill\blacksquare$

Using this lemma and the fact that the fixed points $c_k$ are well-spaced, we can show that the effect of the intervals other than $B_k$ and $B_{k+1}$ on $a_F(s)$ are negligible.
In order to characterize the follower's best response $a_F(s)$, we then need to focus only on the segments adjacent to the observed signal, and in particular figure out the weight of each of these two neighboring intervals in the follower's posterior on $\theta$. We do this in the following lemma.
\begin{Lemma}
\label{lemma::p(th|s)withlog}
Define
\begin{equation}\vs
m_{k+1}=\frac{c_k+c_{k+1}}{2}+\frac{1}{\Delta_{k+1}}\ln\left(\frac{\Prob[\theta\in B_{k}]}{\Prob[\theta\in B_{k+1}]}\right),
\vs\end{equation}
where $\Delta_{k+1}=c_{k+1}-c_k$.
Also, write the signal observed by the follower as $s=m_{k+1}+\delta$. Then, for $0\leq k<m-1$,
\begin{align}\vs
e^{-\Delta_{k\hspace{-1pt}+\hspace{-1pt}1}(\delta+\bar{r}\bar x_{k+1})-\frac{\bar r^2\bar x_{k\hspace{-1pt}+\hspace{-1pt}1}^2}{2}}\hspace{-5pt}\leq\hspace{-3.5pt}\frac{\Prob[\theta\hspace{-2.5pt}\in\hspace{-2.5pt} B_{k}|s]}{\Prob[\theta\hspace{-2.5pt}\in\hspace{-2.5pt} B_{k\hspace{-1.5pt}+\hspace{-1.5pt}1}|s]}\hspace{-2.5pt}\leq\hspace{-2.5pt} e^{\Delta_{k\hspace{-1pt}+\hspace{-1pt}1}(\bar{r}\bar x_{k\hspace{-1pt}+\hspace{-1pt}2}-\delta)+\frac{\bar r^2\bar x_{k\hspace{-1pt}+\hspace{-1pt}2}^2}{2}}\hspace{-5pt}.
\vs\end{align}
For the case involving the tail segment $B_m$,
\begin{equation}\vs
e^{-\Delta_{m}(\delta+\bar{r}\bar x_m)\hspace{-1pt}-\hspace{-1pt}\frac{\bar r^2\bar x_m^2}{2}}\hspace{-5pt}\leq\hspace{-3.5pt}\frac{\Prob[\theta\hspace{-3.5pt}\in\hspace{-3.5pt} B_{m\hspace{-1pt}-\hspace{-1pt}1}|s]}{\Prob[\theta\hspace{-3.5pt}\in\hspace{-3.5pt} B_{m}|s]}\hspace{-2.5pt}\leq\hspace{-2.5pt} 1.16e^{-\Delta_{m}(\delta-\bar{r}\bar x_m)\hspace{-1pt}+\hspace{-1pt}\frac{\bar r^2\bar x_m^2}{2}}\hspace{-5pt}.
\vs\end{equation}
\end{Lemma}
\noindent\textit{Proof.} \hyperref[proof::lemma::p(th|s)withlog]{See the appendix.}$\hfill\blacksquare$

It is worth mentioning that $m_{k+1}$ defined in the above lemma is quite close to the midpoint of $c_k$ and $c_{k+1}$. In particular, using Lemma~\ref{lemma::baseconfig} we can show that $-\frac{1.1}{\Delta_{k+1}}\hhs<\hhs m_{k+1}\hs-\hs\frac{c_k+c_{k+1}}{2}\hhs<\hhs\frac{2.3}{\Delta_{k+1}}$.\footnote{See the proof of Lemma~\ref{lemma::aL} for details.}
We can now characterize the best response of the follower $a_F(s)$ to the leader's strategy $a_L(\theta)\in A_L^m(r_L,\sigma)$ up to the first order.
\begin{Lemma}
\label{lemma::aF(s)}
Let $s=m_{k+1}+\delta$, with $c_k\leq s\leq c_{k+1}$. Then
\begin{align}\vs
a_F(s)&\geq c_k+\frac{\Delta_{k+1}}{1+1.17e^{-\Delta_{k+1}\delta}}-1.01\bar{r}\bar x_{k+2}\nonumber\\
a_F(s)&\leq c_k+\frac{1.17\Delta_{k+1}}{1.17+e^{-\Delta_{k+1}\delta}}+1.01\bar{r}\bar x_{k+2}.
\vs\end{align}
Also,
\begin{align}\vs
\label{eq::daF(s)bounds-precise}
0\hspace{-3pt}\leq\hspace{-3.5pt}\frac{d}{ds}\hspace{-1.5pt}a_F(\hspace{-1.5pt}s\hspace{-1.5pt})\hspace{-3pt}&\leq\hspace{-2.5pt} 1.17e^{-\Delta_{k\hspace{-1pt}+\hspace{-1pt}1}\hspace{-2pt}|\hspace{-1.1pt}\delta\hspace{-1.1pt}|}\hspace{-2.5pt}\Delta_{k\hspace{-1pt}+\hspace{-1pt}1}^2\hspace{-2.5pt}+\hspace{-2.5pt}1.01\bar{r}^2(\hspace{-2pt}\frac{\bar x_{k}\hspace{-2.5pt}+\hspace{-2.5pt}\bar x_{k+1}}{2}\hspace{-2pt})^2\text{\hspace{-2.5pt} for\hspace{-2.5pt} }\delta\hspace{-2.5pt}\leq\hspace{-2.5pt}-0.5\nonumber\\
0\hspace{-3pt}\leq\hspace{-3.5pt}\frac{d}{ds}\hspace{-1.5pt}a_F(\hspace{-1.5pt}s\hspace{-1.5pt})\hspace{-3pt}&\leq\hspace{-2.5pt} 1.17e^{-\Delta_{k\hspace{-1pt}+\hspace{-1pt}1}\hspace{-2pt}|\hspace{-1.1pt}\delta\hspace{-1.1pt}|}\hspace{-2pt}\Delta_{k\hspace{-1pt}+\hspace{-1pt}1}^2\hspace{-2.5pt}+\hspace{-2.5pt}1.01\bar{r}^2(\hspace{-2pt}\frac{\bar x_{k+1}\hspace{-2.5pt}+\hspace{-2.5pt}\bar x_{k+2}}{2}\hspace{-2pt})^2\text{\hspace{-2.5pt} for \hspace{-4.5pt} }\delta\hspace{-2.5pt}>\hspace{-2.5pt}-0.5.
\vs\end{align}
\end{Lemma}
\noindent\textit{Proof.} \hyperref[proof::lemma::aF(s)]{See the appendix.}$\hfill\blacksquare$
\begin{Corollary}
\label{coraF(s)bounds}
A useful consequence of Lemma~\ref{lemma::aF(s)} is that
\begin{align}\vs
\label{eq::aF(s)bounds}
a_F(s)&\geq c_{k+1}-1.17e^{-\Delta_{k+1}\delta}\Delta_{k+1}-1.01\bar{r}\bar x_{k+2}\nonumber\\
a_F(s)&\leq c_k+1.17e^{\Delta_{k+1}\delta}\Delta_{k+1}+1.01\bar{r}\bar x_{k+2},
\vs\end{align}
and
\begin{align}\vs
\label{eq::daF(s)bounds}
0\hspace{-3pt}\leq\hspace{-3.5pt}\frac{d}{ds}a_F(s)\hspace{-3pt}\leq\hspace{-2.5pt} 1.17e^{-\Delta_{k+1}|\delta|}\Delta_{k+1}^2+1.01\bar{r}^2\bar x_{k+2}^2,
\vs\end{align}
where $s=m_{k+1}+\delta$, with $c_k\leq s\leq c_{k+1}$.
\end{Corollary}

Note that the exponential terms in the above bounds vanish quite fast for large $\Delta_{k+1}$ and $|\delta|$. for small $|\delta|$, another useful upper bound on the derivative of $a_F(s)$ is
\begin{equation}\vs
\label{eq::daFgeneric}
\frac{d}{ds}a_F(s)\leq \frac{1}{4}(\Delta_{k+1}+2\bar{r}\bar x_{k+2})^2+0.01\bar r^2\bar x_{k+2}^2.
\vs\end{equation}
\begin{Corollary}
\label{cor1}
Let $s=m_{k+1}+\delta$, with $c_k\leq s\leq c_{k+1}$. Then,
\begin{align}\vs
\label{eq::aF(s)simple bounds}
c_k\hspace{-2pt}-\hspace{-2pt}1.1\bar{r}\bar x_{k+2}\hspace{-2pt}\leq\hspace{-2pt} a_F(s)\hspace{-2pt}\leq\hspace{0pt} &c_k+1.1\bar{r}\bar x_{k+2}\mbox{ for }\delta\hspace{-2pt}<\hspace{-2pt}-\frac{2\sqrt{2\ln\sigma}}{5}\nonumber\\
c_{k+1}\hspace{-2pt}-\hspace{-2pt}1.1\bar{r}\bar x_{k+2}\hspace{-2pt}\leq\hspace{-2pt} a_F(s)\hspace{-2pt}\leq\hspace{0pt} &c_{k+1}+1.1\bar{r}\bar x_{k+2} \mbox{ for } \delta\hspace{-2pt}>\hspace{-2pt}\frac{2\sqrt{2\ln\sigma}}{5}.
\vs\end{align}
\end{Corollary}

Roughly speaking, the above corollary says that, if the observed signal by the follower is far enough from the midpoint of $c_k$ and $c_{k+1}$, then the optimal action of the follower is well-concentrated around $c_k$ or $c_{k+1}$ (whichever that is closer), and changes very slowly according to Lemma~\ref{lemma::aF(s)}.\footnote{Note that $\frac{2\sqrt{2\ln\sigma}}{5}<\frac{\ubar x_1}{5}<\frac{\Delta_{k+1}}{10}$} However, $a_F(s)$ may have very high variations for $s$ close to $m_{k+1}$ as can be seen from Lemma~\ref{lemma::aF(s)}.

The following lemma characterizes $a_F(s)$ when follower makes a tail observation.

\begin{Lemma}
\label{lemma::aF(s)-tail}
Let $s=c_m+\delta$, where $\delta>0$. Then,

\noindent i) for $\delta\leq \bar x_{m+1}$, $c_m-1.01\bar r\bar x_{m+1}\leq a_F(s)\leq c_m+\bar r\bar x_{m+1}$, and $0\leq\frac{d}{ds}a_F(s)\leq 0.8\bar r^2(\frac{\bar x_m+\bar x_{m+1}}{2})^2$.\\
\noindent ii) for $\delta>\bar x_{m+1}$, $c_m-1.01\bar r\bar x_{m}\leq a_F(s)\leq c_m+3r_L\sigma(\delta+1)$, and $0\leq\frac{d}{ds}a_F(s)\leq3r_L^2\sigma^2\delta^2$.
\end{Lemma}
\noindent\textit{Proof.} \hyperref[proof::lemma::aF(s)-tail]{See the appendix.}$\hfill\blacksquare$

Lemma~\ref{lemma::aF(s)} and \ref{lemma::aF(s)-tail} provide the first order characteristics of the best response of the follower to a leader's strategy $a_L(\theta)\in A_L^m(r_L,\sigma)$. We are now ready to analyze the leader's best response $\tilde a_L(\theta)$ to $a_F(s)$ and see if we can keep it in $A_L^m(r_L,\sigma)$. We have $\tilde a_L(\theta)\in\argmax_{a_L}\tilde u_L(\theta,a_L)$, where
\begin{align}\vs
\label{eq::uLBR}
\tilde u_L(\theta,a_L)=&-r_L(\theta\hspace{-2pt}-\hspace{-2pt}a_L)^2\nonumber\\
&-(1\hspace{-2pt}-\hspace{-2pt}r_L)\int_{-\infty}^{\infty}(a_F(s)- a_L)^2\phi(s-a_L)ds.
\vs\end{align}

\begin{Lemma}
\label{lemma::aL}
Consider $\theta\in[c_k,c_{k+1}]$, $0\leq k<m$. Then, there exists a unique $\tilde b_{k+1}\in[c_k,c_{k+1}]$ such that
\begin{align}\vs
\label{eq::aLBR1}
|\tilde a_L(\theta)-c_k|&<5\bar r\bar x_{k+2}\quad\mbox{for }\theta<\tilde b_{k+1},\nonumber\\
|\tilde a_L(\theta)-c_{k+1}|&<5\bar r\bar x_{k+2}\quad\mbox{for }\theta>\tilde b_{k+1}.
\vs\end{align}
\end{Lemma}
\noindent\textit{Proof.} \hyperref[proof::lemma::aL]{See the appendix.}$\hfill\blacksquare$

The points  $\tilde b_{k+1}$ determine the segments of the best response strategy $\tilde a_L(\theta)$.
The leader's best response strategy  clearly has a discontinuity at $\tilde b_{k+1}$. However, as we show in the next lemma, $\tilde a_L(\theta)$ is differentiable at all points $\theta\in[c_k,c_{k+1}]\setminus\{\tilde b_{k+1}\}$.
We further bound the derivative of $\tilde a_L(\theta)$ using the bounds on the follower's strategy and its derivative derived in Lemma~\ref{lemma::aF(s)}, Corollary~\ref{coraF(s)bounds}, and Corollary~\ref{cor1}.
\begin{Lemma}
\label{lemma::daLBR}
Consider $\theta\in[c_k,c_{k+1}]$, $0\leq k<m$, with $\theta\neq\tilde b_{k+1}$. Then,
\begin{align}\vs
\label{eq::daLBR}
\frac{d}{d\theta} \tilde a_L(\theta)&\geq\frac{r_L}{r_L+(1-r_L)(1+0.4\bar r^2\sigma^2)}\nonumber\\
\frac{d}{d\theta} \tilde a_L(\theta)&\leq\frac{r_L}{r_L+(1-r_L)(1-0.45\bar r^2\sigma^2)}.
\vs\end{align}
\end{Lemma}
\noindent\textit{Proof.} \hyperref[proof::lemma::daLBR]{See the appendix.}$\hfill\blacksquare$

Using this lemma and the values
$\ubar{r}=r_L(1-0.5r_L^2\sigma^2)$ and $\bar{r}=r_L(1+0.5r_L^2\sigma^2)$, we can easily verify that $\ubar{r}\leq\frac{d}{d\theta}\tilde a_L(\theta)\leq\bar{r}$. This means that Property~\ref{propty3} is preserved by the best response for $\theta\in[-c_m,c_m]$. We study the tail case later in Lemma~\ref{lemma::aLtail}.
Next lemma identifies $2m+1$ local minima of the MSE term in the leader's payoff, each located in a tiny interval around a fixed point of $a_L(\theta)$, establishing that they are indeed the fixed points of the best response strategy $\tilde a_L(\theta)$.

\begin{Lemma}
\label{lemma::J}
Define
\begin{equation}\vs
\label{eq::J}
\tilde J_L(a_L)=\int_{-\infty}^{\infty}(a_F(s)- a_L)^2\phi(s-a_L)ds.
\vs\end{equation}
Then, $\tilde J_L(a_L)$ is strongly convex over $[c_k-5\bar r\bar x_{k+2},c_k+5\bar r\bar x_{k+2}]$ with $\frac{d^2}{da_L^2}\tilde J_L(a_L)\geq 2(1-0.45\bar r^2\sigma^2)$. Let $\tilde c_k$ be the unique solution of
\begin{equation}\vs
\tilde c_k=\argmin_{a_L\in [c_k-5\bar r\bar x_{k+2},c_k+5\bar r\bar x_{k+2}]\cap [\tilde b_k,\tilde b_{k+1}]} \tilde J_L(a_L).
\vs\end{equation}
Then, $\tilde a_L(\tilde c_k)=\tilde c_k$.
\end{Lemma}
\noindent\textit{Proof.} \hyperref[proof::lemma::J]{See the appendix.}$\hfill\blacksquare$

The fixed point characterized in the above lemma is the unique fixed point in $[\tilde b_k,\tilde b_{k+1}]$ from Property~\ref{propty3}.
Therefore, Property~\ref{propty1} is also preserved under the best response. Next lemma describes the tail properties of $\tilde a_L(\theta)$.
\begin{Lemma}
\label{lemma::aLtail}
If $ \tilde b_m<\theta< \tilde c_m+\sigma\bar x_{m+1}$, then
$\ubar{r}\leq\frac{d}{d\theta}\tilde a_L(\theta)\leq \bar r$.
For $\theta> \tilde c_m+\sigma\bar x_{m+1}$, we have $\tilde a_L(\theta)\leq  \tilde c_m+3r_L(\theta-\tilde c_m)$.
\end{Lemma}
\noindent\textit{Proof.} \hyperref[proof::lemma::aLtail]{See the appendix.}$\hfill\blacksquare$

Now, in order to verify that the updated strategy $\tilde a_L(\theta)$ satisfies Property~\ref{propty2}, we need to bound the displacements in the fixed points $\tilde c_k$ and endpoints $\tilde b_k$.

\begin{Lemma}
\label{lemma::endpoints}
For the endpoints of the intervals corresponding to $\tilde a_L(\theta)$, we have $|\tilde b_{k+1}-\frac{\tilde c_k+\tilde c_{k+1}}{2}|\leq0.1r_L$.
\end{Lemma}
\noindent\textit{Proof.} \hyperref[proof::lemma::endpoints]{See the appendix.}$\hfill\blacksquare$



Bounding the displacement in $\tilde c_k$ can be done in multiple steps: first we need to relate the fixed point of the leader's best response
$\tilde a_L(\theta)$ in interval $\tilde B_k $ to the fixed point of $a_F(s)$ in $B_k$ (i.e., $s_k$), followed by estimating  $s_k$ in terms of $c_k$ and $e_k$ where  $e_k=\E_{N(0,\sigma^2)}[\theta|\theta\in B_k]$, that is, the expected value of $\theta$ over $B_k$. Finally we bound the displacement in $e_k$ with the displacement of the interval endpoints using  properties of truncated normal distribution.

\begin{Lemma}
\label{lemma::ck1}
Let $s_k$ be the fixed point of $a_F(s)$ in the interval $[c_k-5\bar r\bar x_{m+1},c_k+5\bar r\bar x_{m+1}]$, i.e., $a_F(s_k)=s_k$. Then,
\begin{equation}\vs
\label{eq::skck}
|\tilde c_k-s_k|\leq0.42r_L^2(\frac{\bar x_k+\bar x_{k+1}}{2})^2+0.08r_L^2\ubar x_1.
\vs\end{equation}
\end{Lemma}
\noindent\textit{Proof.} \hyperref[proof::lemma::ck1]{See the appendix.}$\hfill\blacksquare$

\begin{Lemma}
\label{lemma::sk}
$s_k$ can be located based on $c_k$ and $e_k$ as
\begin{equation}\vs
\label{eq::sk}
|s_k-(1-r_L)c_k-r_L e_k|\leq1.9 r_L^2\bar x_{k+1}.
\vs\end{equation}
\end{Lemma}
\noindent\textit{Proof.} \hyperref[proof::lemma::sk]{See the appendix.}$\hfill\blacksquare$

Using Lemma~\ref{lemma::endpoints}-\ref{lemma::sk}, we can reach at
\begin{equation}\vs
\label{eq::zeta}
|\tilde c_k-(1-r_L)c_k-r_L\hat e_k|\leq0.42 r_L^2(\frac{\bar x_k+\bar x_{k+1}}{2})^2+2r_L^2\bar x_{k+1},
\vs\end{equation}
where $\hat e_k=\E_{N(0,\sigma^2)}[\theta|\theta\in\hat B_k]$, with $\hat B_k=[\hat b_k,\hat b_{k+1}]$, $\hat b_k=\frac{c_{k-1}+c_{k}}{2}$ and $\hat b_{k+1}=\frac{c_k+c_{k+1}}{2}$.
We can now use \eqref{eq::zeta} and Lemma~\ref{lemma::endpoints} to verify that Property~\ref{propty2} is also preserved by the best response, completing the proof of the invariance of $A_L^m(r_L,\sigma)$ for $m\in M(\sigma)$ and $\sigma\geq 300$ in the regime $\frac{1}{2}\leq r_L\sigma^2\leq1$. This is carried out in the proof of the following theorem.
%

\vspace{-12pt}
\section{Conclusions}
We studied Witsenhausen's counterexample in a leader-follower game setup where the follower makes noisy observations from the leader's action and aims to choose her action as close as possible to that of the leader. Leader who moves first and can see the realization of the state of the world chooses her action to minimize her ex-ante distance from the follower's action as well as the state of the world. We showed the existence of  nonlinear perfect Bayesian equilibria in the regime $\frac{1}{2}\leq r_L\sigma^2\leq1$, where the leader's strategy is a perturbed  near-piecewise-linear version  of an optimal MSE quantizer. We then proved that these equilibria are indeed local minima  of the original Witsenhausen's problem. Incorporating some relevant results from asymptotic quantization theory and lower bounds on the optimal cost of Witsenhausen's problem from the literature, we showed that the proposed local minima include solutions that are at most a constant factor away from the optimal one.

\vspace{-8pt}
\begin{appendix}
{\label{proof::lemma::baseconfig}\noindent\textit{Proof of Lemma~\ref{lemma::baseconfig}.}}
To avoid lengthy expressions, we introduce and work with normalized variables $\hat x_k=\frac{x_k^{\rm Q}}{\sigma}$, $\hat b_k=\frac{b_k^{\rm Q}}{\sigma}$, $\hat c_k=\frac{c_k^{\rm Q}}{\sigma}$, and $\hat B_k=\frac{B_k^{\rm Q}}{\sigma}$. All the expectations are then taken using $N(0,1)$ as the probability measure.
\begin{enumerate}
\item[i)] Let $\rho(h)=\frac{\phi(h)}{1-\Phi(h)}$ denote the \emph{Mill's ratio}. Then, $0\leq1-\rho(\rho-h)\leq(\rho-h)^2$ (see, e.g. \cite{Normal_Horrace} for a proof). This implies that $0\leq1-\hat c_m \hat x_m\leq\hat x_m^2$. Moreover, $\hat x_m$ is decreasing with $m$ and it follows from direct calculation that $\hat x_2=0.48<\frac{1}{2}$, hence completing the proof.

\item[ii)] Since $\hat c_k$ is the centroid of segment $\hat B_k$, we should have
\begin{align}\vs
\label{eq::appen-base1}
\int_{\hat b_k}^{\hat c_k}(\hat c_k-\theta){\phi({\theta})}d\theta=\int_{\hat c_k}^{\hat b_{k+1}}(\theta-\hat c_k){\phi({\theta})}d\theta.
\vs\end{align}
This, together with the fact that $\int_{\hat b_k}^{\hat c_k}(\hat c_k-\theta){\phi({\theta})}d\theta\leq \phi(\hat b_k)\int_{\hat b_k}^{\hat c_k}(\hat c_k-\theta)d\theta$ and $\int_{\hat c_k}^{\hat b_{k+1}}(\theta-\hat c_k){\phi({\theta})}d\theta\geq \phi(\hat b_{k+1})\int_{\hat c_k}^{\hat b_{k+1}}(\theta-\hat c_k)d\theta$ proves the RHS inequality in \eqref{eq::baseii}. To derive the LHS, let $p_k^1=\Prob[\theta|\hat b_k\leq\theta\leq \hat c_k]$ and $p_k^2=\Prob[\theta|\hat c_k\leq\theta\leq \hat b_{k+1}]$. Noting that $\phi({\theta})$ (for $\theta\geq0$) and $\hat c_k-\theta$ are decreasing with $\theta$, we apply algebraic Chebyshev inequality to \eqref{eq::appen-base1} to obtain $p_k^1\hat x_k\leq p_k^2\hat x_{k+1}$. On the other hand,
\begin{align}\vs
\frac{p_k^2}{p_k^1}\leq\frac{\hat x_{k+1}}{\hat x_k}\times\frac{\phi(\hat c_k)}{\phi(\hat b_k)}.
\vs\end{align}
Combining the two, we can derive the LHS in \eqref{eq::baseii}. It then immediately follows from the LHS inequality that $\hat x_{k}\leq \hat x_{k+1}$ for $1\leq k$. Applying the result of part (i) to the RHS inequality we can easily show that $\frac{\hat x_{k+1}}{\hat x_k}\leq e$.

\item[iii)] A useful property here is that $g(x)=\frac{1}{x}\int_{a}^{a+x}\phi(t)dt$ is decreasing in $x$ for $x,a>0$. Using this, we can obtain
\begin{align}\vs
&\frac{\Prob[\theta|\theta\in \hat B_k]}{\Prob[\theta|\theta\in \hat B_j]}\leq\frac{\hat x_k+\hat x_{k+1}}{\hat x_j+\hat x_{j+1}}\times\frac{\phi(\hat b_k)}{\phi({\hat b_j})}\nonumber\\
&\leq\frac{\hat x_k+\hat x_{k+1}}{\hat x_j+\hat x_{j+1}}\times\frac{\hat x_j^2}{\hat x_k^2}\leq\frac{1+e}{2}\frac{\hat x_j}{\hat x_k}\leq\frac{1+e}{2},
\vs\end{align}
where the last two lines follow from part (ii). We have to modify the proof for $k=m$. To extend the proof to the case $k=m$, it suffices to show that ${\Prob[\theta|\theta\in \hat B_m]}\leq{\Prob[\theta|\theta\in \hat B_{m-1}]}$. We first note that
\begin{align}\vs
{\Prob[\theta|\theta\in \hat B_m]}=\frac{\phi({\hat b_m})}{{\hat c_m}}\leq\frac{4}{3}\phi({\hat b_m}){\hat x_m},
\vs\end{align}
where we have used part (i) in the last inequality. The proof now follows from the fact that ${\Prob[\theta|\theta\in \hat B_{m-1}]}>({\hat x_{m-1}+\hat x_m})\phi({\hat b_m})\geq (1+\frac{1}{e}){\hat x_m}\phi({\hat b_m})$.
\item[iv)] We start by showing that
\begin{align}\vs
\label{eq::appen-partiirefiend}
    e^{\frac{1}{2}(\hat x_1+\ldots+\hat x_{k})^2}\leq\frac{\hat x_{k+1}}{\hat x_1}\leq e^{\frac{5}{6}(\hat x_1+\ldots+\hat x_{k+1})^2}.
\vs\end{align}
The LHS easily follows from part (ii), while the RHS requires a more involved analysis as we elaborate below. The idea here is to find an appropriate lower bound for the RHS of \eqref{eq::appen-base1}. Using Jensen's inequality for the function $e^{-x}$, we can obtain
\begin{align}\vs
\int_{\hat c_k}^{\hat b_{k+1}}(\theta-\hat c_k){\phi({\theta})}d\theta\geq \frac{\hat x_{k+1}^2}{2\sqrt{2\pi}}e^{-\frac{1}{2}\int_{\hat c_k}^{\hat b_{k+1}}{\frac{2\theta^2(\theta-\hat c_k)}{\hat x_{k+1}^2}d\theta}}.
\vs\end{align}

\vspace*{-8pt}
Combining this with the same upper bound of $\frac{\hat x_k^2}{2}\phi(\hat b_k)$ as in part (ii) for the LHS of \eqref{eq::appen-base1} and after some simplification we can reach at
\begin{align}\vs
\frac{\hat x_{k+1}^2}{\hat x_k^2}\leq e^{\hat c_k\hat x_k+\frac{2}{3}\hat c_k\hat x_{k+1}-\frac{\hat x_k^2}{2}+\frac{\hat x_{k+1}^2}{4}}.
\vs\end{align}
Substituting $k$ with $1,\ldots,k-1$ and multiplying all these $k$ inequalities we can prove the RHS inequality in \eqref{eq::appen-partiirefiend}.

Incorporating the simple inequality $\frac{\phi(\hat x_1+\ldots+\hat x_k)}{\phi(\hat x_1+\ldots+\hat x_{k+1})}\leq e^{\frac{\hat c_m\hat x_m}{2}}\leq\sqrt{e}$ into \eqref{eq::appen-partiirefiend}, we can find that $\hat x_1\leq \sqrt{2\pi e}\hat x_{k+1}\phi(\hat x_1+\ldots+\hat x_{k+1})$ for $k=1,\ldots,m-1$.
Adding up all these inequalities and $\hat x_1\leq \sqrt{2\pi e}\hat x_{1}\phi(\hat x_1)$ yields
\begin{align}\vs
\frac{m\hat x_1}{\sqrt{2\pi e}}\leq\sum_{k=1}^{m}\hat x_{k}\phi(\hat x_1+\ldots+\hat x_{k})\leq\int_{0}^{\infty}\phi(\theta)d\theta=\frac{1}{2},
\vs\end{align}
proving $\hat x_1\leq \frac{\sqrt{2\pi e}}{2m}$.

Based on the RHS of \eqref{eq::appen-partiirefiend} and following a similar approach we can show that
\begin{align}\vs
\label{eq::appendp41}
\sum_{k=1}^{m-1}\hat x_{k+1}e^{-\frac{5}{6}(\hat x_1+\ldots+\hat x_{k})^2}\leq e(m-1)\hat x_1.
\vs\end{align}
On the other hand,
\begin{align}\vs
\label{eq::appendp42}
\int_{\hat x_1+\ldots+\hat x_m}^{\infty}e^{-\frac{5}{6}\theta^2}d\theta&\leq\frac{6e^{-\frac{5}{6}(\hat x_1+\ldots+\hat x_m)^2}}{10(\hat x_1+\ldots+\hat x_m)}\nonumber\\
&\leq\frac{6e^{-\frac{5}{6}(\hat x_1+\ldots+\hat x_m)^2}\hat x_m}{10(\hat x_1+\ldots+\hat x_m)\hat x_m}\leq1.6\hat x_1,
\vs\end{align}
where the last inequality follows from \eqref{eq::appen-partiirefiend} and part (i). Putting \eqref{eq::appendp41} and \eqref{eq::appendp42} together we can show that
$em\hat x_1\geq\int_{0}^{\infty}e^{-\frac{5}{6}\theta^2}d\theta=\sqrt{\frac{3\pi}{10}}$,
which yields $\hat x_1\geq\frac{\sqrt{3\pi}}{\sqrt{10}em}$.
Using \eqref{eq::appen-partiirefiend} for $k=m-1$, together with $\hat x_m\hat c_m\leq1$ and $\hat x_1\leq \frac{\sqrt{2\pi e}}{2m}$,
we can find
\begin{align}\vs
\label{eq::appen-iv2}
\frac{1}{2\sqrt{1.2}}\geq \hat x_m\sqrt{\ln\frac{2m\hat x_{m}}{\sqrt{2\pi e}}},
\vs\end{align}
using which we can show that
$\hat x_m\leq\frac{1}{1.1\sqrt{\ln m}}$, for $m\geq 5$.
This in turn implies that
\begin{align}\vs
\label{eq::appen-iv3}
\hat c_m\geq \frac{1}{\hat x_m}-\hat x_m\geq 1.1\sqrt{\ln m}-\frac{1}{1.1\sqrt{\ln m}}.
\vs\end{align}
Finally, using the LHS of \eqref{eq::appen-partiirefiend} for $k=m-1$, together with $\hat x_m\leq\frac{1}{1.1\sqrt{\ln m}}$ and $\hat x_1\geq\frac{\sqrt{3\pi}}{\sqrt{10}em}$,  and some manipulation we can obtain
\begin{align}\vs
\label{eq::appen-iv4}
\hat c_m\leq2\sqrt{2\ln m + 1.4}+\frac{2}{1.1\sqrt{\ln m}} \leq2\sqrt{2\ln m + 1.4}+1.45,
\vs\end{align}
for $m\geq5$. Along with $\hat c_m\hat x_m\geq 1-\hat x_m^2$, this leads to
\begin{align}\vs
\label{eq::appen-iv5}
\hat x_m\geq\frac{1}{\frac{\hat c_m}{2}+\sqrt{\frac{\hat c_m^2}{4}+1}}\geq
\frac{1}{2\sqrt{2\ln m+6}}.
\vs\end{align}
\end{enumerate}~$\hfill\blacksquare$

\begin{Remark}
\label{remark::1}
Particular consequences of assuming $m\geq 25$ and $\sigma\geq300$ and that $m$ is such that $x_1^{\rm Q}>2\sqrt{2\ln\sigma}+5$, are frequently used in the proofs concerning best response analysis. We summarize these properties here to avoid confusion in case they are not explicitly mentioned when used in the proofs.
\begin{enumerate}
\item $x_m^{\rm Q}>6x_1^{\rm Q}$ and $c_m^{\rm Q}>13.5x_m^{\rm Q}$. This follows from direct calculation of the optimal $(2m+1)$-level MSE quantizer for $m=25$ and that $\frac{x_m^{\rm Q}}{x_1^{\rm Q}}$ and $\frac{c_m^{\rm Q}}{x_m^{\rm Q}}$ are both increasing with $m$.
\item $x_m^{\rm Q}<0.262\sigma$. This follows from direct calculation of the optimal $(2m+1)$-level MSE quantizer for $m=25$ and that $\frac{x_m^{\rm Q}}{\sigma}$ is decreasing with $m$.
    As a result, $\bar x_m=x_m^{\rm Q}+3<0.272\sigma$ since $\sigma\geq 300$.
\item $x_1^{\rm Q}>11.5$. This follows from $x_1^{\rm Q}>2\sqrt{2\ln\sigma}+5$ and $\sigma\geq300$. As a result, $\ubar x_1=x_1^{\rm Q}-3>8.5$.
\item $m<\frac{\sigma}{2}$. This follows from part iv) of Lemma~\ref{lemma::baseconfig}:
$\frac{2m}{\sigma}\leq \frac{\sqrt{2\pi e}}{x_1^{\rm Q}}<\frac{\sqrt{2\pi e}}{2\sqrt{2\ln\sigma}+5}<1$.
\end{enumerate}
\end{Remark}

{\label{proof::lemma::aFlocal}\noindent\textit{Proof of Lemma~\ref{lemma::aFlocal}.}}
This is an immediate result of Property~\ref{propty2}.~$\hfill\blacksquare$

{\label{proof::lemma::aFlocal-tail}\noindent\textit{Proof of Lemma~\ref{lemma::aFlocal-tail}.}}
We start with the case where $a_L(c_m+\sigma\bar x_{m+1})\leq s\leq c_m+\bar x_{m+1}$. Let $\theta_c=c_m+\sigma\bar x_{m+1}$ and $\delta_c=s-a_L(\theta_c)$. With some manipulation, we can show that
for every $b_m\leq\theta,\theta'\leq\theta_c$,
\begin{equation}\vs
\label{eq::normalapprox}
\frac{\Prob[\theta'|s]}{\Prob[\theta|s]}=\frac{\phi(s-a_L(\theta'))\phi(\frac{\theta'}{\sigma})}{\phi(s-a_L(\theta))\phi(\frac{\theta}{\sigma})}\geq\frac{\phi(\frac{\theta'}{\sigma})\phi(\delta_c+\bar r (\theta_c-\theta'))}{\phi(\frac{\theta}{\sigma})\phi(\delta_c+\ubar r (\theta_c-\theta))}.
\vs\end{equation}
Integrating with respect to $\theta'$ and after some simplification, we arrive at
\begin{align}\vs
\label{eq::disttail1}
\Prob[\theta|s,b_m\leq\theta\leq \theta_c]&\leq \frac{\phi(\delta_c+\ubar r (\theta_c-\theta))}{\phi(\delta_c+\bar r (\theta_c-\theta))}\frac{\bar{\phi}(\theta)}{\bar{\Phi}(\theta_c)-\bar\Phi(b_m)}\nonumber\\
&\leq\xi\frac{\bar{\phi}(\theta)}{\bar{\Phi}(\theta_c)-\bar\Phi(b_m)},
\vs\end{align}
where $\bar{\phi}\sim N(\bar\mu,\bar\nu^2)$, with $\bar\mu=\frac{\bar{r}\sigma^2(\delta_c+\bar{r}\theta_c)}{1+\bar{r}^2\sigma^2}$ and $\bar\nu^2=\frac{\sigma^2}{1+\bar{r}^2\sigma^2}<\sigma^2$, and
\begin{equation}\vs
\label{eq::xi}
\xi=\frac{\phi(\delta_c+\ubar r (\theta_c-b_m))}{\phi(\delta_c+\bar r (\theta_c-b_m))}.
\vs\end{equation}
Using $\delta_c\leq \bar x_{m+1}-\ubar r\sigma\bar x_{m+1}$, it is easy to verify that $\delta_c+r_L(\theta_c-b_m)<\bar x_{m+1}+1$. It then follows that
\begin{align}\vs
\label{eq::xi}
\ln\xi&=(\bar r-\ubar r)(\theta_c-b_m)(\delta_c+r_L(\theta_c-b_m))\nonumber\\
&\leq r_L^3\sigma^2(\sqrt{e}\sigma+1)\bar x_{m}(\sqrt{e}\bar x_{m}+1)<0.001,
\vs\end{align}
for $\sigma\geq 300$ and $m\geq25$ (for which $x_{m}^{\rm Q}<0.262\sigma$),  implying that $\xi<1.01$.

A useful formula is
\begin{align}
\label{eq::varianceid}
\frac{1}{\sigma^2}{\rm Var}_{N(\mu,\sigma^2)}[\theta|c\leq\theta]&=1-\frac{(\E_{N(\mu,\sigma^2)}[\theta|c\leq\theta]-c)(\E_{N(\mu,\sigma^2)}[\theta|c\leq\theta]-\mu)}{\sigma^2}\nonumber\\
&\leq (\frac{\E_{N(\mu,\sigma^2)}[\theta|c\leq\theta]-c}{\sigma})^2,
\end{align}
for $\mu\leq c$. Moreover, the above variance is increasing in $\mu$ for $\mu\leq c$ (see, e.g. \cite{Normal_Horrace}).

Let $e_{\bar \mu}=\E_{N(\bar \mu,\sigma^2)}[\theta|b_m\leq\theta]-b_m$ and
$e_0=\E_{N(0,\sigma^2)}[\theta|b_m\leq\theta]-b_m$. Then, from the increasing property mentioned above, we find
\begin{align}
\label{}
&(b_m+e_{\bar \mu}-\bar \mu)e_{\bar \mu}\leq (b_m+e_0)e_0\Rightarrow\nonumber\\
&e_{\bar \mu}\leq\frac{-(b_m-{\bar \mu})+\sqrt{(b_m-{\bar \mu})^2+4e_0(b_m+e_0)}}{2}\nonumber\\
&\leq\frac{-(b_m-{\bar \mu})+\sqrt{(b_m-{\bar \mu})^2+4\bar x_m(b_m+\bar x_m)}}{2},
\end{align}
and hence,
\begin{align}
\label{}
\frac{e_{\bar \mu}}{\bar x_m}\leq\frac{2(b_m+\bar x_m)}{(b_m-{\bar \mu})+\sqrt{(b_m-{\bar \mu})^2+4\bar x_m(b_m+\bar x_m)}}.
\end{align}
Based on this inequality and that $x_1^{\rm Q}>11.5$ (from the definition of $M(\sigma)$ and that $\sigma\geq 300$), and that $x_m^{\rm Q}>6 x_1^{\rm Q}$ and $b_m^{\rm Q}\geq 12.5x_m^{\rm Q}$ for $m\geq25$, we can show that ${e_{\bar \mu}}\leq 1.134\bar x_m$.
Therefore,
\begin{align}\vs
\label{eq::aFtail1}
&\E[a_L(\theta)|s,b_m\leq\theta\leq \theta_c]\hhs-\hhs c_m\hhs\leq\hhs \xi\bar r \E_{\bar\phi}[\theta-b_m|b_m\leq\theta\leq\theta_c]\nonumber\\
&<\xi\bar r \E_{N(\bar\mu,\sigma^2)}[\theta-b_m|b_m\leq\theta]<
1.01\times1.134\bar r\bar x_m< 0.75\bar r\bar x_{m+1}.
\vs\end{align}

As for the variance, we start with
\begin{align}\vs
\label{eq::VaraFtail0}
&\Var[a_L(\theta)|s,b_m\leq\theta\leq\theta_c]\leq\hhs\xi\bar r^2\Var_{\bar\phi}[\theta|b_m\leq\theta\leq\theta_c]\hhs\leq\xi\bar r^2\Var_{N(\bar\mu,\sigma^2)}[\theta|b_m\hhs\leq\hhs\theta].
\vs\end{align}
Combining this with the bound in \eqref{eq::varianceid}, we get
\begin{align}\vs
\label{eq::VaraFtail1}
\Var[a_L(\theta)|s,b_m\leq\theta\leq\theta_c]\leq\hspace{-3pt}1.01\times 1.134^2\bar r^2\bar x_m^2<\hspace{-3pt} 1.3\bar r^2\bar x_m^2.
\vs\end{align}

Similar results to the above can be derived for the case where $a_L(b_m)\leq s<a_L(\theta_c)$, using $\theta_s$ instead of $\theta_c$, where $s=a_L(\theta_s)$ with $b_m\leq \theta_s<\theta_c$.
The same for the case $s<a_L(b_m)$, following a similar argument with $\bar{\phi_b}\sim N(\bar\mu_b,\bar\nu^2)$, where $\bar\mu_b=\frac{\bar{r}\sigma^2(\bar{r}b_m-\delta_b)}{1+\bar{r}^2\sigma^2}$, $\bar\nu^2=\frac{\sigma^2}{1+\bar{r}^2\sigma^2}$, and $\delta_b=a_L(b_m)-s>0$.


Now we bring into play the tail effect. For every $\theta\geq\theta_c$, we use
\begin{align}\vs
\frac{\Prob[\theta|s]}{\Prob[c_m\leq\theta'\leq\theta_c|s]}=\frac{\phi(s-a_L(\theta))\phi(\frac{\theta}{\sigma})}{\int_{c_m}^{\theta_c}\phi(s-a_L(\theta'))\phi(\frac{\theta'}{\sigma})d\theta'},
\vs\end{align}
using which for $s=c_m+\delta$ with $0\leq\delta\leq\bar x_{m+1}$, we get
\begin{align}\vs
\label{eq::taildistafterthc}
\frac{\Prob[\theta|s]}{\Prob[c_m\leq\theta'\leq\theta_c|s]}\leq\frac{e^{\frac{\delta^2}{2}}\phi(\frac{\theta}{\sigma})}{\sigma(\Phi(\frac{\theta_c}{\sigma})-\Phi(\frac{c_m}{\sigma}))}.
\vs\end{align}
Therefore, for $\theta\geq\theta_c$
\begin{equation}\vs
\label{eq::tailtmp1}
\Prob[\theta|s,\theta\geq c_m]\leq \frac{e^{\frac{\delta^2}{2}}(1-\Phi(\frac{\theta_c}{\sigma}))}{\Phi(\frac{\theta_c}{\sigma})-\Phi(\frac{c_m}{\sigma})}\times\frac{\phi(\frac{\theta}{\sigma})}{\sigma(1-\Phi(\frac{\theta_c}{\sigma}))}.
\vs\end{equation}
Using the inequality $h\leq\rho(h)\leq\frac{h^2+1}{h}$ for $h>0$, we can show that $\frac{1-\Phi(\frac{\theta_c}{\sigma})}{\Phi(\frac{\theta_c}{\sigma})-\Phi(\frac{c_m}{\sigma})}\leq\frac{\phi(\frac{\theta_c}{\sigma})}{\phi(\frac{c_m}{\sigma})}$. This, along with \eqref{eq::tailtmp1} and $0\leq\delta\leq\bar x_{m+1}$ and $\theta_c=c_m+\sigma\bar x_{m+1}$ yields
\begin{align}\vs
\label{eq::taildist}
&\Prob[\theta|s,\theta\in B_m]\leq\Prob[\theta|s,\theta\geq c_m]\nonumber\\
&\leq \frac{e^{-\frac{c_m\bar x_{m+1}}{\sigma}}\phi(\frac{\theta}{\sigma})}{\sigma(1-\Phi(\frac{\theta_c}{\sigma}))}\leq\frac{e^{-{\sigma}}\phi(\frac{\theta}{\sigma})}{\sigma(1-\Phi(\frac{\theta_c}{\sigma}))},
\vs\end{align}
for $\theta\geq\theta_c$, where we have used $c_m\bar x_{m+1}\geq \sigma^2$ (which easily follows from Lemma~\ref{lemma::baseconfig}). Using this along with \eqref{eq::aFtail1}, we can obtain
\begin{align}\vs
\label{tac::tmp3}
&\E[a_L(\theta)|s,\theta\in B_m]-c_m\nonumber\\
&\leq\hhs 0.75\bar r\bar x_{m+1}\hhs+\hhs e^{-{\sigma}}(3{r_L}(\theta_c-c_m+\sigma(\rho(\frac{\theta_c}{\sigma})-\frac{\theta_c}{\sigma})))\nonumber\\
&\leq\hhs 0.75\bar r\bar x_{m+1}\hhs+\hhs 3r_L\sigma e^{-{\sigma}}\hhs(\bar x_{m+1}\hhs+\hhs\frac{\sigma}{c_m\hhs+\hhs\sigma\bar x_{m+1}})\hhs\leq\hhs \bar r\bar x_{m+1}.
\vs\end{align}
Also, $\E[a_L(\theta)|s,\theta\in B_m]\geq a_L(b_m)\geq c_m-\bar r\bar x_m$.
To bound the variance, let $\kappa=\E[a_L(\theta)|s,b_m\leq\theta\leq\theta_c]$ ($\kappa>a_L(b_m)$). Then,
\begin{align}\vs
\label{eq::VaraFtail}
&\Var[a_L(\theta)|s,\theta\in B_m]\leq\E[(a_L(\theta)-\kappa)^2|s,\theta\in B_m]\nonumber\\
&\leq 1.3\bar r^2\bar x_m^2+e^{-{\sigma}}\E_{N(0,\sigma^2)}[(a_L(\theta)-a_L(b_m))^2|\theta\geq\theta_c]\nonumber\\
&\leq1.3\bar r^2\bar x_m^2+e^{-{\sigma}}\E_{N(0,\sigma^2)}[(3r_L(\theta-c_m)+\bar r\bar x_m)^2|\theta\geq\theta_c]\nonumber\\
&\leq1.3\bar r^2\bar x_m^2+9\bar r^2e^{-{\sigma}}((\sigma\bar x_{m+1}+1)^2+\sigma^2)\nonumber\\
&\leq0.75\bar r^2(\frac{\bar x_m+\bar x_{m+1}}{2})^2,
\vs\end{align}
on noting $\sigma\geq 300$.

For $s<c_{m-1}$, we use the fact that $\Var[\theta|s,b_m\leq\theta\leq\theta_c]\leq(\frac{\theta_c-b_m}{2})^2$ to get
\begin{equation}\vs
\label{eq::VaraFtailsleq-cm}
\Var[a_L(\theta)|s,b_m\leq\theta\leq\theta_c]\leq \bar r^2(\frac{\theta_c-b_m}{2})^2<0.3.
\vs\end{equation}
As for the effect of $\theta>\theta_c$, we can easily see that for $s<c_m$ \eqref{eq::taildist} becomes
\begin{equation}\vs
\label{eq::taildist2}
\Prob[\theta|s,\theta\in B_m]\leq \frac{e^{-\sigma-\frac{\bar x_{m+1}^2}{2}}\phi(\frac{\theta}{\sigma})}{\sigma(1-\Phi(\frac{\theta_c}{\sigma}))}.
\vs\end{equation}
We can use this to bound the variance similar to \eqref{eq::VaraFtail}:
\begin{equation}\vs
\label{eq::VaraFtaileffect}
\Var[a_L(\theta)|s,\theta\in B_m]\hhs\leq\hhs\left.
                                       \begin{cases}
                                         0.75\bar r^2(\frac{\bar x_m+\bar x_{m+1}}{2})^2,\hhs &\hhs {c_{m-1}\leq s\leq c_m}\\
                                         \frac{1}{3}, \hhs& \hhs{s<c_{m-1}}
                                       \end{cases}
                                     \right.
\vs\end{equation}

For the case where $s>c_m+\bar x_{m+1}$ (i.e., $\delta>\bar x_{m+1}$), let $\theta_s=c_m+\sigma\delta$. Then, similar to \eqref{eq::taildist} we can obtain
\begin{equation}\vs
\label{eq::taildistslarge}
\Prob[\theta|s,\theta\in B_m]\leq \frac{e^{-\sigma}\phi(\frac{\theta}{\sigma})}{\sigma(1-\Phi(\frac{\theta_s}{\sigma}))},
\vs\end{equation}
for $\theta\geq\theta_s$. Using this and similar to \eqref{tac::tmp3}, we can reach at
\begin{align}\vs
\E[a_L(\theta)|s,\theta\in B_m]< c_m+3r_L\sigma(\delta+1).
\vs\end{align}
To bound the variance, similar to \eqref{eq::VaraFtail} we can show
\begin{align}\vs
\label{eq::VaraFtailslarge}
\Var[a_L(\theta)|s,\theta\in B_m]<2.5r_L^2\sigma^2\delta^2,
 \vs\end{align}
which completes the proof.
$\hfill\blacksquare$


{\label{proof::lemma::p(theta|s)}\noindent\textit{Proof of Lemma~\ref{lemma::p(theta|s)}.}}
First, we use the properties of the base configuration listed in Lemma~\ref{lemma::baseconfig} to show that
\begin{align}\vs
\label{priors}
\frac{\Prob[\theta\in B_{k-j}]}{\Prob[\theta\in B_{k}]}\leq e^{2j+1}.
\vs\end{align}
Let $\tilde B_k\hs\hhs=\hs\hhs[\tilde b_k,\tilde b_{k+1}]\hs\hhs=\hs\hhs[\frac{c_{k-1}^{\rm Q}+c_{k}^{\rm Q}}{2}\hs\hhs+\hs\hhs3,\frac{c_{k}^{\rm Q}+c_{k+1}^{\rm Q}}{2}]$, and $\tilde B_{k-j}\hs\hhs=\hs\hhs[\tilde b_{k-j},\tilde b_{k-j+1}]\hs\hhs=\hs\hhs[\frac{c_{k-j-1}^{\rm Q}+c_{k-j}^{\rm Q}}{2}\hs\hhs-\hs\hhs3,\frac{c_{k-j}^{\rm Q}+c_{k-j+1}^{\rm Q}}{2}]$. Then, it is straightforward to verify that
\begin{align}\vs
&\frac{\Prob[\theta\in B_{k-j}]}{\Prob[\theta\in B_{k}]}\leq\frac{\Prob[\theta\in \tilde B_{k-j}]}{\Prob[\theta\in\tilde B_{k}]}.
\vs\end{align}
Similar to part iii) of Lemma~\ref{lemma::baseconfig}, we can write
\begin{align}\vs
&\frac{\Prob[\theta\in \tilde B_{k-j}]}{\Prob[\theta\in\tilde  B_{k}]}\leq \frac{\phi(\frac{\tilde b_{k-j+1}}{\sigma})}{\phi(\frac{\tilde b_{k+1}+\delta}{\sigma})}\times\frac{\tilde b_{k+1}-\tilde b_{k}+\delta}{\tilde b_{k+1}-\tilde b_{k}},
\vs\end{align}
where $\delta=\max\{0,(\tilde b_{k-j+1}-\tilde b_{k-j})-(\tilde b_{k+1}-\tilde b_{k})\}\leq 6$. As a result,
\begin{align}\vs
\label{eq::Probratio}
\frac{\Prob[\theta\in B_{k-j}]}{\Prob[\theta\in B_{k}]}&\leq\frac{\phi(\frac{c_{k-j}^{\rm Q}+c_{k-j+1}^{\rm Q}}{2\sigma})}{\phi(\frac{c_{k}^{\rm Q}+c_{k+1}^{\rm Q}+12}{2\sigma})}\times\frac{x_{k+1}^Q+x_{k}^Q-3+\delta}{x_{k+1}^Q+x_{k}^Q-3}\nonumber\\
&\leq 1.3e^{\frac{1}{\sigma^2}(c_m^{\rm Q}+3)(c_m^{\rm Q}-c_{m-j}^{\rm Q}+6)},
\vs\end{align}
where we have used the fact that $x_1^Q>11.5$.
Based on the above inequality, together with $c_m^{\rm Q}-c_{m-j}^{\rm Q}\leq 2jx_m^Q$, and $\frac{c_m^{\rm Q} x_m^{\rm Q}}{\sigma^2}\leq 1$, and that $\frac{c_m^{\rm Q}}{\sigma}\leq 2\sqrt{2\ln\sigma+1.5}+1.5$ (which follows from Lemma~\ref{lemma::baseconfig} and that $m<\sigma$ as stated in Remark~\ref{remark::1}), and a bit of manipulation we can show that $\frac{\Prob[\theta\in B_{k-j}]}{\Prob[\theta\in B_{k}]}\leq e^{2j+1}$.
%
%
%

Now, to prove the lemma for $k<m$, we write
\begin{align}\vs
\label{tac::tmp2}
&\frac{\Prob[\theta\in B_{k-j}|s]}{\Prob[\theta\in B_{k}|s]}\leq \frac{\Prob[\theta\in B_{k-j}]\phi(\delta+(c_k-c_{k-j})-\bar{r}\bar x_m)}{\Prob[\theta\in B_{k}]\phi(\delta+\bar{r}\bar x_m)}\nonumber\\
&\leq \frac{\Prob[\theta\in B_{k-j}]}{\Prob[\theta\in B_{k}]} e^{-\frac{(c_k-c_{k-j})^2}{2}+\bar r\bar x_m(c_k-c_{k-j})-\delta(c_k-c_{k-j}-2\bar r\bar x_m)}\nonumber\\
&\leq e^{-\frac{(c_k-c_{k-j})^2}{2}+2j+1+2j\bar r\bar x_m^2}\leq e^{-\frac{(c_k-c_{k-j})^2}{2}+3j+1},
\vs\end{align}
using {$\bar r\bar x_m^2< \frac{1}{2}$} (which follows from $x_m^Q<0.262\sigma$ for $m\geq 25$).
The case $k=m$ needs separate treatment. Define $\hat B_m=[b_m,c_m+\bar x_m]$. Then, it is easy to verify that \eqref{priors} still holds if we replace $B_m$ with $\hat B_m$. Therefore, the proof in this case follows from an argument similar to above on noting that $\Prob[\theta\in \hat B_{m}|s]\leq \Prob[\theta\in B_{m}|s]$.
$\hfill\blacksquare$

{\label{proof::lemma::p(th|s)withlog}\noindent\textit{Proof of Lemma~\ref{lemma::p(th|s)withlog}.}}
If $s\geq c_k+\bar r\bar x_{k+2}$, then
\begin{align}\vs
&\frac{\Prob[\theta\in B_{k}|s]}{\Prob[\theta\in B_{k+1}|s]}\hhs\leq\hhs \frac{\Prob[\theta\in B_{k}]\phi(m_{k+1}+\delta-c_k-\bar{r}\bar x_{k+2})}{\Prob[\theta\in B_{k+1}]\phi(m_{k+1}+\delta-c_{k+1}-\bar{r}\bar x_{k+2})}\nonumber\\
&\hhs\leq\hhs \frac{\Prob[\theta\in B_{k}]}{\Prob[\theta\in B_{k+1}]}e^{\Delta_{k+1}(\frac{c_k+c_{k+1}}{2}+\bar{r}\bar x_{k+2}-\delta)}\hhs\leq\hhs e^{\Delta_{k+1}(\bar{r}\bar x_{k+2}-\delta)},
\vs\end{align}
where the last inequality follows from the definition of $m_{k+1}$. However, for the case where $s<c_k+\bar r\bar x_{k+2}$, the upper bound on the likelihood $\Prob[s|\theta\in B_k]$ in the first inequality may be less and hence is replaced by 1, which will thereby lead to
\begin{equation}\vs
\frac{\Prob[\theta\in B_{k}|s]}{\Prob[\theta\in B_{k+1}|s]}\leq e^{\Delta_{k+1}(\bar{r}\bar x_{k+2}-\delta)+\frac{\bar r^2\bar x_{k+2}^2}{2}}.
\vs\end{equation}
The other side of the inequality can be proved similarly.

For the case $k=m-1$, the lower bound $\phi(m_{k+1}+\delta-c_{k+1}-\bar{r}\bar x_{k+2})$ for the likelihood $\Prob[s|\theta\in B_{k+1}]$ is not valid anymore. To fix this, as in Lemma~\ref{lemma::p(theta|s)}, we use $\hat B_m=[b_m,c_m+\bar x_m]$ instead of $B_m$ to obtain,
\begin{align}\vs
\frac{\Prob[\theta\in B_{m-1}|s]}{\Prob[\theta\in B_{m}|s]}&\leq\frac{\Prob[\theta\in B_{m-1}|s]}{\Prob[\theta\in \hat B_{m}|s]}\nonumber\\
&\leq \frac{\Prob[\theta\in B_{m}]}{\Prob[\theta\in \hat B_{m}]}e^{\Delta_{m}(\bar{r}\bar x_m-\delta)+\frac{\bar r^2\bar x_m^2}{2}}.
\vs\end{align}
On the other hand, we can show that
$\frac{\Prob[\theta\in B_{m}]}{\Prob[\theta\in \hat B_{m}]}<1.16$, which completes the proof. It is easy to see that the inequality in LHS stays as before for $k=m-1$.~$\hfill\blacksquare$

{\label{proof::lemma::aF(s)}\noindent\textit{Proof of Lemma~\ref{lemma::aF(s)}.}}
As the first step we bound the effect of intervals other than $B_k\cup B_{k+1}$.
Let $\eta=\E[a_L(\theta)|s,\theta\in B_{k}\cup B_{k+1}]$ and $\eta_{j}=\E[a_L(\theta)|s,\theta\in B_{j}]$ for $-m\leq j\leq m$. Using Lemma~\ref{lemma::p(theta|s)}, we can write
\begin{align}\vs
\label{aFtmp1}
&\sum_{j=1}^{k+m} \frac{\Prob[\theta\in B_{k-j}|s]}{\Prob[\theta\in B_{k}|s]}(\eta-\eta_{k-j})\nonumber\\
&\leq\hhs\sum_{j=1}^{k+m} (c_k-c_{k-j}+2\bar x_m+2\bar{r}\bar x_m)e^{-\frac{(c_k-c_{k-j})^2}{2}+3j+1}\nonumber\\
&\leq\hhs \sum_{j=1}^{k+m}(2j\ubar x_1+2\bar x_m+2\bar r\bar x_m)e^{-2j^2\ubar x_1^2+3j+1}\nonumber\\
&\leq\hhs 4e^{-2\ubar x_1^2+4}\bar x_m \sum_{j=1}^{k+m} je^{-2(j^2-1)\ubar x_1^2+3(j-1)}\nonumber\\
&\leq\hhs 4e^{-2\ubar x_1^2+4}\bar x_m \sum_{j=1}^{\infty} e^{-(j-1)^2}\leq\hhs \frac{5.6e^4}{\sigma^{15}}<10^{-10}\bar r^2\ubar x_1^2,
\vs\end{align}
where we have  used the identity $\sum_{j=0}^{\infty}e^{-j^2}\approx1.386$.
Similarly, we can bound the effect of non-neighboring intervals on the variance:
\begin{align}\vs
\sum_{j=1}^{k+m} \frac{\Prob[\theta\in B_{k-j}|s]}{\Prob[\theta\in B_{k}|s]}(\eta-\eta_{k-j})^2\leq \frac{5.6e^4}{\sigma^{14}}.
\vs\end{align}
On the other hand,
\begin{align}\vs
&\sum_{j=1}^{k+m} \frac{\Prob[\theta\in B_{k-j}|s]}{\Prob[\theta\in B_{k}|s]}\Var[a_L(\theta)|s,\theta\in B_{k-j}]\nonumber\\
&\leq\sum_{j=1}^{k+m} \frac{1}{3}e^{-\frac{(c_k-c_{k-j})^2}{2}+3j+1}\leq \frac{0.5e^{4}}{\sigma^{16}}.
\vs\end{align}
Combining the two, we obtain
\begin{align}\vs
\label{aFtmp2}
\sum_{j=-m\atop j\notin\{k,k+1\}}^{m}\hs\hs \Prob[\theta\hhs\in\hhs B_{j}|s]((\eta-\eta_{j})^2\hhs+\hhs\Var[a_L(\theta)|s,\theta\hhs\in\hhs B_{j}])\hhs\leq\hhs10^{-10}\bar r^2\ubar x_1^2.
\vs\end{align}
Therefore, the effect of intervals other than $B_k$ and $B_{k+1}$ on $a_F(s)$ (and its derivative given by $\Var[a_L|s]$) is quite negligible. Now, focusing on these two intervals (i.e., $B_k$ and $B_{k+1}$), we have
\begin{align}\vs
&\E[a_L(\theta)|s,\theta\in B_k\cup B_{k+1}]\nonumber\\
&=p\E[a_L(\theta)|s,\theta\in B_k]+(1-p)\E[a_L(\theta)|s,\theta\in B_{k+1}]\nonumber\\
&\leq p(c_k+\bar{r}\bar x_{k+1})+(1-p)(c_k+\Delta_{k+1}+\bar{r}\bar x_{k+2})\nonumber\\
&\leq c_k+(1-p)\Delta_{k+1}+\bar{r}\bar x_{k+2},
\vs\end{align}
where
$p=\frac{\Prob[\theta\in B_{k}|s]}{\Prob[\theta\in B_{k}\cup B_{k+1}|s]}$.
The proof for the upper bound on $a_F(s)$ now follows from Lemma~\ref{lemma::p(th|s)withlog}. The proof for the lower bound on $a_F(s)$ is similar.
Now, as for the derivative, we first note that $\frac{d}{ds}a_F(s)=\Var[a_L|s]$. Again, focusing on $B_k\cup B_{k+1}$, we can write
\begin{align}\vs
&\Var[a_L|s,\theta\in B_k\cup B_{k+1}]\nonumber\\
&\leq p \Var[a_L|s,\theta\in B_k]+(1-p)\Var[a_L|s,\theta\in B_{k+1}]\nonumber\\
&+p(1-p)(\E[a_L|s,\theta\in B_k]-\E[a_L|s,\theta\in B_{k+1}])^2.
\vs\end{align}
The rest easily follows from Lemma~\ref{lemma::aFlocal} and Lemma~\ref{lemma::p(th|s)withlog}.
$\hfill\blacksquare$

{\label{proof::lemma::aF(s)-tail}\noindent\textit{Proof of Lemma~\ref{lemma::aF(s)-tail}.}}
Exploiting the term $e^{-\delta(c_m-c_{m-r}-2\bar r\sigma)}$ in \eqref{tac::tmp2} (for $k=m$),
it is easy to observe that the same upper bounds given by \eqref{aFtmp1} and \eqref{aFtmp2} hold for the effect of intervals other than $B_m$ on $a_F(s)$ provided
\begin{equation}\vs
e^{-\delta(c_m-c_{m-r}-2\bar r\bar x_m)}(\E[a_L(\theta)|s,\theta\in B_m]-c_m)<2\bar x_m.
\vs\end{equation}
Verifying the above inequality is quite straightforward using Lemma~\ref{lemma::aFlocal-tail}, and specially noting that
$\E[a_L(\theta)|s,\theta\in B_m]-c_m<3r_L\sigma(\delta+1)$ for $\delta>\bar x_{m+1}$.
The proof of the lemma is now an immediate consequence of Lemma~\ref{lemma::aFlocal-tail}.
$\hfill\blacksquare$

{\label{proof::lemma::aL}\noindent\textit{Proof of Lemma~\ref{lemma::aL}.}}
We start by showing that
\begin{equation}\vs
\label{appen::JLck}
\tilde J_L(c_k)=\int_{-\infty}^{\infty}(a_F(s)- c_k)^2\phi(s-c_k)ds\leq 1.1\bar r^2\bar x_{k+2}^2.
\vs\end{equation}
Using the upper bound on $a_F(s)$ given in Corollary~\ref{coraF(s)bounds}, we can write
\begin{align}\vs
\label{appen::JLcktmp1}
&\int_{c_k}^{m_{k+1}}\hs(a_F(s)- c_k)^2\phi(s-c_k)ds\nonumber\\
&\leq\hs \int_{0}^{m_{k+1}-c_k}\hspace{-30pt}(1.01\bar r\bar x_{k+2}+1.17\Delta_{k+1}e^{-\delta\Delta_{k+1}})^2\phi(m_{k+1}-c_k-\delta)d\delta.
\vs\end{align}
A useful inequality here is
\begin{align}\vs
\label{appen::JLcktmp11}
\int_{0}^{\Lambda}\hs e^{-\delta\Delta}\phi(\Lambda-\delta)d\delta\leq\phi(\Lambda)\int_{0}^{\Lambda}\hs e^{-\delta(\Delta-\Lambda)}d\delta\leq\frac{\phi(\Lambda)}{\Delta-\Lambda}.
\vs\end{align}
Another useful property is that
\begin{align}\vs
\label{appen::JLcktmp12}
\frac{\Delta_{k+1}}{2}+\frac{2.3}{\Delta_{k+1}}\geq m_{k+1}-c_k\geq\frac{\Delta_{k+1}}{2}-\frac{1.1}{\Delta_{k+1}}\geq\ubar x_k.
\end{align}
The proof of the RHS is similar to part iii) of Lemma~\ref{lemma::baseconfig}:
\begin{align}\vs
&\frac{\Prob[\theta|\theta\in B_{k+1}]}{\Prob[\theta|\theta\in B_k]}\leq\frac{b_{k+2}-b_{k+1}}{b_{k+1}-b_{k}}\leq\frac{0.5(c_{k+2}-c_{k})+0.2r_L}{0.5(c_{k+1}-c_{k-1})-0.2r_L}\nonumber\\
&\leq\frac{x_{k+2}^{\rm Q}+x_{k+1}^{\rm Q}+3}{x_{k+1}^{\rm Q}+x_{k}^{\rm Q}-3}\leq\frac{e(1+e)x_{k}^{\rm Q}+3}{(1+e)x_{k}^{\rm Q}-3}.
\vs\end{align}
As a result,
\begin{align}\vs
\label{appen::Deltamk}
m_{k+1}-c_k\geq\frac{\Delta_{k+1}}{2}-\ln\left(\frac{e(1+e)x_{k}^{\rm Q}+3}{(1+e)x_{k}^{\rm Q}-3}\right)\frac{1}{\Delta_{k+1}}\geq
\frac{\Delta_{k+1}}{2}-\frac{1.1}{\Delta_{k+1}}\geq\ubar x_k,
\vs\end{align}
on noting that $x_k^{\rm Q}\geq 11.5$ and $\Delta_{k+1}\geq 17$.
This also implies that $\Delta_{k+1}\leq 2\Delta_{k+1}^m+0.14$, where we define $\Delta_{k+1}^m=m_{k+1}-c_k$.  The LHS of \eqref{appen::JLcktmp12} follows from an analysis similar to \eqref{eq::Probratio} for the adjacent intervals of $B_k$ and $B_{k+1}$:
\begin{align}\vs
\frac{\Prob[\theta\in B_{k}]}{\Prob[\theta\in B_{k+1}]}&\leq1.3\times\frac{\phi(\frac{c_{k-1}^{\rm Q}+c_{k}^{\rm Q}+6}{2\sigma})}{\phi(\frac{c_{k}^{\rm Q}+c_{k+1}^{\rm Q}}{2\sigma})}\nonumber\\
&\leq 1.3e^{\frac{1}{2\sigma^2}(2x_m^{\rm Q}-3)(2c_m^{\rm Q}+3)}\leq1.3 e^2,
\vs\end{align}
from which the proof follows from the definition of $m_{k+1}$.

Incorporating \eqref{appen::JLcktmp11} and \eqref{appen::JLcktmp12} into \eqref{appen::JLcktmp1} and after some manipulation, we can find
\begin{align}\vs
\label{appen::JLcktmp2}
&\int_{c_k}^{m_{k+1}}\hs(a_F(s)- c_k)^2\phi(s-c_k)ds\leq \frac{1.01^2\bar r^2\bar x_{k+2}^2}{2}\nonumber\\
&+ \left(\frac{1.37(2\Delta_{k+1}^m+0.14)^2}{2(2\Delta_{k+1}^m+0.14)-\Delta_{k+1}^m}+\frac{2.37\bar r\bar x_{k+2}(2\Delta_{k+1}^m+0.14)}{(2\Delta_{k+1}^m+0.14)-\Delta_{k+1}^m}\right)\phi(\Delta_{k+1}^m)\nonumber\\
&\leq \frac{1.01^2\bar r^2\bar x_{k+2}^2}{2}+ \left(\frac{1.37(2\Delta_{k+1}^m+0.14)}{1.5}+\frac{2.37\bar r\bar x_{k+2}}{0.5}\right)\phi(\Delta_{k+1}^m)\nonumber\\
&\leq \frac{1.01^2\bar r^2\bar x_{k+2}^2}{2}+ \left(\frac{1.37(2\ubar x_1+0.14)}{1.5}+4.74\bar r\bar x_{k+2}\right)\phi(\ubar x_1)\nonumber\\
&\leq \frac{1.01^2\bar r^2\bar x_{k+2}^2}{2}+ (1.83\ubar x_1+0.1+4.74\bar r\bar x_{k+2})\phi(\ubar x_1)\nonumber\\
&\leq \frac{1.01^2\bar r^2\bar x_{k+2}^2}{2}+ \frac{2\ubar x_1e^{-2(\ubar x_1-2)}}{\sqrt{2\pi e^4}\sigma^4},
\vs\end{align}
where we have also used $\ubar x_k\geq \ubar x_1\geq 2\sqrt{2\ln\sigma}+2$.

For $s\in[m_{k+1},c_{k+1}]$, we have $a_F(s)\leq c_{k+1}+1.1\bar r \bar x_{k+2}$ according to Corollary~\ref{cor1}. Therefore,
\begin{align}\vs
\label{appen::JLcktmp3}
&\int_{m_{k+1}}^{c_{k+1}}\hs\hs(a_F(s)\hhs-\hhs c_k)^2\phi(s\hhs-\hhs c_k)ds\hhs\leq
(\Delta_{k+1}+1.1\bar r \bar x_{k+2})^2\Phi(c_k\hhs-\hhs m_{k+1})\nonumber\\
&\leq (2\Delta_{k+1}^m+0.14+1.1\bar r \bar x_{k+2})^2\frac{\phi(\Delta_{k+1}^m)}{\Delta_{k+1}^m}\nonumber\\
&\leq\hhs (2\ubar x_1\hhs+0.14+\hhs1.1\bar r \bar x_{k+2})^2\frac{\phi(\ubar x_1)}{\ubar x_1}
\hhs\leq\hhs \frac{5\ubar x_1e^{-2(\ubar x_1-2)}}{\sqrt{2\pi e^4}\sigma^4}.
\vs\end{align}
Using $a_F(s)\leq s+1.1\bar x_m$, we can show
\begin{align}\vs
\label{appen::JLcktmp4}
&\int_{c_{k+1}}^{\infty}\hs(a_F(s)- c_k)^2\phi(s-c_k)ds\hhs\leq\hhs\frac{\sqrt{2\pi}}{2}\phi(2\ubar x_1)(1.1\bar x_m\hhs+\hhs 2\ubar x_1\hhs+\hhs1)^2\nonumber\\
&\leq\frac{\sqrt{2\pi}}{2}\phi(4\sqrt{2\ln\sigma})(1.1\bar x_m+\ubar x_1+1)^2\leq \frac{(1.1\bar x_m+\ubar x_1+1)^2}{2\sigma^{16}}<10^{-4}\bar r^2\bar x_{k+2}^2.
\vs\end{align}
Combining \eqref{appen::JLcktmp2}, \eqref{appen::JLcktmp3}, and \eqref{appen::JLcktmp4}, we can arrive at
\begin{align}\vs
\int_{-\infty}^{\infty}(a_F(s)- c_k)^2\phi(s-c_k)ds\hhs&\leq 1.01^2\bar r^2\bar x_{k+2}^2\hhs+2\times 10^{-4}\bar r^2\bar x_{k+2}^2\hhs
+\frac{14\ubar x_1e^{-2(\ubar x_1-2)}}{\sqrt{2\pi e^4}\sigma^4}\nonumber\\
&\leq 1.1\bar r^2\bar x_{k+2}^2.
\vs\end{align}
{For the case of $k=m$,} $\int_{c_m}^{\infty}(a_F(s)- c_m)^2\phi(s-c_m)ds$ needs a different treatment. First we note that
\begin{align}\vs
&\int_{c_m+\bar x_{m+1}}^{\infty}\hs\hs(a_F(s)- c_m)^2\phi(s-c_m)ds\hhs\leq\hhs\int_{\bar x_{m+1}}^{\infty}\hs\hs9r_L^2\sigma^2(\delta+1)^2\phi(\delta)d\delta\nonumber\\
&\leq9r_L^2\sigma^2\phi(\bar x_{m+1})\frac{(\bar x_{m+1}+2)^2}{2}< 10^{-4}\bar r^2\bar x_{m+1}^2.
\vs\end{align}
Also,
\begin{equation}\vs
\int_{c_m}^{c_m+\bar x_{m+1}}\hs(a_F(s)- c_m)^2\phi(s-c_m)ds\hhs<\hhs0.5\bar r^2\bar x_{m+1}^2,
\vs\end{equation}
using which it is straightforward to verify that \eqref{appen::JLck} holds for $k=m$ as well (define $\bar x_{m+2}=\bar x_{m+1}$ for consistency).

Let $\theta=c_k+\epsilon$, with $0\leq\epsilon\leq\frac{\Delta_{k+1}}{2}$. We first show that $\tilde a_L(\theta)$ lies in a $5\bar r\bar x_{k+2}$-vicinity of either $c_k$ or $c_{k+1}$.
We begin with the case where $\tilde a_L(\theta)\in[c_k,c_{k+1}]$. Let $\tilde a_L(\theta)=c_k+\epsilon'$, with $5\bar r\bar x_{k+2}\leq\epsilon'\leq{\Delta_{k+1}}-5\bar r\bar x_{k+2}$.
We can use Corollary~\ref{cor1} to obtain a lower bound for $\tilde J_L(\tilde a_L)$:
\begin{align}\vs
\label{JLtildeUB}
\tilde J(\tilde a_L)\hhs\geq\hhs (\epsilon'\hhs-\hhs1.1\bar r\bar x_{k+2})^2(1\hhs-\hhs\Phi(\frac{4\sqrt{2\ln\sigma}}{5}))\hhs\geq\hhs\frac{5(\epsilon'\hhs-1.1\bar r\bar x_{k+2})^2}{16\sqrt{\pi\ln\sigma}\sigma^{\frac{16}{25}}},
\vs\end{align}
where the last inequality follows from the property that ${1-\Phi(x)}\geq\frac{x{\phi(x)}}{1+x^2}$.
Putting this together with $\tilde u_L(\theta,\tilde a_L)\geq \tilde u_L(\theta,c_k)\geq-r_L(\frac{\Delta_{k+1}}{2})^2-1.1(1-r_L)\bar r^2\bar x_{k+2}^2$, it is easy to show that $\epsilon'-1.1\bar r\bar x_{k+2}<\frac{\Delta_{k+1}}{\sqrt{\sigma}}$ for $\sigma\geq300$. A  second use of Corollary~\ref{cor1} now yields $\tilde J_L(\tilde a_L)\geq (\epsilon'-1.1\bar r\bar x_{k+2})^2\Phi(\frac{\Delta_{k+1}}{4})\geq 0.99(\epsilon'-1.1\bar r\bar x_{k+2})^2$ noting $\Delta_{k+1}\geq 2\ubar x_1>17$. Therefore,
\begin{equation}\vs
\label{appen::tmpJLck}
\tilde u_L(\theta,\tilde a_L)\leq -r_L(\epsilon'-\epsilon)^2-0.99(\epsilon'-1.1\bar r\bar x_{k+2})^2.
\vs\end{equation}
On the other hand, using \eqref{appen::JLck}, we get $\tilde u_L(\theta,c_k)\geq -r_L\epsilon^2-0.99(1-r_L)(1.1\bar{r}\bar x_{k+2})^2$. The RHS in \eqref{appen::tmpJLck} is maximized for $\epsilon^*=\frac{r_L\epsilon+0.99(1-r_L)\times1.1\bar r\bar x_{k+2}}{r_L+0.99(1-r_L)}$. Having $\tilde u_L(\theta,\tilde a_L)\geq \tilde u_L(\theta,c_k)$ requires that $|\epsilon'-\epsilon^*|<|0-\epsilon^*|$, that is,
$\epsilon'\leq 2\epsilon^*$ (recall the assumption by contradiction that $5\bar r\bar x_{k+2}\leq\epsilon'$). This then requires $\epsilon'<5\bar r\bar x_{k+2}$.

For the case $\tilde a_L(\theta)\notin[c_k,c_{k+1}]$, suppose $\tilde a_L(\theta)<c_k$ (the other case is similar), and let $\tilde a_L(\theta)=c_k-\epsilon'$. Then, it follows from $\tilde u_L(\theta,\tilde a_L)\geq \tilde u_L(\theta,c_k)$ that
$-r_L\epsilon^2-1.1(1-r_L)\bar r^2\bar x_{k+2}^2\leq -r_L(\epsilon+\epsilon')^2$,
from which it easily follows that $\epsilon'<1.1$. Now an argument similar to the case $\tilde a_L(\theta)\in[c_k,c_{k+1}]$ shows that
$\tilde J_L(\tilde a_L)\geq 0.99(\epsilon'-1.1\bar r\bar x_{k+2})^2$. Combining this with $|\theta-c_k|<|\theta-\tilde a_L|$ and $\tilde u_L(\theta,\tilde a_L)\geq \tilde u_L(\theta,c_k)$, we get $1.1\bar r^2\bar x_{k+2}^2> 0.99(\epsilon'-1.1\bar r\bar x_{k+2})^2$ resulting in $\epsilon'<2.5\bar r\bar x_{k+2}$.

Similar to Lemma 7 in \cite{Witsen_68}, we can show that $\tilde a_L(\theta)$ is increasing. The fact that $\tilde a_L(\theta)$ is increasing implies that it cannot swing between the two neighborhoods.
Therefore, there exists a unique $\tilde b_{k+1}$ separating the two regimes of $[c_k-5\bar r\bar x_{k+2}, c_k+5\bar r\bar x_{k+2}]$ and $[c_{k+1}-5\bar r\bar x_{k+2}, c_{k+1}+5\bar r\bar x_{k+2}]$, thus completing the proof.
$\hfill\blacksquare$

{\label{proof::lemma::daLBR}\noindent\textit{Proof of Lemma~\ref{lemma::daLBR}.}}
We start by
\begin{align}\vs
\label{eq::dJ}
&\frac{d}{d a_L}\int_{-\infty}^{\infty}(a_F(s)- a_L)^2\phi(s-a_L)ds\nonumber\\
&=\frac{d}{d a_L}\int_{-\infty}^{\infty}(a_F(s+a_L)- a_L)^2\phi(s)ds\nonumber\\
&=2\int_{-\infty}^{\infty}( a_L-a_F(s+a_L))( 1-\frac{d}{ds}a_F(s+a_L))\phi(s)ds\nonumber\\
&=2\int_{-\infty}^{\infty}( a_L-a_F(s))( 1-\frac{d}{ds}a_F(s))\phi(s-a_L)ds.
\vs\end{align}
Changing the order of differentiation and integration requires
the map $a_L\mapsto(a_F(s+a_L)- a_L)^2$ to be (i) continuously differentiable, and (ii) its partial derivative be bounded by an integrable function in an open interval around $a_L$.\footnote{See, e.g., Proposition 14.2.2 in \cite{Gasquet_Analysis}.} These are easy to verify, noting that the
follower's strategy $a_F(s)$ (i) is analytic, from which the continuity of partial derivatives is immediate, and (ii) $|a_F(s)|<|s|+\sigma$ and $|\frac{d}{ds}a_F(s)|<s^2+c_m^2$, from which integrability of the partial derivatives follows from the finiteness of the moments of normal distribution.
We can similarly, verify the following identity derived using the integration by part for any given $a_L\in\mathbb{R}$:
\begin{align}\vs
\label{eq::integratebypart}
&\int_{-\infty}^{\infty}(a_L-a_F(s))(s-a_L)\phi(s-a_L)ds=-\int_{-\infty}^{\infty}(a_L-a_F(s))d\phi(s-a_L)\nonumber\\
&=\int_{-\infty}^{\infty}\phi(s-a_L)d(a_L-a_F(s))=-\int_{-\infty}^{\infty}\frac{d}{ds}a_F(s)\phi(s-a_L)ds.
\vs\end{align}

Using \eqref{eq::J}, the first order condition for the optimal leader's response $\tilde a_L$ gives
\begin{equation}\vs
\label{appen::tildeaL}
r_L(\tilde a_L-\theta)+(1-r_L)\hs\int_{-\infty}^{\infty}\hspace{-12pt}(\tilde a_L-a_F(s))(1-\frac{d}{ds}a_F(s))\phi(s-\tilde a_L)ds=0.
\vs\end{equation}
Define the real analytic function $\Theta:\mathbb{R}\to\mathbb{R}$ as
\begin{equation}\vs
\label{eq::Theta}
\Theta(x)= x+\frac{(1-r_L)}{r_L}\hs\int_{-\infty}^{\infty}\hspace{-12pt}(x-a_F(s))(1-\frac{d}{ds}a_F(s))\phi(s-x)ds.
\vs\end{equation}
From \eqref{appen::tildeaL} it then follows that  $\Theta(\cdot)$ is the left inverse of the leader's best response strategy $\tilde a_L(\cdot)$, that is, $\Theta(\tilde a_L(\theta))=\theta$ for all $\theta\in\mathbb{R}$.

Differentiating $\Theta(x)$ we get
\begin{align}\vs
\label{tmp2}
&{r_L}\frac{d}{dx}\Theta(x)=r_L+\hhs(1-r_L)\hs\int_{-\infty}^{\infty}(1-\frac{d}{ds}a_F(s))\phi(s-x)ds\nonumber\\
&+\hhs(1-r_L)\hs\int_{-\infty}^{\infty}\hhs(x-a_F(s))(1-\hhs\frac{d}{ds}a_F(s))(s-\hhs x)\phi(s\hhs-\hhs x)ds\nonumber\\
&=r_L+\hhs(1-r_L)\hs\left(1-\int_{-\infty}^{\infty}\hhs(2+(x-a_F(s))(s\hhs-\hhs x))\hhs\frac{d}{ds}a_F(s)\phi(s\hhs-\hhs x)ds\right),
\vs\end{align}
where the last equality is obtained using \eqref{eq::integratebypart}.
Next, we use \eqref{tmp2} to bound $\frac{d}{dx}\Theta(x)$ for $x=\tilde a_L\in [c_k-5\bar r\bar x_{k+2},c_k+5\bar r\bar x_{k+2}]$.

Using \eqref{eq::daF(s)bounds} to bound $\frac{d}{ds}a_F(s)$ and then applying \eqref{appen::JLcktmp11}, we can obtain
\begin{align}\vs
\label{appen::ddJL1}
&\int_{\tilde a_L}^{m_{k+1}}(\frac{d}{ds}a_F(s)-1.01\bar r^2\bar x_{k+2}^2)\phi(s-\tilde a_L)ds\nonumber\\
&\leq 1.17\Delta_{k+1}^2\int_{0}^{m_{k+1}-\tilde a_L}e^{-\Delta_{k+1}\delta}\phi(m_{k+1}-\tilde a_L-\delta)d\delta\nonumber\\
&\leq\frac{1.17\Delta_{k+1}^2\phi(m_{k+1}-\tilde a_L)}{\Delta_{k+1}-m_{k+1}+\tilde a_L}.\vs\end{align}
Another useful inequality similar to \eqref{appen::JLcktmp11} is
\begin{align}\vs
\label{}
\int_{0}^{\Lambda}\hs e^{-\delta\Delta}\phi(\Lambda+\delta)d\delta\leq\frac{\phi(\Lambda)}{\Delta+\Lambda},
\vs\end{align}
which together with the bound on $\frac{d}{ds}a_F(s)$ given in \eqref{eq::daF(s)bounds} yields
\begin{align}\vs
\label{appen::ddJL2}
&\int_{m_{k+1}}^{c_{k+1}}(\frac{d}{ds}a_F(s)-1.01\bar r^2\bar x_{k+2}^2)\phi(s-\tilde a_L)ds\nonumber\\
&\leq 1.17\Delta_{k+1}^2\int_{0}^{c_{k+1}-m_{k+1}}e^{-\Delta_{k+1}\delta}\phi(m_{k+1}-\tilde a_L+\delta)d\delta\nonumber\\
&\leq\frac{1.17\Delta_{k+1}^2\phi(m_{k+1}-\tilde a_L)}{\Delta_{k+1}+m_{k+1}-\tilde a_L}.\vs\end{align}
Let $\Delta_{k+1}^L=m_{k+1}-\tilde a_L$. Then, using \eqref{appen::Deltamk}, we can show that for $\tilde a_L\in [c_k-5\bar r\bar x_{k+2},c_k+5\bar r\bar x_{k+2}]$, $\Delta_{k+1}^L\geq \ubar x_k$. Also, $\Delta_{k+1}\leq 2\Delta_{k+1}^L+0.14$. This, together with
\eqref{appen::ddJL1} and \eqref{appen::ddJL2}, yields
\begin{align}\vs
\label{eq::deltalk}
&\int_{\tilde a_L}^{c_{k+1}}(\frac{d}{ds}a_F(s)-1.01\bar r^2\bar x_{k+2}^2)\phi(s-\tilde a_L)ds\nonumber\\
&\leq \frac{2.34\Delta_{k+1}^3\phi(m_{k+1}-\tilde a_L)}{\Delta_{k+1}^2-(m_{k+1}-\tilde a_L)^2}\leq\frac{2.34(2\Delta_{k+1}^L+0.14)^3\phi(\Delta_{k+1}^L)}{(2\Delta_{k+1}^L+0.14)^2-(\Delta_{k+1}^L)^2}\nonumber\\
&\leq \frac{2.34(2\ubar x_1+0.14)^3\phi(\ubar x_1)}{(2\ubar x_1+0.14)^2-\ubar x_1^2}\leq {(6.24\ubar x_1+1)\phi(\ubar x_1)}\leq \frac{(6.24\ubar x_1+1)e^{-2(\ubar x_1-2)}}{\sqrt{2\pi e^4}\sigma^4},
\vs\end{align}
where we have again used $\ubar x_1\geq 2\sqrt{2\ln\sigma}+2$.
On the other hand,
\begin{align}\vs
&\int_{c_{k+1}}^{c_m+\bar x_{m+1}}(\frac{d}{ds}a_F(s)-1.01\bar r^2\bar x_{k+2}^2)\phi(s-\tilde a_L)ds\nonumber\\
&\leq\Big(\frac{1}{4}(\Delta_{m+1}+2\bar r\bar x_{m+1})^2+0.01 \bar r^2\bar x_{m+1}^2-1.01\bar r^2\bar x_{k+2}^2\Big) \frac{\phi(c_{k+1}-\tilde a_L)}{c_{k+1}-\tilde a_L}\nonumber\\
&\leq\bar x_{m+1}^2((1+\bar r)^2+0.01 \bar r^2) \frac{\phi(2\ubar x_1)}{2\ubar x_1}\nonumber\\
&\leq\frac{\bar x_{m+1}^2((1+\bar r)^2+0.01 \bar r^2)}{2\sqrt{2\pi}\sigma^{16}\ubar x_1}.
\vs\end{align}
Combining the two, we get
\begin{align}\vs
\label{daLtmp1}
&\int_{c_k}^{c_m+\bar x_{m+1}}(\frac{d}{ds}a_F(s)-1.01\bar r^2 \bar x_{k+2}^2)\phi(s-\tilde a_L)ds\nonumber\\
&<\frac{(6.24\ubar x_1+1)e^{-2(\ubar x_1-2)}}{\sqrt{2\pi e^4}\sigma^4}+\frac{\bar x_{m+1}^2((1+\bar r)^2+0.01 \bar r^2)}{2\sqrt{2\pi}\sigma^{16}\ubar x_1}
<0.01\bar r^2\bar x_{k+2}^2,
\vs\end{align}
noting that $\sigma\geq 300$.

{For $s>c_m+\bar x_{m+1}$,} using Lemma~\ref{lemma::aF(s)-tail} and the fact that $\tilde a_L(\theta)\leq\tilde a_L(c_m)< c_m+5r_L\sigma$, and similar machinery to the above, we get
\begin{align}\vs
\label{k=mdaL1}
&\int_{c_m+\bar x_{m+1}}^{\infty}\frac{d}{ds}a_F(s)\phi(s-\tilde a_L)ds\nonumber\\
&\leq3r_L^2\sigma^2\int_{0}^{\infty}(\delta+\bar x_{m+1})^2\phi(\delta+\bar x_{m+1}-5r_L \bar x_{m+1})d\delta\nonumber\\
&\leq\hhs\frac{3\sqrt{2\pi}}{2}r_L(\bar x_{m+1}\hhs+\hhs1)^2\phi(\bar x_{m+1}\hhs-\hhs5 r_L\bar x_{m+1})\hhs<\hhs 10^{-4}\bar r^2\hhs\bar x_{k+2}^2.
\vs\end{align}
This, along with \eqref{daLtmp1} yields
\begin{equation}\vs
\label{k<mdaL}
\int_{-\infty}^{\infty}\frac{d}{ds}a_F(s)\phi(s-\tilde a_L)ds<1.05\bar r^2\bar x_{k+2}^2.
\vs\end{equation}
{The case $k=m$} is even easier on noting that $\frac{d}{ds}a_F(s)\leq 0.8\bar r^2\bar x_{m+1}^2$ over the whole interval $s\in[c_m,c_m+\bar x_{m+1}]$.
%

Now, to bound the other term assume $a_F(\tilde a_L)\leq \tilde a_L$ (the other case is similar).
From \eqref{lemma::aF(s)}, we can see
\begin{align}
a_F(m_{k+1})\geq c_k+\frac{\Delta_{k+1}}{2.17}-1.01\bar r\bar x_{k+2}>c_k+5\bar r\bar x_{k+2}>\tilde a_L.
\end{align}
Using this, together with $\tilde a_L-a_F(\tilde a_L)\leq 5\bar r\bar x_{k+2}+1.1\bar r\bar x_{k+2}<6.1\bar r\bar x_{k+2}$, we get
\begin{align}\vs
\label{skcktmp1}
&\int_{-\infty}^{\infty}\frac{d}{ds}a_F(s)(\tilde a_L-a_F(s))(s-\tilde a_L)\phi(s-\tilde a_L)ds\nonumber\\
&\leq
(\tilde a_L-a_F(\tilde a_L))\hhs\int_{\tilde a_L}^{m_{k+1}}\hs\frac{d}{ds}a_F(s)(s-\tilde a_L)\phi(s-\tilde a_L)ds\nonumber\\
&\leq
6.1\bar r\bar x_{k+2}\left(\frac{1.01\bar r^2\bar x_{k+2}^2}{\sqrt{2\pi}}+(m_{k+1}-\tilde a_L)\hhs\int_{\tilde a_L}^{m_{k+1}}\hs(\frac{d}{ds}a_F(s)-1.01\bar r^2\bar x_{k+2}^2)\phi(s-\tilde a_L)ds\right)\nonumber\\
&\leq
6.1\bar r\bar x_{k+2}\left(\frac{1.01\bar r^2\bar x_{k+2}^2}{\sqrt{2\pi}}+\frac{1.17\Delta_{k+1}^2(m_{k+1}-\tilde a_L)\phi(m_{k+1}-\tilde a_L)}{\Delta_{k+1}-m_{k+1}+\tilde a_L}\right)\nonumber\\
&\leq6.1\bar r\bar x_{k+2}\left(\frac{1.01\bar r^2\bar x_{k+2}^2}{\sqrt{2\pi}}+{1.17(2\ubar x_1+0.14)^2\phi(\ubar x_1)}\right)\nonumber\\
&<2.46\bar r^3\bar x_{k+2}^3+\frac{7.15\bar r\bar x_{k+2}(2\ubar x_1+0.14)^2e^{-2(\ubar x_1-2)}}{\sqrt{2\pi e^4}\sigma^4}<0.1\bar r^2\bar x_{k+2}^2,
\vs\end{align}
for $\sigma\geq 300$.
{For the case $k=m$} and $a_F(\tilde a_L)\leq \tilde a_L$,
\begin{align}\vs
&\int_{-\infty}^{\infty}\hspace{-3pt}\frac{d}{ds}a_F(s)(\tilde a_L-a_F(s))(s-\tilde a_L)\phi(s-\tilde a_L)ds\nonumber\\
&\leq\hspace{-3pt}
(\tilde a_L-a_F(\tilde a_L))\hspace{-3pt}\int_{\tilde a_L}^{\infty}\hspace{-5pt}\frac{d}{ds}a_F(s)(s-\tilde a_L)\phi(s-\tilde a_L)ds.
\vs\end{align}
We break this into two parts. First, using an approach similar to \eqref{k=mdaL1}, we can obtain
\begin{align}\vs
\label{skcktmp2}
\int_{c_m+\bar x_{m+1}}^{\infty}\frac{d}{ds}a_F(s)(s-\tilde a_L)\phi(s-\tilde a_L)ds< 10^{-4}\bar r^2\bar x_{m+1}^2.
\vs\end{align}
The other part is $\int_{\tilde a_L}^{c_m+\bar x_{m+1}}\hs\frac{d}{ds}a_F(s)(s-\tilde a_L)\phi(s-\tilde a_L)ds$. If $\tilde a_L\geq c_m$, then
\begin{align}\vs
\label{skcktmp2b}
\int_{\tilde a_L}^{c_m+\bar x_{m+1}}\hs\frac{d}{ds}a_F(s)(s-\tilde a_L)\phi(s-\tilde a_L)ds\hhs&\leq\hhs\frac{0.8\bar r^2\bar x_{m+1}^2}{\sqrt{2\pi}},
\vs\end{align}
where we have used that $\frac{d}{ds}a_F(s)\leq0.8\bar r^2\bar x_{m+1}^2$ for $s\in [c_m,c_m+\bar x_{m+1}]$ according to Lemma~\ref{lemma::aF(s)-tail}.
Hence we get
\begin{align}\vs
&\int_{-\infty}^{\infty}\frac{d}{ds}a_F(s)(\tilde a_L-a_F(s))(s-\tilde a_L)\phi(s-\tilde a_L)ds\nonumber\\
&\leq
6.1\bar r\bar x_{m+1}(0.32\bar r^2\bar x_{m+1}^2+10^{-4}\bar r^2\bar x_{m+1}^2)<0.1\bar r^2\bar x_{m+1}^2.
\vs\end{align}
We need to revise \eqref{skcktmp2b} when $\tilde a_L<c_m$. On one hand $\tilde a_L-a_F(\tilde a_L)<c_m-a_F(\tilde a_L)<1.1\bar r\bar x_{m+1}$ (instead of $6.1\bar r\bar x_{m+1}$). We also use a loose bound such as\footnote{Note that $\tilde a_L>a_F(\tilde a_L)>c_m-1.01\bar r\bar x_{m+1}$.}
\begin{align}\vs
\label{skcktmp2c}
\int_{\tilde a_L}^{c_m+\bar x_{m+1}}\hs\frac{d}{ds}a_F(s)(s-\tilde a_L)\phi(s-\tilde a_L)ds\hhs&\leq\hhs\frac{1.1\bar r^2\bar x_{m+1}^2}{\sqrt{2\pi}},
\vs\end{align}
in order to get
\begin{align}\vs
&\int_{-\infty}^{\infty}\frac{d}{ds}a_F(s)(\tilde a_L-a_F(s))(s-\tilde a_L)\phi(s-\tilde a_L)ds\nonumber\\
&\leq
1.1\bar r\bar x_{m+1}(\frac{1.1\bar r^2\bar x_{m+1}^2}{\sqrt{2\pi}}+10^{-4}\bar r^2\bar x_{m+1}^2)<0.1\bar r^2\bar x_{m+1}^2.
\vs\end{align}

Finally, putting all together we get
\begin{align}\vs
\label{tac::tmp11}
&\int_{-\infty}^{\infty}\frac{d}{ds}a_F(s)(2+(\tilde a_L-a_F(s))(s-\tilde a_L))\phi(s-\tilde a_L)ds\nonumber\\
&< (2\times1.05+0.1)\bar r^2\bar x_{m+1}^2<2.2e\bar r^2\bar x_m^2<0.45\bar r^2\sigma^2,
\vs\end{align}
where we use that $\bar x_m<0.272\sigma$ for $m\geq 25$ and $\sigma\geq 300$.

Now for the other side, on noting that $a_F(s)\leq s+1.1\bar x_{m+1}$ and that for $s\in[c_{k+1},c_m+\bar x_{m+1}]$ we have $\frac{d}{ds}a_F(s)\leq((1+\bar r)^2+1.01\bar r^2)\bar x_{m+1}^2$, we can write
\begin{align}\vs
\label{daLtmp3}
&\int_{c_{k+1}}^{c_m+\bar x_{m+1}}\frac{d}{ds}a_F(s)((a_F(s)-\tilde a_L)(s-\tilde a_L)-2)\phi(s-\tilde a_L)ds\nonumber\\
&<\hs\hhs((1+\bar r)^2+1.01\bar r^2)\bar x_{m+1}^2\int_{c_{k+1}}^{\infty}(s+1.1\bar x_{m+1}-\tilde a_L)(s-\tilde a_L)\phi(s-\tilde a_L)ds\nonumber\\
&<\hs\hhs\sqrt{2\pi}((1+\bar r)^2+1.01\bar r^2)\bar x_{m+1}^2\phi(c_{k+1}-\tilde a_L)\int_{0}^{\infty}(\delta+c_{k+1}-\tilde a_L+1.1\bar x_{m+1})(\delta+c_{k+1}-\tilde a_L)\phi(\delta)d\delta\nonumber\\
&<\hs\hhs\sqrt{2\pi}((1\hhs+\hhs\bar r)^2\hhs+\hhs1.01\bar r^2)\bar x_{m+1}^2\phi(c_{k+1}-\tilde a_L)\hs\hhs\left(\frac{1\hhs+\hhs(c_{k+1}\hhs-\hhs\tilde a_L)(c_{k+1}\hhs-\hhs\tilde a_L\hhs+\hhs1.1\bar x_{m+1})}{2}\hhs+\hhs\frac{2(c_{k+1}\hhs-\hhs\tilde a_L)\hhs+\hhs1.1\bar x_{m+1}}{\sqrt{2\pi}}\right)\nonumber\\
&<\hs\hhs\sqrt{2\pi}((1+\bar r)^2+1.01\bar r^2)\bar x_{m+1}^2\phi(2\ubar x_1)\left(\frac{1+2\ubar x_1(2\ubar x_1+1.1\bar x_{m+1})}{2}+\frac{4\ubar x_1+1.1\bar x_{m+1}}{\sqrt{2\pi}}\right)\nonumber\\
&<\hs\hhs10\bar x_{m+1}^4\phi(4\sqrt{2\ln\sigma})<10^{-4}\bar r^2\bar x_{m+1}^2.
\vs\end{align}

For $s\geq c_m+\bar x_{m+1}$, $a_F(s)\leq c_m+3r_L\sigma(s-c_m+1)<s$, and $\frac{d}{ds}a_F(s)\leq3r_L(s-c_m)^2$. An approach similar to \eqref{k=mdaL1} leads to
\begin{align}\vs
\label{daLtmp2}
&\int_{c_m+\bar x_{m+1}}^{\infty}\hs\frac{d}{ds}a_F(s)((a_F(s)-\tilde a_L\hhs)(s-\tilde a_L\hhs)-2)\phi(s-\tilde a_L)ds<\hhs 10^{-4}\bar r^2\bar x_{m+1}^2.
\vs\end{align}

If  $c_k\leq s\leq m_{k+1}-\frac{{2\ln\sigma}}{\Delta_{k+1}}$, then, according to Corollary~\ref{coraF(s)bounds},
\begin{align}
&a_F(s)<c_k+1.17\frac{\Delta_{k+1}}{\sigma^2}+1.01 \bar r \bar x_{k+2}<c_k+2.34r_L\Delta_{k+1}+1.01 \bar r \bar x_{k+2}<c_k+5.8\bar r \bar x_{k+2}.
\end{align}
Therefore, $(a_F(s)-\tilde a_L)(s-\tilde a_L)<11\bar r\bar x_{k+2}^2<2$. On the other hand, it is easy to find
\begin{align}\vs
\label{daLtmp11}
&\int_{m_{k+1}-\frac{{2\ln\sigma}}{\Delta_{k+1}}}^{c_{k+1}}\hhs\frac{d}{ds}a_F(s)((a_F(s)\hhs-\hhs\tilde a_L)(s-\hhs\tilde a_L)\hhs-\hhs2)\phi(s\hhs-\hhs\tilde a_L)ds\nonumber\\
&\leq\left((\frac{\Delta_{k+1}}{2}+\bar r\bar x_{k+2})^2+0.01\bar r^2 \bar x_{k+2}^2\right)\Delta_{k+1}
\phi(m_{k+1}-\tilde a_L-\frac{{2\ln\sigma}}{\Delta_{k+1}}).
\vs\end{align}
Let $\Delta_{k+1}^L=m_{k+1}-\tilde a_L$. Then, from the above inequality (and similar to \eqref{eq::deltalk}) we can get
\begin{align}\vs
\label{daLtmp12}
&\int_{m_{k+1}-\frac{{2\ln\sigma}}{\Delta_{k+1}}}^{c_{k+1}}\hhs\frac{d}{ds}a_F(s)((a_F(s)\hhs-\hhs\tilde a_L)(s-\hhs\tilde a_L)\hhs-\hhs2)\phi(s\hhs-\hhs\tilde a_L)ds\nonumber\\
&\leq((\Delta_{k+1}^L+0.07+\bar r\bar x_{k+2})^2+0.01\bar r^2 \bar x_{k+2}^2)(2\Delta_{k+1}^L+0.14)
\phi(\Delta_{k+1}^L-\frac{{2\ln\sigma}}{\Delta_{k+1}})\nonumber\\
&\leq((\ubar x_1+0.07+\bar r\bar x_{k+2})^2+0.01\bar r^2 \bar x_{k+2}^2)(2\ubar x_1+0.14)
\phi(\ubar x_1-\frac{{2\ln\sigma}}{4\sqrt{2\ln\sigma}})\nonumber\\
&\leq\hhs2(\ubar x_1+0.1)^3\phi(\ubar x_1-\hhs\frac{2\sqrt{2\ln\sigma}}{8})\leq\hhs2(2\sqrt{2\ln\sigma}+2.1)^3\phi(2\sqrt{2\ln\sigma}+2-\hhs\frac{2\sqrt{2\ln\sigma}}{8})\nonumber\\
&\leq\hhs\frac{54\ln\sigma\sqrt{\ln\sigma}e^{-\frac{14}{8}(\ubar x_1-2)}}{\sqrt{2\pi e^4}\sigma^3}<\hhs0.15\bar r^2\sigma^2,
\vs\end{align}
for $\sigma\geq300$. Overall, we arrive at
\begin{equation}\vs
\label{aLtmp2}
\int_{-\infty}^{+\infty}\hs\frac{d}{ds}a_F(s)((a_F(s)-\tilde a_L)(s-\tilde a_L)-2)\phi(s-\tilde a_L)ds\hhs<\hhs 0.4\bar r^2\sigma^2.
\vs\end{equation}
Similar argument holds for $k=m$.

It follows from the above analysis that
\begin{align}\vs
\label{eq::dtheta(s)}
&r_L+\hhs(1-r_L)(1-0.45\bar r^2\sigma^2)\leq {r_L}\frac{d}{dx}\Theta(x)\leq r_L+\hhs(1-r_L)(1+0.4\bar r^2\sigma^2),
\vs\end{align}
for $x\in [c_k-5\bar r\bar x_{k+2},c_k+5\bar r\bar x_{k+2}]$, where $\Theta(x)$ defined in \eqref{eq::Theta} is the real analytic left-inverse of the leader's best response strategy $\tilde a_L(\theta)$, that is, $\Theta(\tilde a_L(\theta))=\theta$ for all $\theta\in\mathbb{R}$. Moreover, $\tilde a_L(\theta)$ is an increasing function, as already discussed in the proof of the previous lemma.

We next claim that $\tilde a_L(\theta)$ is continuous over both $(c_k,\tilde b_{k+1})$ and $(\tilde b_{k+1}, c_{k+1})$. Consider $\theta\in(c_k,\tilde b_{k+1})$. From Lemma~\ref{lemma::aL}, we know that
$\tilde a_L(\theta)\in [c_k-5\bar r\bar x_{k+2},c_k+5\bar r\bar x_{k+2}]$. This is a tiny interval around a fixed point of the leader's original strategy $a_L(\theta)$ (according to which the follower's best response $a_F(s)$ is derived). We next show that leader's utility $\tilde u_L(\theta,a_L)$ given in \eqref{eq::uLBR} is strongly concave (in $a_L$) over this interval.
The good news is that we have already bounded $\frac{\partial^2}{\partial a_L^2}\tilde u_L(\theta,a_L)$
for $ a_L\in [c_k-5\bar r\bar x_{k+2},c_k+5\bar r\bar x_{k+2}]$ while bounding the derivative of $\Theta(x)$ over this interval. In fact, from the definition of both functions given in \eqref{eq::uLBR} and \eqref{eq::Theta}, it is easy to verify that
\begin{align}
\frac{\partial^2}{\partial a_L^2}\tilde u_L(\theta,a_L)=-2{r_L}\frac{d}{dx}\Theta(x)_{|x=a_L}.
\end{align}
Strong concavity of $\tilde u_L(\theta,a_L)$ for $a_L\in [c_k-5\bar r\bar x_{k+2},c_k+5\bar r\bar x_{k+2}]$ then follows from \eqref{eq::dtheta(s)}. Recall that for  $\theta\in(c_k,\tilde b_{k+1})$, leader's best response
$\tilde a_L(\theta)\in[c_k-5\bar r\bar x_{k+2},c_k+5\bar r\bar x_{k+2}]$, meaning that $\tilde a_L(\theta)$ is the unique maximizer of $\tilde u_L(\theta,a_L)$ over $[c_k-5\bar r\bar x_{k+2},c_k+5\bar r\bar x_{k+2}]$. This uniqueness property implies the continuity of $\tilde a_L(\theta)$ over $\theta\in(c_k,\tilde b_{k+1})$, noting that both left and right limits of $\tilde a_L(\theta)$ are maximizers of leader's utility over $[c_k-5\bar r\bar x_{k+2},c_k+5\bar r\bar x_{k+2}]$, and hence have to be identical.

Therefore, leader's best response $\tilde a_L(\theta)$ is continuous over both $(c_k,\tilde b_{k+1})$ and $(\tilde b_{k+1}, c_{k+1})$, with its graph coinciding its analytic left inverse $\Theta(x)$ over both intervals. We can thus use \eqref{eq::dtheta(s)} to bound its derivative over each of these intervals, hence completing the proof.
$\hfill\blacksquare$

{\label{proof::lemma::J}\noindent\textit{Proof of Lemma~\ref{lemma::J}.}}
We have already calculated the second derivative of $\tilde J_L$ when deriving the partial derivative in \eqref{tmp2}.
The same argument implies $\frac{d^2}{da_L^2}\tilde J_L(a_L)\geq 2(1-0.45\bar r^2\sigma^2)$ for $a_L\in[c_k-5\bar r\bar x_{k+2},c_k+5\bar r\bar x_{k+2}]$. This, implies that $\tilde J_L$ is strongly convex over $[c_k-5\bar r\bar x_{k+2},c_k+5\bar r\bar x_{k+2}]$. It's unique minimizer $\tilde c_k$ minimizes both losses in the leader's payoff, hence it is the fixed point of $\tilde a_L(\theta)$, that is $\tilde a_L(\tilde c_k)=\tilde c_k$.
$\hfill\blacksquare$

{\label{proof::lemma::aLtail}\noindent\textit{Proof of Lemma~\ref{lemma::aLtail}.}}
The case $\tilde b_m<\theta\leq c_m$ follows from Lemma~\ref{lemma::daLBR}, so we only need to consider $\theta>c_m$.
As the first step, we derive some useful lower bounds on $\tilde J_L(a_L)$ for $a_L=c_m+\epsilon$ with $\epsilon\geq0$. In particular, we claim that for $\epsilon>\bar r\bar x_{m+1}$ and $\sigma\geq300$, we have
\begin{align}\vs
\label{JLtildelb}
\tilde J_L(a_L)\geq 0.99(\epsilon-\bar r\bar x_{m+1})^2.
\vs\end{align}
We consider two cases: If $\epsilon\leq\frac{3}{4}\bar x_{m+1}$, then $a_F(s)-c_m\leq \bar r\bar x_{m+1}$ for $s\leq a_L+\frac{1}{4}\bar x_{m+1}$, and hence $\tilde J_L(a_L)\geq (\epsilon-\bar r\bar x_{m+1})^2\Phi(\frac{\bar x_{m+1}}{4})$. If $\epsilon>\frac{3}{4}\bar x_{m+1}$, then
$a_L-a_F(a_L+\frac{\bar x_{m+1}}{4}-1)\geq (1-4r_L\sigma)\epsilon$, and hence $\tilde J_L(a_L)\geq (1-4r_L\sigma)^2\Phi(\frac{\bar x_{m+1}}{4}-1)\epsilon^2$. These two observations result in the lower bound in \eqref{JLtildelb} noting $\sigma\geq300$.

Consider now $\theta=c_m+\delta$, $0\leq\delta$. We claim that $\tilde a_L(\theta)< c_m+2.2r_L(\delta+\bar x_{m+1})$. Let $\tilde a_L=c_m+\epsilon$. Then,
\begin{equation}\vs
\label{appen::tmpJLaL}
\tilde u_L(\theta,\tilde a_L)\leq -r_L(\epsilon-\delta)^2-0.99(1-r_L)(\epsilon-\bar r\bar x_{m+1})^2.
\vs\end{equation}
Using an approach similar to Lemma~\ref{lemma::aL}, we can show that $\epsilon<2\epsilon^*$ where
\begin{equation}\vs
\epsilon^*=\frac{r_L\delta+0.99(1-r_L) \bar r\bar x_{m+1}}{r_L+0.99(1-r_L)},
\vs\end{equation}
from which we easily get $\epsilon<2.02r_L(\delta+\bar x_{m+1})$. Using this, the proof of the second part of the lemma is quite straightforward.

As for $\frac{d}{d\theta}\tilde a_L(\theta)$, the case $\tilde a_L(\theta)\leq  c_m+5\bar r\bar x_{m+1}$ is covered in Lemma~\ref{lemma::daLBR}. Hence, we study the case $\tilde a_L(\theta)>c_m+5\bar r\bar x_{m+1}$.
Noting that
\begin{align}
\label{eq::tailaLbound}
\tilde a_L(\theta)-c_m\leq2.2r_L(\sigma\bar x_{m+1}+5\bar r\bar x_{m+1}+\bar x_{m+1})<2.4,
\end{align}
and similar to \eqref{k=mdaL1}, we can obtain
\begin{equation}\vs
\int_{c_m+\bar x_{m+1}}^{\infty}\frac{d}{ds}a_F(s)\phi(s-\tilde a_L)ds< 10^{-4}\bar r^2\bar x_{m+1}^2.
\vs\end{equation}
Also,
\begin{equation}\vs
\int_{c_m}^{c_m+\bar x_{m+1}}\frac{d}{ds}a_F(s)\phi(s-\tilde a_L)ds\leq0.8\bar r^2\bar x_{m+1}^2\Phi(\tilde a_L-c_m).
\vs\end{equation}
On the other hand, similar to \eqref{daLtmp1}, we have
\begin{align}\vs
&\int_{-\infty}^{c_m}\frac{d}{ds}a_F(s)\phi(s-\tilde a_L)ds<10^{-4}\bar r^2\bar x_{m+1}^2+1.01\bar r^2\bar x_{m+1}^2\Phi(c_m-\tilde a_L).
\vs\end{align}
Therefore, 
\begin{equation}\vs
\int_{-\infty}^{\infty}\frac{d}{ds}a_F(s)\phi(s-\tilde a_L)ds<0.91\bar r^2\bar x_{m+1}^2<2.48\bar r^2\bar x_m^2,
\vs\end{equation}
since $\bar x_{m+1}=\sqrt{e}\bar x_m$.

To bound the other term, using $\tilde a_L-a_F(s)<2.2r_L(\sigma\bar x_{m+1}+\bar x_{m+1}+5r_L\sigma)+1.5\bar r\bar x_{m+1}<2.5$ for $s\geq c_m+\bar x_{m+1}$ and similar to \eqref{k=mdaL1}, we can show that
\begin{equation}
\int_{c_m+\bar x_{m+1}}^{\infty}\hs\frac{d}{ds}a_F(s)(\tilde a_L-a_F(s))(s\hhs-\hhs\tilde a_L)\phi(s\hhs-\hhs\tilde a_L)ds\hhs<\hhs 10^{-4}\bar r^2\bar x_{m+1}^2.
\end{equation}
For $s\leq c_m+\bar x_{m+1}$, we have $a_F(s)<\tilde a_L$. Therefore,
\begin{align}
&\int_{-\infty}^{c_m+\bar x_{m+1}}\frac{d}{ds}a_F(s)(\tilde a_L-a_F(s))(s-\tilde a_L)\phi(s-\tilde a_L)ds&\nonumber\\
&\leq
\int_{\tilde a_L}^{c_m+\bar x_{m+1}}\frac{d}{ds}a_F(s)(\tilde a_L-a_F(s))(s- \tilde a_L)\phi(s-\tilde a_L)ds\nonumber\\
&\leq\frac{1}{\sqrt{2\pi}}(\bar r\sigma\bar x_{m+1}+1.01\bar r\bar x_{m+1})\times0.8\bar r^2(\frac{\bar x_m+\bar x_{m+1}}{2})^2<0.26\bar r^2\bar x_m^2.
\end{align}
Putting all together, we get
\begin{align}
\label{tmptaildaL2}
&\int_{-\infty}^{\infty}\hhs\frac{d}{ds}a_F(s)(2+(\tilde a_L-a_F(s))(s-\tilde a_L))\phi(s-\tilde a_L)ds<5.22\bar r^2\bar x_{m}^2<0.4\bar r^2\sigma^2,
\end{align}
since $\bar x_m<0.272\sigma$ for $m\geq 25$ and $\sigma\geq300$.

For the other side, similar to \eqref{daLtmp3}, \eqref{daLtmp2}, and \eqref{daLtmp11}, we have
\begin{align}\vs
&\int_{-\infty}^{c_{m-1}}\hs\frac{d}{ds}a_F(s)((a_F(s)\hhs-\hhs\tilde a_L)(s\hhs-\hhs\tilde a_L)-2)\phi(s-\tilde a_L)ds<10^{-4}\bar r^2\bar x_{m+1}^2,\nonumber\\
&\int_{c_m+\bar x_{m+1}}^{\infty}\hhs\frac{d}{ds}a_F(s)((a_F(s)-\tilde a_L)(s-\tilde a_L)-2)\phi(s-\tilde a_L)ds<10^{-4}\bar r^2\bar x_{m+1}^2,\nonumber\\
&\int_{c_{m-1}}^{m_{m}+{\sqrt{\ln\sigma}}}\hs\frac{d}{ds}a_F(s)((a_F(s)-\tilde a_L)(s-\tilde a_L)-2)\phi(s-\tilde a_L)ds<0.1\bar r^2\sigma^2.
\end{align}
Also, same as \eqref{daLtmp3},
\begin{align}\vs
&\int_{-\infty}^{c_{m-1}}\hs\frac{d}{ds}a_F(s)((a_F(s)\hhs-\hhs\tilde a_L)(s\hhs-\hhs\tilde a_L)-2)\phi(s-\tilde a_L)ds<10^{-4}\bar r^2\bar x_{m+1}^2.
\vs\end{align}
Consider now $s\in[m_{m}+\sqrt{\ln\sigma},c_m+\bar x_{m+1}]$. It is easy to verify that
$(a_F(s)-\tilde a_L)(s-\tilde a_L)<2$ for $s\in[\tilde a_L,c_m+\bar x_{m+1}]$. On the other hand, using Lemma~\ref{lemma::aF(s)} and \ref{lemma::aF(s)-tail}, we can verify that $\frac{d}{ds}a_F(s)\leq 1.1 \bar r^2\bar x_{m+1}^2$ for $s\in[m_m+\sqrt{\ln\sigma},\tilde a_L]$. Moreover, it follows from Corollary~\ref{cor1} that $a_F(s)\geq c_m-1.1\bar r\bar x_{m+1}$ for $s\in[m_m+\sqrt{\ln\sigma},\tilde a_L]$, implying that $\tilde a_L-a_F(s)<2.4+1.1\bar r\bar x_{m+1}<2.5$ (see \eqref{eq::tailaLbound}).
This yields
\begin{align}
&\int_{m_{m}+\sqrt{\ln\sigma}}^{\tilde a_L}\frac{d}{ds}a_F(s)((a_F(s)-\tilde a_L)(s-\tilde a_L)-2)\phi(s-\tilde a_L)ds\nonumber\\
&<2.5\times1.1\bar r^2\bar x_{m+1}^2\int_{m_{m}+\sqrt{\ln\sigma}}^{\tilde a_L}(\tilde a_L-s)\phi(s-\tilde a_L)ds\nonumber\\
&<\frac{2.5\times1.1\bar r^2\bar x_{m+1}^2}{\sqrt{2\pi}}<1.1\bar r^2\bar x_{m+1}^2.
\end{align}
From the above analysis it follows that
\begin{align}
&\int_{-\infty}^{\infty}\hs\frac{d}{ds}a_F(s)((a_F(s)-a_L)(s-a_L)-2)\phi(s-a_L)ds\nonumber\\
&<0.1\bar r^2\sigma^2+1.2\bar r^2\bar x_{m+1}^2<(0.1+1.2\times e\times0.272^2)\bar r^2\sigma^2<0.4\bar r^2\sigma^2.
\end{align}
Using this and \eqref{tmptaildaL2}, and following exact same steps as in the proof of Lemma~\ref{lemma::daLBR}, we can show that $\ubar{r}\leq\frac{d}{d\theta}\tilde a_L(\theta)\leq \bar r$.~$\hfill\blacksquare$

{\label{proof::lemma::endpoints}\noindent\textit{Proof of Lemma~\ref{lemma::endpoints}.}}
We start by finding an upper bound for $\tilde J_L(\tilde c_k)$. We use $\tilde J_L(\tilde c_k)\leq \tilde J_L(s_k)$, where we recall that $a_F(s_k)=s_k$. Noting that the upper bound on the derivative of $a_F(s)$ in \eqref{eq::daF(s)bounds} is increasing for $s_k\leq s\leq m_{k+1}$, we can write
\begin{align}\vs
&\int_{s_k}^{m_{k+1}}(a_F(s)-s_k)^2\phi(s-s_k)ds=\int_{s_k}^{m_{k+1}}(a_F(s)-a_F(s_k))^2\phi(s-s_k)ds\leq\nonumber\\
&\int_{s_k}^{m_{k+1}}(s-s_k)^2
(1.17e^{-\Delta_{k+1}(m_{k+1}-s)}\Delta_{k+1}^2+1.01\bar{r}^2\bar x_{k+2}^2)^2\phi(s-s_k)ds.
\end{align}
Following a similar machinery as the one already used in, e.g., \eqref{eq::deltalk}, we can arrive at
\begin{align}\vs
&\int_{s_k}^{m_{k+1}}(a_F(s)-a_F(s_k))^2\phi(s-s_k)ds\leq\hhs
\frac{1.01^2\bar r^4\bar x_{k+2}^4}{2}\nonumber\\
&+2.4\bar r^2\bar x_{k+2}^2\Delta_{k+1}^2 \frac{(m_{k+1}-s_k)^2\phi(m_{k+1}-s_k)}{\Delta_{k+1}-m_{k+1}+s_k}\nonumber\\
&+1.37\Delta_{k+1}^4 \frac{(m_{k+1}-s_k)^2\phi(m_{k+1}-s_k)}{2\Delta_{k+1}-m_{k+1}+s_k}\nonumber\\
&\leq\hhs\frac{1.01^2\bar r^4\bar x_{k+2}^4}{2}+(9.6\bar r^2\bar x_{k+2}^2\ubar x_1^2(\ubar x_1+0.07)+7.31\ubar x_1^2(\ubar x_1+0.07)^3)\phi(\ubar x_1)\nonumber\\
&<\frac{7.5\ubar x_1^5 e^{-2(\ubar x_1-2)}}{\sqrt{2\pi e^4}\sigma^4},
\vs\end{align}
where the last inequality follows from $\ubar x_1\geq2\sqrt{2\ln\sigma}+2$.
Similarly, using \eqref{eq::daFgeneric} and Lemma~\ref{lemma::aF(s)-tail} we can show that
\begin{equation}\vs
\int_{m_{k+1}}^{\infty}(a_F(s)-s_k)^2\phi(s-s_k)ds<\frac{2\ubar x_1^5e^{-2(\ubar x_1-2)}}{\sqrt{2\pi e^4}\sigma^4}.
\vs\end{equation}
These two yield $\tilde J_L(s_k)<\frac{9.5\ubar x_1^5e^{-2(\ubar x_1-2)}}{\sqrt{2\pi e^4}\sigma^4}<0.2r_L^2$ using $\ubar x_1> 8.75$ for $\sigma\geq300$.
It is easy to verify that the same hold when $k=m$. The analysis in this case is even simpler on noting that $\frac{d}{ds}a_F(s)\leq 0.8\bar r^2\bar x_{m+1}^2$ over the whole interval $s\in[c_m,c_m+\bar x_{m+1}]$.

Applying the Envelope's theorem to \eqref{eq::uLBR}, we get
\begin{equation}\vs
\frac{d}{d\theta}\tilde u_L(\theta,\tilde a_L(\theta))=2r_L(\theta-\tilde a_L(\theta)).
\vs\end{equation}
Integrating this, along the inequality below
\begin{equation}\vs
\ubar{r}\theta+(1-\ubar{r})\tilde c_k\leq \tilde a_L(\theta)\leq \bar{r}\theta+(1-\bar{r})\tilde c_k,
\vs\end{equation}
we get
\begin{equation}\vs
r_L(1-\bar{r})(\theta-\tilde c_k)^2\leq \tilde u_L(\theta,\cdot)-\tilde u_L(\tilde c_k,\cdot)\leq r_L(1-\ubar{r})(\theta-\tilde c_k)^2,
\vs\end{equation}
where we use $\tilde u_L(\theta,\cdot)$ as a short-note for $\tilde u_L(\theta,\tilde a_L(\theta))$.
Now, we note that at the endpoint $\theta=\tilde b_{k+1}$, the above should also hold for $\tilde c_{k+1}$.
that is,
\begin{align}\vs
r_L(1-\bar{r})(\tilde b_{k+1}-\tilde c_k)^2\leq \tilde u_L(\tilde b_{k+1},\cdot)&-\tilde u_L(\tilde c_k,\cdot)\leq r_L(1-\ubar{r})(\tilde b_{k+1}-\tilde c_k)^2\nonumber\\
r_L(1-\bar{r})(\tilde b_{k+1}-\tilde c_{k+1})^2\leq \tilde u_L(\tilde b_{k+1},\cdot)&-\tilde u_L(\tilde c_{k+1},\cdot)\leq r_L(1-\ubar{r})(\tilde b_{k+1}-\tilde c_{k+1})^2.
\vs\end{align}
Using this and noting $0<\tilde u_L(\tilde c_{k},\cdot)<0.2r_L^2$ and $0<\tilde u_L(\tilde c_{k+1},\cdot)<0.2r_L^2$, we can arrive at
\begin{align}\vs
(1\hhs-\hhs r_L)|(\tilde b_{k+1}\hhs-\hhs\tilde c_{k+1})^2\hhs-\hhs(\tilde b_{k+1}\hhs-\hhs\tilde c_{k})^2|\hhs<\hhs0.2r_L\hhs+\hhs\frac{r_L^2}{2}(\tilde c_{k+1}\hhs-\hhs\tilde c_{k})^2,
\vs\end{align}
using which the rest of the proof is straightforward.
$\hfill\blacksquare$

{\label{proof::lemma::ck1}\noindent\textit{Proof of Lemma~\ref{lemma::ck1}.}}
Using Corollary~\ref{cor1}, it is straightforward to show that $|s_k-c_k|<1.1\bar r\bar x_{m+1}$.
Evaluating the derivative of $\tilde J_L(a_L)$ at $a_L=s_k$, we get
\begin{equation}\vs
\frac{d}{da_L}\tilde J_L(s_k)\hhs=\hhs2\hhs\int_{-\infty}^{\infty}\hhs(a_F(s_k)-a_F(s))(1-\frac{d}{ds}a_F(s))\phi(s-s_k)ds,
\vs\end{equation}
where we have also used $s_k=a_F(s_k)$.
We consider the case $s_k\leq\tilde c_k$ yielding $\frac{d}{da_L}\tilde J_L(s_k)\leq0$ (the other case is quite similar). Noting that for $s<s_k$, $\frac{d}{ds}a_F(s)>1$ requires $s<m_{k}+0.5$, we can obtain
\begin{align}\vs
-&\frac{d}{da_L}\tilde J_L(s_k)\leq2\int_{s_k}^{\infty}(a_F(s)-a_F(s_k))\phi(s-s_k)ds\nonumber\\
&+2\int_{-\infty}^{m_k+0.5}\hhs(a_F(s_k)-a_F(s))(1-\frac{d}{ds}a_F(s))\phi(s-s_k)ds.
\vs\end{align}
With a bit of manipulation similar to the ones used in Lemma~\ref{lemma::aL}-\ref{lemma::endpoints}, we can obtain
\begin{align}\vs
&\int_{s_k}^{\infty}\hs(a_F(s)-a_F(s_k))\phi(s-s_k)ds\nonumber\\
&\leq\hhs
\frac{1.01\bar r^2(\bar x_{k}+\bar x_{k+1})^2}{4\sqrt{2\pi}}+\frac{4.68\ubar x_1^2e^{-2(\ubar x_1-2)}}{\sqrt{2\pi e^4}\sigma^4}.
\vs\end{align}

As for the second term, similarly
\begin{align}\vs
&\int_{-\infty}^{m_{k}+0.5}(a_F(s_k)-a_F(s))(1-\frac{d}{ds}a_F(s))\phi(s-s_k)ds\nonumber\\
&<(1+\bar r)^2\bar x_{m+1}^2\phi(2\ubar x_1-1.1\bar r\bar x_{m+1})+(1+\bar r)^2\ubar x_1^2\phi(\ubar x_1-0.5-1.1\bar r\bar x_{m+1}).
\vs\end{align}
By putting the above inequalities together and noting $\ubar x_1\geq2\sqrt{2\ln\sigma}+2$ and $\sigma\geq300$, it is a matter of some machinery to verify that
\begin{equation}\vs
0\leq\hhs-\frac{d}{d a_L}\tilde J_L(s_k)\hhs<\frac{1.01\bar r^2(\bar x_{k}+\bar x_{k+1})^2}{2\sqrt{2\pi}}+0.1\bar r^2\ubar x_1.
\vs\end{equation}
On the other hand, as we showed before $\frac{d^2}{da_L^2}\tilde J_L(s_k)\geq 2(1-0.45\bar r^2\sigma^2)$, meaning that $s_k$ is at most $\frac{1.01\bar r^2(\bar x_{k}+\bar x_{k+1})^2}{4\sqrt{2\pi}(1-0.45\bar r^2\sigma^2)}+\frac{0.1\bar r^2\ubar x_1}{2(1-0.45\bar r^2\sigma^2)}<0.42 r_L^2(\frac{\bar x_k+\bar x_{k+1}}{2})^2+0.08 r_L^2\ubar x_1$ away from the minimizer at which the first derivative is zero. This completes the proof.
$\hfill\blacksquare$

{\label{proof::lemma::sk}\noindent\textit{Proof of Lemma~\ref{lemma::sk}.}}
We start by finding the fixed point of $\E[a_L(\theta)|s,\theta\in B_k]$, that is $\E[a_L(\theta)|\hat s_k,\theta\in B_k]=\hat s_k$. $\hat s_k$ is the solution of the following equation
\begin{equation}\vs
\label{skcktmp3}
\int_{b_k}^{b_{k+1}} (a_L(\theta)-\hat s_k)\phi(a_L(\theta)-\hat s_k)\phi(\frac{\theta}{\sigma})d\theta=0.
\vs\end{equation}
Noting that $\hat s_k$ is close to $c_k$ and in particular $|a_L(\theta)-\hat s_k|<1$ for $\theta\in B_k$, and that $x\phi(x)$ is increasing for $|x|<1$, together with the fact that $a_L(\theta)\leq a_L(b_k)+\bar r(\theta-b_k)$,
we obtain
\begin{equation}\vs
\label{tmp4}
\int_{b_k}^{b_{k+1}}\hs\hs (\bar r(\theta-b_k)-(\hat s_k-a_L\hhs(b_k)))\phi(\bar r(\theta-b_k)-(\hat s_k-a_L(b_k)))\phi(\frac{\theta}{\sigma})d\theta\hhs\geq\hhs0.
\vs\end{equation}
Therefore by finding the solution of
\begin{equation}\vs
\label{ttmp4}
\int_{b_k}^{b_{k+1}} \bar r(\theta-y)\phi(\bar r(\theta- y))\phi(\hhs\frac{\theta}{\sigma}\hhs)d\theta=0,
\vs\end{equation}
we can upper-bound the fixed point $\hat s_k$ as $\hat s_k\leq a_L(b_k)+\bar r(y-b_k)$. Simplifying \eqref{ttmp4} yields
\begin{equation}\vs
\label{eq::y}
y=\E_{\bar\psi_y}[\theta|\theta\in B_k],
\vs\end{equation}
where $\bar\psi_y\sim N(\frac{\bar r^2\sigma^2}{1+\bar r^2\sigma^2}y,\frac{\sigma^2}{1+\bar r^2\sigma^2})$. A quick bound for $y$ can be obtained from $\E_{\bar\psi_y}[\theta|\theta\in B_k]\leq e_k+\frac{\bar r^2\sigma^2}{1+\bar r^2\sigma^2}y$, which yields $y-e_k\leq \bar r^2\sigma^2 e_k<0.1\bar x_k$.
An alternative representation of \eqref{eq::y} is
\begin{align}\vs
&\frac{y}{\sqrt{1+\bar r^2\sigma^2}}=\E_{N(0,\sigma^2)}[\theta|b_k\hhs-\hhs\Delta b_k(y)\hhs\leq\hhs\theta\hhs\leq\hhs b_{k+1}\hhs-\hhs\Delta b_k(y)\hhs+\hhs\frac{\bar r^2\sigma^2(b_{k+1}\hhs-\hhs b_k)}{1+\sqrt{1\hhs+\hhs\bar r^2\sigma^2}}],
\vs\end{align}
where $\Delta b_k(y)=\frac{\bar r^2\sigma^2}{1+\bar r^2\sigma^2+\sqrt{1+\bar r^2\sigma^2}}y+\frac{\bar r^2\sigma^2}{1+\sqrt{1+\bar r^2\sigma^2}}(y-b_k)>0$. A useful property of normal distribution is that
$(\frac{\partial}{\partial a}+\frac{\partial}{\partial b})\E_{N(0,\sigma^2)}[\theta|a\leq\theta\leq b]=1-\Var_{N(0,\sigma^2)}[\theta|a\leq\theta\leq b]\geq1-\frac{(b-a)^2}{(1+\sqrt{3})^2\sigma^2}$. Also, $\frac{\partial}{\partial b}\E_{N(0,\sigma^2)}[\theta|a\leq\theta\leq b]\leq\frac{1}{2}$. Applying these to the above equation we can obtain
\begin{align}\vs
0\leq\frac{y}{\sqrt{1+\bar r^2\sigma^2}}+\Delta b_k(y)-e_k\leq&
\frac{\bar r^2\sigma^2(b_{k+1}-b_k)}{2(1+\sqrt{1+\bar r^2\sigma^2})}+\frac{(b_{k+1}-b_k)^2}{(1+\sqrt{3})^2\sigma^2}\Delta b_k(y),
\vs\end{align}
where $e_k=\E_{N(0,\sigma^2)}[\theta|\theta\in B_k]$.
Simplifying this, along with  $y(b_{k+1}-b_k)<2.2\sigma^2$ (which follows from $c_m^{\rm Q}(c_m^{\rm Q}-b_m^{\rm Q})\leq1$), we can arrive at
\begin{align}\vs
|y-e_k|\leq 0.4{\bar r^2\sigma^2(b_{k+1}-b_k)}.
\vs\end{align}
Recalling $\hat s_k\leq a_L(b_k)+\bar r(y-b_k)$, and that $a_L(b_k)\leq c_k-\ubar r(c_k-b_k)\leq c_k-r_L(c_k-b_k)+(r_L-\ubar r)\bar x_k$, we  can reach at
\begin{align}\vs
\label{skcktmp4}
\hat s_k\leq (1-r_L)c_k+r_Le_k+r_L^2\bar x_k+0.4\bar r^3\sigma^2(\bar x_{k}+\bar x_{k+1}).
\vs\end{align}

Following a similar argument to lower-bound $\hat s_k$, we can show that $|\hat s_k-(1-r_L)c_k-r_Le_k|\leq r_L^2\bar x_k+0.4\bar r^3\sigma^2(\bar x_{k}+\bar x_{k+1})$.

To find the fixed point of $a_F(s)$ in $B_k$ (that is $a_F(s_k)=s_k$), we first note that $\hat s_k$ lies within the interval $[c_k-2\bar r\bar x_{m+1},c_k+2\bar r\bar x_{m+1}]$. Moreover, for $s\in [c_k-2\bar r\bar x_{m+1},c_k+2\bar r\bar x_{m+1}]$, using Lemma~\ref{lemma::aF(s)} we can obtain $\frac{d}{ds}a_F(s)<0.01$. This along with$|a_F(s_k)-s_k|<1.5\bar r\bar x_{m+1}$ implies that for
$s\in [c_k-2\bar r\bar x_{m+1},c_k+2\bar r\bar x_{m+1}]$, we have $a_F(s)\in [c_k-2\bar r\bar x_{m+1},c_k+2\bar r\bar x_{m+1}]$. Therefore,
\begin{equation}\vs
\label{skcktmp5}
|s_k-\hat s_k|\hhs<\hhs\frac{|\hat s_k-a_F(\hat s_k)|}{1-0.01}=\hhs\frac{|\E[a_L(\theta)|\hat s_k,\theta\in B_k]-\E[a_L(\theta)|\hat s_k]|}{1-0.01}.
\vs\end{equation}
Assume $\hat s_k\geq c_k$ (the other case is similar).
We have already shown as part of Lemma~\ref{lemma::aF(s)} that while observing $s\in[c_k,c_{k+1}]$, the effect of intervals other than $B_k\cup B_{k+1}$ on $a_F(s)$ is negligible (as given by \eqref{aFtmp1}). Similarly, and by using Lemma~\ref{lemma::p(th|s)withlog}, we can show that
\begin{align}\vs
&{\Prob[\theta\hhs\in\hhs B_{k+1}|\hat s_k]}(\E[a_L(\theta)|\hat s_k,\theta\hhs\in\hhs B_{k+1}]\hhs-\hhs\E[a_L(\theta)|\hat s_k,\theta\hhs\in\hhs B_k])\leq\hspace{-3pt}1.17 e^{-\Delta_{k+1}\delta}\hspace{-2pt}(\Delta_{k+1}\hhs+\hspace{-2pt}2\bar r\bar x_{m+1}),
\vs\end{align}
where $\delta=m_{k+1}-\hat s_k$. Combining this and \eqref{aFtmp1}, we can arrive at
\begin{equation}\vs
|\E[a_L(\theta)|\hat s_k,\theta\in B_k]-\E[a_L(\theta)|\hat s_k]|<\hhs10^{-4} r_L^2\bar x_1.
\vs\end{equation}
After all, we have 
\begin{align}\vs
\label{skcmtmplast}
|s_k-(1-r_L)c_k-r_Le_k|&\hhs<\hhs1.02r_L^2\bar x_k\hhs+\hhs{0.41 \bar r^3\sigma^2(\bar x_{k}+\bar x_{k+1})}<1.9r_L^2\bar x_{k+1}.
\vs\end{align}

The proof has to be modified for the tail case $k=m$, since in the tail $|a_L(\theta)-\hat s_m|<1$ does not hold for all $\theta\in B_m$. To handle this, we define $\E[a_L(\theta)|\hat s_m,\theta\in \hat B_m]=\hat s_m$, where $\hat B_m=[b_m,b_m+\sigma^2-\sigma]$. It is easy to verify that $\bar r(\sigma^2-\sigma)<1$ and hence $|a_L(\theta)-\hat s_m|<1$ holds in this interval. We also need to find an alternative to the variance inequality $\Var_{N(0,\sigma^2)}[\theta|\theta\i B_k]\leq \frac{(b_{k+1}-b_k)^2}{(1+\sqrt{3})^2\sigma^2}$. Here we use  $\Var_{N(0,\sigma^2)}[\theta|\theta\in B_m]\leq (e_m-b_m)^2$. Also, we can show $\Delta b_m(y)<\bar r e_m$, from which and a bit manipulation we can reach at $\Var_{N(0,\sigma^2)}[\theta|\theta\geq b_m-\Delta b_m(y)]\leq \frac{(e_m-b_m)^2}{(1-\bar r)^2}$.
We can then show that the results derived above also holds for $\hat B_m$. More precisely,
\begin{align}\vs
&|\hat s_m-(1-r_L)c_m-r_L\E_{N(0,\sigma^2)}[\theta|\theta\in  \hat B_m]|<r_L^2\bar x_m+0.4\bar r^3\sigma^2(\bar x_m+\bar x_{m+1}).
\vs\end{align}
Let $\theta_b=b_m+\sigma^2-\sigma$. We can easily observe that
\begin{align}\vs
e_m-\E_{N(0,\sigma^2)}[\theta|\theta\in  \hat B_m]&\leq\Prob[\theta\geq\theta_b|\theta\hhs\in\hhs B_m](\E_{N(0,\sigma^2)}[\theta|\theta\geq\theta_b]\hhs-\hhs\E_{N(0,\sigma^2)}[\theta|\theta\hhs\in\hhs  \hat B_m])\nonumber\\
&\leq\frac{\phi(\frac{\theta_b}{\sigma})}{\phi(\frac{b_m}{\sigma})}(\theta_b+\frac{\sigma^2}{\theta_b}-b_m)< 10^{-4}r_L^3\sigma.
\vs\end{align}
We also need to bound $|\E[a_L(\theta)|\hat s_m,\theta\in \hat B_m]-\E[a_L(\theta)|\hat s_m]|$. We start with $|\E[a_L(\theta)|\hat s_m,\theta\in \hat B_m]-\E[a_L(\theta)|\hat s_m, \theta\in B_m]|$. First, note that
\begin{equation}\vs
a_L(\theta_b)\geq c_m+\ubar r(\sigma^2-2\sigma)>c_m+3\bar r\sigma,
\vs\end{equation}
which along $|\hat s_m-c_m|<\bar r\sigma$ implies that $|a_L(\theta_b)-\hat s_m|>2\bar r\sigma$. Consequently, $|\hat s_m-a_L(\theta')|<|a_L(\theta_b)-\hat s_m|$ for all $\theta'\in\hat B_m$.
Similar to \eqref{eq::taildist}, we can hence derive
\begin{equation}\vs
\Prob[\theta|\hat s_m,\theta\in B_m]\leq \frac{e^{-\frac{(\sigma-1)^2}{2}}\phi(\frac{\theta}{\sigma})}{\sigma(1-\Phi(\frac{\theta_b}{\sigma}))},
\vs\end{equation}
for $\theta\geq\theta_b$.
Using this along with  $\E[a_L(\theta)|\hat s_m,\theta\in \hat B_m]\geq c_m-\bar r\sigma$, we can then obtain
\begin{align}\vs
\E[a_L(\theta)|\hat s_m,\theta\in B_m]-\E[a_L(\theta)|\hat s_m,\theta\in \hat B_m]< 10^{-4}r_L^3\sigma.
\vs\end{align}
Therefore, the same steps as in the non-tail case can be followed in this case as well, resulting in \eqref{skcmtmplast} to  also hold for $k=m$.
$\hfill\blacksquare$

{\label{proof::theorem::invariance}\noindent\textit{Proof of Theorem~\ref{theorem::invariance}.}}
We can find by direct calculation of the optimal quantizer for $m=25$ that $\frac{x_1^{\rm Q}}{\sigma}>0.041$. On the other hand, $\frac{2\sqrt{2\ln\sigma}+5}{\sigma}<0.04$ for $\sigma\geq300$.
This implies that $25\in M(\sigma)$ for $\sigma\geq300$ and hence is nonempty.

We next use \eqref{eq::zeta} to verify that Property~\ref{propty2} is also preserved by the best response, completing the proof of the invariance of $A_L^m$. 
It suffices to show that
\begin{equation}\vs
\label{eq::prop2}
|\hat e_k-c_k^{\rm Q}|+0.42 r_L(\frac{\bar x_k+\bar x_{k+1}}{2})^2+2r_L\bar x_{k+1}\leq2.9.
\vs\end{equation}
To bound $|\hat e_k-c_k^{\rm Q}|$, we first note that
$\frac{\partial}{\partial\hat b_k}\hat e_k+\frac{\partial}{\partial\hat b_{k+1}}\hat e_k=1-\frac{\Var_{N(0,\sigma^2)}[\theta|\theta\in\hat B_k]}{\sigma^2}$. Using the inequality $\Var[\theta|\theta\in\hat B_k]>(\hat e_k-\E[\theta|\hat b_k\leq\theta\leq\hat e_k])(\E[\theta|\hat e_k\leq\theta\leq\hat b_{k+1}]-\hat e_k)$ and Lemma~\ref{lemma::baseconfig} we can show that
\begin{align}\vs
\Var_{N(0,\sigma^2)}[\theta|\theta\in\hat B_k]>\frac{\ubar x_k\ubar x_{k+1}}{2(1+\sqrt{e})},
\vs\end{align}
yielding $|\hat e_k-c_k^{\rm Q}|<2.9(1-\frac{\ubar x_k\ubar x_{k+1}}{2(1+\sqrt{e})\sigma^2})$. Therefore to have \eqref{eq::prop2}, it suffices to have
\begin{equation}\vs
\label{}
\frac{0.42 (\bar x_k+\bar x_{k+1})^2}{4\ubar x_k\ubar x_{k+1}}+\frac{2\bar x_{k+1}}{\ubar x_k\ubar x_{k+1}}\leq\frac{2.9}{2(1+\sqrt{e})},
\vs\end{equation}
which is satisfied for sufficiently large values of $\sigma$ (e.g., $\sigma\geq300$) on noting that $\frac{0.42 (e+1)^2}{4e}\leq\frac{2.9}{2(1+\sqrt{e})}$.
This completes the proof of the invariance of $A_L^m$.

The existence of an equilibrium with $a_L^*(\theta)\in A_L^m$ and $a_F^*(s)=\E_\delta[a_L^*|s]$ follows from an argument similar to the one used in \cite{Witsen_68} for the existence of an optimal solution.
Let $(a_L^n(\theta),a_F^n(s)=\E_\delta[a_L^n|s])$ be a maximizing sequence for the ex-ante expected payoff of the leader, that is,
\begin{align}
\label{eq::equilibrium}
&\lim_{n\to\infty}\E_\theta[u_L(\theta,a_L^n(\theta),a_F^n(s))]=\sup \{\E_\theta[u_L(\theta,a_L(\theta),a_F(s))]|a_L(\theta)\in A_L^m,a_F(s)=\E_\delta[a_L|s]\}.
\end{align}
The first step is to show that this supremum is attained for a pair of strategies
$(a_L^*(\theta),a_F^*(s)=\E_\delta[a_L^*|s])$ with $a_L^*(\theta)\in A_L^m$. Strategies  $a_L(\theta)\in A_L^m$ are  increasing and bounded ($|a_L(\theta)|<|\theta|+1$). Using Lemma~8 in \cite{Witsen_68} (which is a variation of the Helly's selection principle based on the diagonalisation argument), there exists a subsequence $a_L^{n_k}(\theta)$ converging pointwise to a limit strategy $a_L(\theta)$, and so do the interval endpoints and fixed points, that is, $\{b_j^{n_k}\}\to\{b_j\}$ and $\{c_j^{n_k}\}\to\{c_j\}$. Relabel $(a_L^{n_k}(\theta),a_F^{n_k}(s))$ as $(a_L^{n}(\theta),a_F^{n}(s))$. From $a_L^n(\theta)\to a_L(\theta)$, it is easy to see that $a_L(\theta)$ satisfies Property 1-2. Property 3, however, concerns the derivative of $a_L$ and thus cannot be deduced directly from pointwise convergence. We now turn into the follower's best response sequence $a_F^n(s)$, claiming that $a_F^n(s)\to a_F(s)$. The proof is again similar to the approach used in Theorem~1 in \cite{Witsen_68}: First note that
\begin{align}
\label{eq::supremum}
a_F^n(s)=\E_\delta[a_L^n|s]=\frac{\int_{-\infty}^{\infty}{a_L^n(\theta)\phi(s-a_L^n(\theta))\phi(\frac{\theta}{\sigma})d\theta}}{\int_{-\infty}^{\infty}{\phi(s-a_L^n(\theta))\phi(\frac{\theta}{\sigma})d\theta}}.
\end{align}
For every $s\in\mathbb{R}$, functions $\phi(s-x)$ and $x\phi(s-x)$ are continuous and bounded functions of $x$. Therefore, $\phi(s-a_L^n(\theta))\to \phi(s-a_L(\theta))$ and $a_L^n(\theta)\phi(s-a_L^n(\theta))\to a_L(\theta)\phi(s-a_L(\theta))$ pointwise in $\theta$, for all $s$. This proves $a_F^n(s)\to a_F(s)$ using bounded convergence theorem.
We can similarly show that $u_L(\theta,a_L^n(\theta),a_F^n(s))\to u_L(\theta,a_L(\theta),a_F(s))$ pointwise in $\theta$ and subsequently
$\E_\theta[u_L(\theta,a_L^n(\theta),a_F^n(s))]\to \E_\theta[u_L(\theta,a_L(\theta),a_F(s))]$, that is, the supremum   in \eqref{eq::equilibrium} is attained by the pair of the strategies $(a_L(\theta),a_F(s)=\E_\delta[a_L|s])$, though not ruling out the possibility of $a_L(\theta)\notin A_L^m$.

Exploiting the fact that $a_F$ is analytic and that $\frac{d}{ds}a_F(s)=\Var[a_L|s]$ and following an argument similar to above, we can then show that $\frac{d}{ds}a_F^n(s)\to\frac{d}{ds}a_F(s)$ pointwise in $s$.
Therefore, the best response characteristics of the follower to a strategy in $A_L^m$ (as given in Lemma 6-7 and Corollary 1-2) which also involves bounds on the derivative, hold for $a_F(s)$.
These conditions force the best response of the leader to $a_F(s)$, denoted by $\tilde a_L(\theta)$, to lie in $A_L^m$.
On the other hand,
\begin{align}
\label{eq::attain2}
\E_\theta[u_L(\theta,\tilde a_L(\theta),\tilde a_F(s))]&\geq \E_\theta[u_L(\theta,\tilde a_L(\theta),a_F(s))]\geq\E_\theta[u_L(\theta,a_L(\theta),a_F(s))],
\end{align}
where $\tilde a_F(s)$ is the follower's best response to $\tilde a_L(\theta)\in A_L^m$. Therefore, the pair of the strategies $(\tilde a_L(\theta),\tilde a_F(s)=\E_\delta[\tilde a_L|s])$ where $\tilde a_L(\theta)\in A_L^m$ attains the supremum in \eqref{eq::supremum}, that is, it is a maximizer for the expected payoff of the leader over $A_L^m$.

Finally, from \eqref{eq::attain2} we should have
$\E_\theta[u_L(\theta,\tilde a_L(\theta),\tilde a_F(s))]= \E_\theta[u_L(\theta,\tilde a_L(\theta),a_F(s))]$. This implies that
$\tilde a_F(s)=a_F(s)$ almost surely, otherwise
replacing $a_F(s)$ with leader's best
response $\tilde a_F(s)$ would result in a higher expected payoff for the leader. Noting that both $\tilde a_F(s)$ and $a_F(s)$
are analytic, almost surely equivalence implies being identical. This implies that the pair of the strategies $(\tilde a_L(\theta),a_F(s))$ are best responses to each other, and hence correspond to an equilibrium of the game.~$\hfill\blacksquare$

{\label{proof::lemma::localoptima}\noindent\textit{Proof of Lemma~\ref{lemma::localoptima}.}}
We need to prove that there does not exist an infinitesimal variation of $(a_L^*,a_F^*)$, namely $(a_L^\delta,a_F^\delta)$, for which $U(a_L^\delta,a_F^\delta)<U(a_L^*,a_F^*)$. Noting that $U(a_L^\delta,\E[a_L^\delta|s])\leq U(a_L^\delta,a_F^\delta)$, we only need to consider the strategies in which the follower's action is the expected action of the leader given the observation $s$ (i.e., $\E[a_L^\delta|s]$). The idea is to show that for a sufficiently small $\delta_L>0$ and any strategy $a_L^\delta$ with $\|a_L^\delta-a_L^*\|_\infty<\delta_L$, the best response image obtained from $a_L^\delta\to a_F^\delta\to \tilde a_L^\delta$ lies in $A_L^m$. The proof then follows from the fact that $U(\tilde a_L^\delta,\E[\tilde a_L^\delta|s])\leq U(a_L^\delta,\E[a_L^\delta|s])$, and that $(a_L^*,a_F^*)$ is the  minimizer of $U$ over all pair of strategies $(a_L,a_F)$ with $a_L\in A_L^m$ (see the proof of Theorem~\ref{theorem::invariance}).

To prove the inclusion of $\tilde a_L^\delta$ in $A_L^m$, it suffices to show that all the properties for the follower's best response to a strategy in $A_L^m$ given specifically by Lemma~\ref{lemma::aF(s)}-\ref{lemma::aF(s)-tail} and Lemma~\ref{lemma::sk} also hold for $a_F^\delta$, noting that these are all we need to deduce Property~\ref{propty1}-\ref{propty3} for the leader's best response (which define the set $A_L^m$).
What is left is then to show that the properties for $a_F^*$ given by Lemma~\ref{lemma::aF(s)}-\ref{lemma::aF(s)-tail} and Lemma~\ref{lemma::sk}
also hold for $a_F^\delta(s)=\E[a_L^\delta|s]$ for sufficiently small $\delta$.
The proof easily follows from a couple of simple observations. First, it is straightforward to  verify that all the bounds given for $a_F^*$ in the aforementioned lemmas are indeed strict. Therefore, by recasting the corresponding inequalities as continuous functions of $\delta_L$ we can ensure that all of them will still hold for sufficiently small $\delta_L$. We elaborate on this in more details in what follows.

We start by verifying that  Lemma~\ref{lemma::aFlocal}-\ref{lemma::aFlocal-tail} also hold for $a_L^\delta$ for small enough $\delta_L$. In Lemma~\ref{lemma::aFlocal},
\begin{align}\vs
a_L^\delta(\theta)-c_k^*&\leq\delta_L+\bar r(b_{k+1}^*-c_k^*)\leq \delta_L+0.1\bar r r_L+\bar r \frac{c_{k+1}^*-c_k^*}{2}\nonumber\\
&\leq \delta_L+0.1\bar r r_L+\bar r (x_{k+1}^{\rm Q}+2.9)\leq\bar r \bar x_{k+1},
\vs\end{align}
for small enough $\delta_L$, where we recall that $\bar x_{k+1}=x_{k+1}^{\rm Q}+3$. Similarly we can show that $a_L^\delta(\theta)-c_k^*\geq-\bar r \bar x_k$, hence Lemma~\ref{lemma::aFlocal} also holds for $a_L^\delta$.
Next, we study the effect of $\delta_L$ in Lemma~\ref{lemma::aFlocal-tail}. As for \eqref{eq::normalapprox}, using
\begin{align}\vs
\frac{\phi(s-a_L^\delta(\theta'))}{\phi(s-a_L^\delta(\theta))}\geq\frac{\phi(s-a_L^*(\theta'))}{\phi(s-a_L^*(\theta))}e^{-\delta_L^2-2\delta_L\bar x_{m+1}},
\vs\end{align}
the RHS of the inequality will be multiplied by $e^{-\delta_L^2-2\delta_L\bar x_{m+1}}$. As a result, the value of $\xi$ in \eqref{eq::xi} will be multiplied by $e^{\delta_L^2+2\delta_L\bar x_{m+1}}$. \eqref{eq::aFtail1} will then become
\begin{align}\vs
\E[a_L^\delta(\theta)|s,b_m^*\leq\theta\leq\theta_c^*]-c_m^*&\leq\delta_L+ 1.128\times1.025e^{\delta_L^2+2\delta_L\bar x_{m+1}}\bar r\bar x_m<0.75\bar r\bar x_{m+1},
\vs\end{align}
for sufficiently small $\delta_L$. For the bound on variance in \eqref{eq::VaraFtail0}-\eqref{eq::VaraFtail1}, let $\mu^*=\E[a_L^*(\theta)|s,b_m^*\leq\theta\leq \theta_c^*]$. Then,
\begin{align}\vs
\Var&[a_L^\delta(\theta)|s,b_m^*\leq\theta\leq \theta_c^*]\leq\Var[a_L^*(\theta)|s,b_m^*\leq\theta\leq \theta_c^*]\nonumber\\
&+\delta_L^2+2\delta_L\sqrt{\Var[a_L^*(\theta)|s,b_m^*\leq\theta\leq \theta_c^*]}.
\vs\end{align}
Hence, \eqref{eq::VaraFtail1} becomes
\begin{align}\vs
\Var[a_L^\delta(\theta)|s,b_m^*\leq\theta\leq \theta_c^*]\leq&1.025 e^{\delta_L^2+2\delta_L\bar x_{m+1}}\times1.128^2\bar r^2(\bar x_m+\frac{4\bar x_m^2\bar\mu}{\sigma^2})^2\nonumber\\
&+\delta_L^2+2.2\delta_L\bar r\bar x_m< 1.2\bar r^2\bar x_m^2,
\vs\end{align}
for sufficiently small $\delta_L$. As for the modification required in the tail effect, $e^{\frac{\delta^2}{2}}$ in \eqref{eq::taildistafterthc} has to be replaced with $e^{\frac{(\delta+\delta_L)^2}{2}}$, using which we can verify that \eqref{eq::taildist} still holds for small enough $\delta_L$. The rest of the changes are similar.

Lemma~\ref{lemma::p(theta|s)}-\ref{lemma::p(th|s)withlog} are based on Lemma~\ref{lemma::aFlocal}-\ref{lemma::aFlocal-tail}, and Lemma~\ref{lemma::aF(s)}-\ref{lemma::aF(s)-tail} are derived using Lemma~\ref{lemma::aFlocal}-\ref{lemma::p(th|s)withlog}, hence also hold for $a_F^\delta$.
Finally, in Lemma~\ref{lemma::sk} which is about the fixed points of the follower's strategy, noting $a_L^\delta(b_k^*)\leq a_L^*(b_k^*)+\delta_L$, we need to add $\delta_L$ to the RHS of \eqref{skcktmp4}. Using this, we can easily verify that this lemma also holds for $a_F^\delta$. Therefore, all the properties required for the follower's strategy to deduce Property~\ref{propty1}-\ref{propty3} for the leader's best response are satisfied for $a_F^\delta$ for sufficiently small $\delta_L$, indicating that $\tilde a_L^\delta$ lies in $A_L^m$. This completes the proof.~$\hfill\blacksquare$

{\label{proof::lemma::DLF}\noindent\textit{Proof of Lemma~\ref{lemma::DLF}.}}
Using
\begin{align}\vs
\liminf_{m\to\infty}\frac{\frac{(x_1^{\rm Q})^2}{\sqrt{3}}}{D_L^{\rm Q}}=1,
\vs\end{align}
it suffices to show that $D_F^{\rm Q}\leq4\sqrt{\frac{2}{e}}\frac{(2-r_L)^2}{(1-r_L)^2}\phi(\frac{x_1^{\rm Q}}{\sqrt{2}})+r_L^2 D_L^{\rm Q}$.
Consider an interval $B_k^{\rm Q}$ and some $\theta\in B_k^{\rm Q}$. For any $j>k$ (similarly for $j<k$), we have
\begin{align}\vs
\frac{(c_j^{\rm Q}-a_L^{\rm Q}(\theta))^2}{(b_j^{\rm Q}-a_L^{\rm Q}(\theta))^2}\leq \frac{(2x_j^{\rm Q}-r_Lx_j^{\rm Q})^2}{(x_j^{\rm Q}-r_Lx_j^{\rm Q})^2}=\frac{(2-r_L)^2}{(1-r_L)^2}.
\vs\end{align}
Using this, we can obtain
\begin{align}\vs
&\int_{s\notin B_k^{\rm Q}}(a_F^{\rm Q}(s)- a_L^{\rm Q}(\theta))^2\phi(s-a_L(\theta))ds\nonumber\\
&\leq\frac{(2-r_L)^2}{(1-r_L)^2}\int_{s\notin B_k^{\rm Q}}(s- a_L^{\rm Q}(\theta))^2\phi(s-a_L^{\rm Q}(\theta))ds.
\vs\end{align}
Combining this with the inequality
$\int_{a}^{\infty}x^2\phi(x)dx\leq2\max(xe^{-\frac{x^2}{4}})\phi(\frac{a}{\sqrt{2}})$,
we arrive at
\begin{align}\vs
\int_{s\notin B_k^{\rm Q}}(a_F^{\rm Q}(s)- a_L^{\rm Q}(\theta))^2\phi(s-a_L^{\rm Q}(\theta))ds\leq4\sqrt{\frac{2}{e}}\phi(\frac{x_1^{\rm Q}}{\sqrt{2}}).
\vs\end{align}
On the other hand,
\begin{align}\vs
&\int_{s\in B_k^{\rm Q}}(a_F^{\rm Q}(s)- a_L^{\rm Q}(\theta))^2\phi(s-a_L(\theta))ds\nonumber\\
&\leq(c_k^{\rm Q}-a_L^{\rm Q}(\theta))^2=r_L^2(a_F^{\rm Q}(\theta)-\theta)^2,
\vs\end{align}
implying that
\begin{align}\vs
&\int_{-\infty}^{\infty}\int_{s\in B_k^{\rm Q}}(a_F^{\rm Q}(s)- a_L^{\rm Q}(\theta))^2\phi(s-a_L^{\rm Q}(\theta))\frac{\phi(\frac{\theta}{\sigma})}{\sigma}dsd\theta\leq
\nonumber\\
&r_L^2\int_{-\infty}^{\infty}(\theta-a_F^{\rm Q}(\theta))^2 \frac{\phi(\frac{\theta}{\sigma})}{\sigma}d\theta=r_L^2 D_L^{\rm Q},
\vs\end{align}
which completes the proof.~$\hfill\blacksquare$

{\label{proof::lemma::Sahai}\noindent\textit{Proof of Lemma~\ref{lemma::Sahai}.}}
First we note that the minimum value of the cost functional $U(a_L,a_F)$ with $r_L\sigma^2=1$ is asymptotically the same as the optimal cost of Witsenhausen's problem for $k^2\sigma^2=1$.  Using the inequality given by (17) in the proof of Theorem~4 in \cite{Grover_TAC_13}, we can obtain
\begin{align}\vs
U^*(\sigma)>\min_{P^*>0.5}\{k^2P^*+\frac{1}{15}e^{-12P^*}\},
\vs\end{align}
noting that in the scalar version of Witsenhausen's problem we have $m=1$. Minimizing the RHS above we can find
$U^*(\sigma)>\frac{\ln\sigma}{6\sigma^2}+\frac{1-\ln1.25}{12\sigma^2}$, which completes the proof.~$\hfill\blacksquare$

{\label{proof::theorem::performance}\noindent\textit{Proof of Theorem~\ref{theorem::performance}.}}
The first part of the theorem follows directly from Lemma~\ref{lemma::localoptima}. Using $M(\sigma)=\{m\in\mathbb{N}|x_1^{\rm Q}>2\sqrt{2\ln\sigma}+4, m\geq25\}$ and that $m\frac{x_1^{\rm Q}}{\sigma}\approx\frac{\sqrt{6\pi}}{4}$ for large $m$, we get
\begin{align}\vs
M(\sigma)\approx\{m\in\mathbb{N}|25<m<\frac{\sqrt{6\pi}\sigma}{8\sqrt{2\ln\sigma}+20}\}.
\vs\end{align}
Denote with $x_1^*$ the minimizer of the asymptotic upper bound on $U(a_L^*,a_F^*)$ given by Lemma~\ref{lemma::DLF}. Then, it is easy to verify that $\lim_{\sigma\to\infty}\frac{x_1^*}{2\sqrt{2\ln\sigma}}=1$, which clearly intersects $M(\sigma)$ for large $\sigma$. The corresponding asymptotic cost is $\approx\frac{8r_L\ln\sigma}{\sqrt{3}}$, hence proving \eqref{eq::th2}. We can use  Lemma~\ref{lemma::DLF} and Lemma~\ref{lemma::Sahai} to show that the equilibrium corresponding to $\argmin_{m\in M(\sigma)} U(a_L^*,a_F^*)$ is at most
$\frac{\frac{8r_L\ln\sigma}{\sqrt{3}}}{\frac{\ln\sigma}{6\sigma^2}}=16\sqrt{3}<27.8$ away from the optimal cost as $\sigma\to\infty$.~$\hfill\blacksquare$


\end{appendix}

\vspace{-14pt}

\bibliographystyle{IEEEtran}

\begin{thebibliography}{10}
\providecommand{\url}[1]{#1}
\csname url@samestyle\endcsname
\providecommand{\newblock}{\relax}
\providecommand{\bibinfo}[2]{#2}
\providecommand{\BIBentrySTDinterwordspacing}{\spaceskip=0pt\relax}
\providecommand{\BIBentryALTinterwordstretchfactor}{4}
\providecommand{\BIBentryALTinterwordspacing}{\spaceskip=\fontdimen2\font plus
\BIBentryALTinterwordstretchfactor\fontdimen3\font minus
  \fontdimen4\font\relax}
\providecommand{\BIBforeignlanguage}[2]{{%
\expandafter\ifx\csname l@#1\endcsname\relax
\typeout{** WARNING: IEEEtran.bst: No hyphenation pattern has been}%
\typeout{** loaded for the language `#1'. Using the pattern for}%
\typeout{** the default language instead.}%
\else
\language=\csname l@#1\endcsname
\fi
#2}}
\providecommand{\BIBdecl}{\relax}
\BIBdecl

\bibitem{Witsen_68}
H.~S. Witsenhausen, ``A counterexample in stochastic optimum control,''
  \emph{SIAM J. Control}, vol.~6, pp. 131--147, 1968.

\bibitem{Basar_TAC_76}
T.~Ba{\c{s}}ar, ``On the optimality of nonlinear designs in the control of
  linear decentralized systems,'' \emph{IEEE Transactions on Automatic
  Control}, vol.~21, p. 797, 1976.

\bibitem{Sobel_82}
V.~P. Crawford and J.~Sobel, ``Strategic information transmission,''
  \emph{Econometrica}, vol.~50, pp. 431--1451, 1982.

\bibitem{Serdar_TAC_17}
S.~Saritas, S.~Y{\"u}ksel, , and S.~Gezici, ``Quadratic multi-dimensional
  signaling games and affine equilibria,'' \emph{IEEE Transactions on Automatic
  Control}, vol.~62, pp. 605--619, 2017.

\bibitem{Ho_TAC_80}
Y.~Ho and T.~Chang, ``Another look at the nonclassical information structure
  problem,'' \emph{IEEE Transactions on Automatic Control}, vol.~25, pp.
  537--540, 1980.

\bibitem{Basar_TAC_87}
R.~Bansal and T.~Ba{\c{s}}ar, ``Stochastic teams with nonclassical information
  revisited: When is an affine law optimal?'' \emph{IEEE Transactions on
  Automatic Control}, vol.~32, pp. 554--559, 1987.

\bibitem{Sanjoy_99}
S.~Mitter and A.~Sahai, ``Information and control: {W}itsenhausen revisited,''
  \emph{Lecture Notes in Control and Information Sciences}, vol. 241, pp.
  281--293, 1999.

\bibitem{Verdu_CDC_11}
Y.~Wu and S.~Verd{\'u}, ``{W}itsenhausen's counterexample: A view from optimal
  transport theory,'' in \emph{Proceedings of the 50th IEEE Conf. Decision
  Control (CDC)}, 2011, pp. 5732--5737.

\bibitem{Grover_IT_15}
P.~Grover, A.~B. Wagner, and A.~Sahai, ``Information embedding and the triple
  role of control,'' \emph{IEEE Transactions on Information Theory}, vol.~61,
  pp. 1539--1549, 2015.

\bibitem{Basar_SIAM_15}
A.~Gupta, S.~Y{\"u}ksel, T.~Ba{\c{s}}ar, and C.~Langbort, ``On the existence of
  optimal policies for a class of static and sequential dynamic teams,''
  \emph{SIAM J. Control Optim.}, vol.~53, pp. 1681--1712, 2015.

\bibitem{Basar_Birkhauser_13}
S.~Y{\"u}ksel and T.~Ba{\c{s}}ar, ``Stochastic networked control systems:
  Stabilization and optimization under information constraints,''
  \emph{Birkh{\"a}user}, 2013.

\bibitem{Basar_CDC_08}
T.~Ba{\c{s}}ar, ``Variations on the theme of the {W}itsenhausen
  counterexample,'' in \emph{Proceedings of the 47th IEEE Conf. Decision
  Control (CDC)}, 2008, pp. 1614--1619.

\bibitem{Rotkowitz2006}
M.~Rotkowitz, ``Linear controllers are uniformly optimal for the {W}itsenhausen
  counterexample,'' in \emph{Proceedings of the 45th IEEE Conf. Decision
  Control (CDC)}, 2006, pp. 553--558.

\bibitem{Grover_IJSCC_2010}
P.~Grover and A.~Sahai, ``Vector {W}itsenhausen counterexample as assisted
  interference suppression,'' \emph{Int. J. Syst., Control Commun (IJSCC) Spec.
  iss. Inf Process. Decision Making in Distributed Control Syst.}, vol.~2, pp.
  197--237, 2010.

\bibitem{Grover_TAC_13}
P.~Grover, S.~Y. Park, and A.~Sahai, ``Approximately optimal solutions to the
  finite-dimensional {W}itsenhausen counterexample,'' \emph{IEEE Transactions
  on Automatic Control}, vol.~58, pp. 2189--2204, 2013.

\bibitem{Yuksel_arxiv_16}
N.~Saldi, S.~Y{\"u}ksel, and T.~Linder, ``Finite model approximations and
  asymptotic optimality of quantized policies in decentralized stochastic
  control,'' \emph{IEEE Transactions on Automatic Control}, vol.~62, pp.
  2360--2373, 2017.

\bibitem{Shamma_CDC}
N.~Li, J.~R. Marden, and J.~S. Shamma, ``Learning approaches to the
  {W}itsenhausen counterexample from a view of potential games,'' in
  \emph{Proceedings of the 48th IEEE Conf. Decision Control (CDC)}, 2009.

\bibitem{parisini_TAC_01}
M.~Baglietto, T.~Parisini, and R.~Zoppoli, ``Numerical solutions to the
  {W}itsenhausen counterexample by approximating networks,'' \emph{IEEE
  Transactions on Automatic Control}, vol.~46, pp. 1471--1477, 2001.

\bibitem{Lee_TAC_01}
J.~T. Lee, E.~Lau, and Y.-C. L.Ho, ``The {W}itsenhausen counterexample: A
  hierarchical search approach for nonconvex optimization problems,''
  \emph{IEEE Transactions on Automatic Control}, vol.~46, pp. 382--397, 2001.

\bibitem{Rose_ISIT_2014}
M.~Mehmetoglu, E.~Akyol, and K.~Rose, ``A deterministic annealing approach to
  {W}itsenhausen's counterexample,'' in \emph{ISIT}, 2014.

\bibitem{Vives_90}
X.~Vives, ``Complementarities and games: New developments,'' \emph{Journal of
  Economic Literature}, vol.~43, pp. 437--479, 2005.

\bibitem{Grover_WiOPT_09}
P.~Grover, A.~Sahai, and S.~Y. Park, ``Finite-dimensional {W}itsenhausen
  counterexample,'' in \emph{Proceedings of the 7th International Conference on
  Modeling and Optimization in Mobile, Ad Hoc, and Wireless Networks}, 2009,
  pp. 604--613.

\bibitem{Panter_IRE_51}
P.~F. Panter and W.~Dite, ``Quantization distortion in pulse count modulation
  with nonuniform spacing of levels,'' \emph{Proceedings of the I.R.E.},
  vol.~49, pp. 44--48, 1951.

\bibitem{Loyd_IT_82}
S.~Lloyd, ``Least squares quantization in {PCM},'' \emph{IEEE Transactions on
  Information Theory}, vol.~28, pp. 129--137, 1982.

\bibitem{Na_IT_01}
S.~Na and D.~L. Neuhoff, ``On the support of {MSE}-optimal, fixed-rate, scalar
  quantizers,'' \emph{IEEE Transactions on Information Theory}, vol.~47, pp.
  2972--2982, 2001.

\bibitem{Tirole_91}
D.~Fudenberg and J.~Tirole, \emph{Game Theory}.\hskip 1em plus 0.5em minus
  0.4em\relax MIT Press, 1991.

\bibitem{Normal_Horrace}
W.~C. Horrace, ``Moments of the truncated normal distribution,'' \emph{Journal
  of Productivity Analysis}, vol.~43, pp. 133--138, 2015.

\bibitem{Gasquet_Analysis}
C.~Gasquet and P.~Witomski, \emph{Fourier analysis and applications}.\hskip 1em
  plus 0.5em minus 0.4em\relax Springer, 1999.

\end{thebibliography}

\end{document}